\documentclass[11pt,a4paper]{article}

\usepackage{url,amsmath,amssymb,latexsym,pstricks,mathrsfs,comment,amsthm,graphicx,tikz,tikz-cd,enumerate,accents,pgffor,cite,wrapfig,multicol,float,cases,calc,bibspacing,geometry,bbm,arydshln}
\usepackage[colorlinks]{hyperref}
\usepackage[T1]{fontenc}
\usepackage{ wasysym }
\usepackage{stmaryrd}

\usetikzlibrary{calc}

\allowdisplaybreaks

\newcommand{\nc}{\newcommand}
\nc{\rnc}{\renewcommand}

\let\oldproofname=\proofname
\rnc{\proofname}{\rm\bf{\oldproofname}}

\rnc{\arraystretch}{1.2}

\nc{\C}{\mathscr C}
\rnc{\O}{\mathcal O}
\nc{\I}{\mathcal I}
\rnc{\S}{\mathcal G} 
\nc{\B}{\mathcal B} 
\rnc{\P}{\mathcal P} 
\nc{\PB}{\mathcal{PB}} 
\nc{\N}{\mathbb N}
\nc{\E}{\mathbb E}
\nc{\F}{\mathbb F}
\nc{\X}{\mathbb X}
\nc{\Y}{\mathbb Y}
\nc{\G}{\mathbb G}
\nc{\EX}{\mathcal E_X}
\nc{\FX}{\mathcal F_X}
\nc{\GX}{\S_X}
\nc{\FLX}{\FX^L}
\nc{\FRX}{\FX^R}
\nc{\GLX}{\GX^L}
\nc{\GRX}{\GX^R}
\nc{\Q}{\mathcal Q_X}

\nc{\om}{\omega}
\nc{\al}{\alpha}
\nc{\be}{\beta}
\nc{\ga}{\gamma}
\nc{\de}{\delta}
\nc{\ve}{\varepsilon}
\nc{\lam}{\lambda}
\nc{\si}{\sigma}
\nc{\ka}{\kappa}

\nc{\Ga}{\Gamma}
\nc{\Om}{\Omega}

\nc{\dom}{\operatorname{Dom}} 
\nc{\codom}{\operatorname{Codom}}
\nc{\rank}{\operatorname{rank}}
\nc{\relrank}[2]{\rank(#1\hspace{0.5truemm}{:}\hspace{0.5truemm}#2)}
\nc{\Fix}{\operatorname{Fix}}
\nc{\fix}{\operatorname{fix}}
\nc{\Supp}{\operatorname{Supp}}
\nc{\supp}{\operatorname{supp}}
\nc{\defect}{\operatorname{def}}
\nc{\Defect}{\operatorname{Def}}
\nc{\Codef}{\operatorname{Codef}}
\nc{\codef}{\operatorname{codef}}
\nc{\sh}{\operatorname{sh}}
\nc{\Sh}{\operatorname{Sh}}
\nc{\Fail}{\operatorname{Fail}}
\nc{\fail}{\operatorname{fail}}
\nc{\SR}{\operatorname{SR}}
\nc{\sing}{\operatorname{sing}}
\nc{\cosing}{\operatorname{cosing}}

\nc{\set}[2]{\{ {#1} : {#2} \}} 
\nc{\bigset}[2]{\big\{ {#1}: {#2} \big\}} 

\rnc{\emptyset}{\varnothing}
\nc{\1}{\mathbb I}
\nc{\lar}[1]{ \xrightarrow {\ #1\ }}
\nc{\mt}{\mapsto}
\nc{\sm}{\setminus}
\nc{\sub}{\subseteq}
\nc{\suq}{\subsetneq}
\nc{\la}{\langle}
\nc{\ra}{\rangle}

\nc{\OR}{\qquad\text{or}\qquad}
\nc{\COMMA}{,\quad}
\nc{\COMMa}{,\ \ }
\nc{\AND}{\qquad\text{and}\qquad}
\nc{\ANDSIM}{\qquad\text{and similarly}\qquad}
\rnc{\iff}{\ \Leftrightarrow\ }
\rnc{\implies}{\ \Rightarrow\ }

\nc{\bit}{\begin{itemize}}
\nc{\eit}{\end{itemize}}
\nc{\bmc}{\begin{multicols}}
\nc{\emc}{\end{multicols}}
\nc{\itemit}[1]{\item[\emph{(#1)}]}
\nc{\itemnit}[1]{\item[(#1)]}
\nc{\pfitem}[1]{\medskip \noindent (#1).}
\nc{\pfcase}[1]{\medskip \noindent {\bf Case #1.}}
\nc{\pf}{\begin{proof}}
\nc{\epf}{\end{proof}}
\nc{\epfres}{\hfill$\Box$}
\nc{\epfreseq}{\tag*{$\Box$}}

\rnc{\c}{@{}c@{}}
\nc{\cend}{@{}c@{\hspace{1.5truemm}}}
\nc{\sqcend}{@{}c@{\hspace{1.0truemm}}}
\nc{\cstart}{@{\hspace{1.5truemm}}c@{}}

\nc{\partn}[4]{
\Big(  
{\scriptsize \rnc*{\arraystretch}{1} \begin{array} {\c|\cend}
#1 \:&\: #2  \\ \cline{2-2}
#3 \:&\: #4
\rule[0mm]{0mm}{2.7mm}
\end{array}  }
\hspace{-1.5 truemm} \Big) 
}

\nc{\sqpartn}[4]{
\Big[
{\scriptsize \rnc*{\arraystretch}{1} \begin{array} {\c|\sqcend}
#1 \:&\: #2  \\ \cline{2-2}
#3 \:&\: #4
\rule[0mm]{0mm}{2.7mm}
\end{array}  }
\hspace{-0.8 truemm} \Big]
}

\nc{\partnII}[8]{
\Big(  \hspace{-1.5 truemm}
{\scriptsize \rnc*{\arraystretch}{1} \begin{array} {\cstart|\c|\c|\cend}
#1 \:&\: #2 \:&\: #3 \:&\: #4 \\ 
#5 \:&\: #6 \:&\: #7 \:&\: #8
\rule[0mm]{0mm}{2.7mm}
\end{array}  }
\hspace{-1.5 truemm} \Big) 
}

\nc{\partnIII}[6]{
\Big(  \hspace{-1.5 truemm}
{\scriptsize \rnc*{\arraystretch}{1} \begin{array} {\cstart|\c|\cend}
#1 \:&\: #2 \:&\: #3 \\ 
#4 \:&\: #5 \:&\: #6
\rule[0mm]{0mm}{2.7mm}
\end{array}  }
\hspace{-1.5 truemm} \Big) 
}

\nc{\partnIV}[6]{
\Big(  \hspace{-1.5 truemm}
{\scriptsize \rnc*{\arraystretch}{1} \begin{array} {\cstart|\c|\c|\c|\cend}
\cdots \:&\: #1 \:&\: #2 \:&\: #3 \:&\: \cdots \\ 
\cdots \:&\: #4 \:&\: #5 \:&\: #6 \:&\: \cdots
\rule[0mm]{0mm}{2.7mm}
\end{array}  }
\hspace{-1.5 truemm} \Big) 
}

\nc{\partnV}[4]{
\Big(  
{\scriptsize \rnc*{\arraystretch}{1} \begin{array} {\c|\cend}
#1 \:&\: #2  \\ 
#3 \:&\: #4
\rule[0mm]{0mm}{2.7mm}
\end{array}  }
\hspace{-1.5 truemm} \Big) 
}

\nc{\partnVI}[2]{
\Big(
{ \scriptsize \rnc*{\arraystretch}{1}
\begin{array} {\c}
 #1  \\ 
 #2
\rule[0mm]{0mm}{2.7mm}
\end{array} 
}
\Big)
}

\nc{\partnVII}[8]{
\Big(  \hspace{-1.5 truemm}
{\scriptsize \rnc*{\arraystretch}{1} \begin{array} {\cstart|\c|\c|\cend}
#1 \:&\: #2 \:&\: #3 \:&\: #4 \\ \cline{2-4}
#5 \:&\: #6 \:&\: #7 \:&\: #8
\rule[0mm]{0mm}{2.7mm}
\end{array}  }
\hspace{-1.5 truemm} \Big) 
}

\nc{\colv}[3]{\fill[#3] (#1,#2)circle(.17);}
\nc{\uv}[1]{\fill (#1,2)circle(.13);}
\nc{\lv}[1]{\fill (#1,0)circle(.13);}
\nc{\uvg}[1]{\fill[lightgray] (#1,2)circle(.17);}
\nc{\lvg}[1]{\fill[lightgray] (#1,0)circle(.17);}
\nc{\uvs}[1]{{\foreach \x in {#1} { \uv{\x}}}}
\nc{\lvs}[1]{{\foreach \x in {#1} { \lv{\x}}}}

\nc{\buv}[1]{\fill (#1,2)circle(.18);}
\nc{\buvs}[1]{{
\foreach \x in {#1}
{ \buv{\x}}
}}
\nc{\blv}[1]{\fill (#1,0)circle(.18);}
\nc{\blvs}[1]{{
\foreach \x in {#1}
{ \blv{\x}}
}}

\nc{\uarcs}[1]{
{\foreach \x/\y in {#1}
{ \uarc{\x}{\y} }
}
}

\nc{\darcs}[1]{
{\foreach \x/\y in {#1}
{ \darc{\x}{\y} }
}
}

\nc{\darcxhalf}[3]{\draw(#1,0)arc(180:90:#3) (#1+#3,#3)--(#2,#3) ;}
\nc{\darchalf}[2]{\darcxhalf{#1}{#2}{.4}}
\nc{\uarcxhalf}[3]{\draw(#1,2)arc(180:270:#3) (#1+#3,1.5-#3)--(#2,1.5-#3) ;}
\nc{\uarchalf}[2]{\uarcxhalf{#1}{#2}{.4}}

\nc{\stline}[2]{\draw(#1,2)--(#2,0);}
\nc{\stlines}[1]{{\foreach \x/\y in {#1} { \stline{\x}{\y} }}}

\nc{\darcx}[3]{\draw(#1,0)arc(180:90:#3) (#1+#3,#3)--(#2-#3,#3) (#2-#3,#3) arc(90:0:#3);}
\nc{\darc}[2]{\darcx{#1}{#2}{.4}}
\nc{\uarcx}[3]{\draw(#1,2)arc(180:270:#3) (#1+#3,2-#3)--(#2-#3,2-#3) (#2-#3,2-#3) arc(270:360:#3);}
\nc{\uarc}[2]{\uarcx{#1}{#2}{.4}}

\numberwithin{equation}{section}

\newtheorem{thm}[equation]{Theorem}
\newtheorem{lemma}[equation]{Lemma}
\newtheorem{cor}[equation]{Corollary}

\theoremstyle{definition}

\newtheorem{rem}[equation]{Remark}

\begin{document}

\title{Idempotents and one-sided units in infinite partial Brauer monoids\vspace{-0.3cm}}
\author{James East\\
{\it\small Centre for Research in Mathematics; School of Computing, Engineering and Mathematics,}\\
{\it\small Western Sydney University, Locked Bag 1797, Penrith NSW 2751, Australia.}\\
{\tt\small J.East\,@\,WesternSydney.edu.au}}
\date{}

\maketitle

\vspace{-1.0cm}

\begin{abstract}
We study monoids generated by various combinations of idempotents and one- or two-sided units of an infinite partial Brauer monoid.  This yields a total of eight such monoids, each with a natural characterisation in terms of relationships between parameters associated to Brauer graphs.  We calculate the relative ranks of each monoid modulo any other such monoid it may contain, and then apply these results to determine the Sierpi\'nski rank of each monoid, and ascertain which ones have the semigroup Bergman property.  We also make some fundamental observations about idempotents and units in arbitrary monoids, and prove some general results about relative ranks for submonoids generated by these sets.



\end{abstract}

\thanks{Dedicated to Dr Des FitzGerald on the occasion of his 70th birthday.}

\tableofcontents

\section{Introduction}\label{sect:intro}

An idempotent in an algebraic structure with a product is an element $x$ satisfying $x=x^2$.  Idempotents have long played an important role in semigroup theory and other branches of mathematics, and there exist many interesting results.  For example, Erdos showed in 1967 that any singular square matrix over a field is a product of idempotent matrices \cite{Erdos1967}; this followed in the footsteps of an earlier result of Howie \cite{Howie1966}, which showed that any non-bijective mapping of a finite set to itself is a product of idempotent mappings.  In the same paper, Howie also characterised the products of idempotent mappings on an infinite set; a crucial role was played by certain parameters that quantify how far a mapping is from being injective or surjective.  

The above-mentioned papers have generated a substantial literature that is still growing today, with many subsequent studies uncovering intriguing connections to finite combinatorics or infinite cardinal arithmetic.  To list a select few examples: Fountain and Lewin simulaneously extended the Erdos and Howie results above to endomorphism monoids of independence algebras \cite{FL1992,FL1993}; Gray showed (among many other things) that every singular $n\times n$ matrix of rank at most $r$ over a field is a product of idempotent matrices of rank~$r$, and calculated the minimal number of (idempotent) matrices required to generate all such matrices \cite{Gray2007}; Howie and his collaborators conducted further studies on mappings of finite sets \cite{Howie1978,HLM1990,HM1990,GH1987}; more recently, others have considered idempotent-generation in finite and infinite diagram monoids \cite{EF2012,East2011,MM2007,EG2017,DEG2017,East2014}.  For more background on the role of idempotents in semigroup theory, including applications to many branches of mathematics not mentioned here, we refer to the introductions of \cite{DGR2017,DGY2015,EG2017} for thorough discussions.

The above-mentioned article of Fountain and Lewin \cite{FL1993} also considered products of idempotents and units (a unit of a monoid is an element $x$ with a two-sided inverse: $ax=xa=1$ for some $a$).  In fact, in order to describe the submonoid of the endomorphism monoid of an infinite dimensional independence algebra, the submonoid generated by idempotents and units was first described.  Monoids generated by idempotents and units have also been studied in many other contexts; see for example \cite{East2014,EF2012,HHR1998,CH1974,East2012,FitzGerald2010}.  Of particular immediate relevance is the article of Higgins, Howie and Ru\v skuc \cite{HHR1998}, in which \emph{one-sided} units in the monoid~$P$ of all partial mappings of an infinite set to itself were also considered (a one-sided unit of a monoid is an element~$x$ with a one-sided inverse: $ax=1$ for some $a$, or~$xb=1$ for some $b$, or possibly both).  Denoting by $S$, $I$, $G$ and $E$ the sets of all surjective, injective, bijective and idempotent mappings, respectively, they considered all products of these sets: for example, it was shown that the set~$IS=\set{fg}{f\in I,\ g\in S}$ is equal to all of $P$.  All other products of two or more of these sets were calculated, and the semigroup~$M=\la S,I,G,E\ra$ generated by all four sets was described.  It is important to note here that~$M$ is not a subsemigroup of $P$ itself, but rather of the \emph{power semigroup} of~$P$; the latter consists of all subsets of $P$, with the semigroup operation being set product.  Subsemigroups of $P$ generated by unions of the above sets were not explicitly considered in \cite{HHR1998}, but descriptions of them may be deduced from results therein: for example,~$\la S\cup I\ra=P$ and $\la E\cup G\ra=EG=GE$ consists of all so-called \emph{semi-balanced} mappings.  It was also shown that two (but no fewer) elements of $P$ may be added to $E\cup G$ in order to obtain a generating set for $P$.  This last result can be stated in terms of relative ranks: the \emph{relative rank} \cite{HRH1998} of a semigroup $T$ modulo a subset $A\sub T$, denoted $\relrank TA$, is the minimum size of a subset $U\sub T$ such that $T=\la A\cup U\ra$; thus, the aforementioned result from \cite{HHR1998} states that $\relrank{P}{E\cup G}=2$.   This extends other results of the same authors \cite{HRH1998}, which calculate relative ranks in monoids of (full) mappings modulo the sets of idempotents or units.  It follows from the proof of \cite[Lemma 4.2]{HHR1998} that the sets $I$ and $S$ are precisely the right and left units of $P$, respectively.  A few results from \cite{HHR1998} concerning $P$ were established by proving general results about arbitrary monoids; several others may also be deduced from further general results we prove in Section \ref{sect:monoids} below.  We also note that Mitchell and P\'eresse \cite{MP2011} have (among other things) calculated the relative ranks of the monoids of all (full) injective or surjective mappings on an infinite set modulo the bijective mappings; thus, this is an instance of calculating relative ranks of the left (or right) units of a monoid (the monoid of \emph{all} mappings in this case) modulo the two-sided units.  

In \cite{EF2012}, the idempotent-generated submonoid of an infinite partition monoid was described, as well as the submonoid generated by the idempotents and units (see \cite{East2011} for the finite case).  Partition monoids, and other \emph{diagram monoids} such as Brauer and Temperley-Lieb monoids, arise in many branches of mathematics, including knot theory, theoretical physics and representation theory \cite{HR2005,Brauer1937,Jones1987,Kauffman1990,Jones1994_2,Martin1994, TL1971}; see also the introductions of \cite{DEEFHHL1,EG2017} for a discussion of the fruitful relationship between diagram monoids and semigroup theory.  One-sided units did not feature in \cite{EF2012}, but they were used implicitly in \cite{East2014}, where it was shown that every element of an infinite partition monoid is a product of a right unit by a left unit (in that order, but not the other).  Other results of \cite{East2014} included the calculation of the relative ranks of an infinite partition monoid modulo its (two-sided) units and/or idempotents.  Applications of these results included proofs that infinite partition monoids have the \emph{(semigroup) Bergman property}, and also finite \emph{Sierpi\'nski rank}.  A semigroup $S$ has the Bergman property \cite{MMR2009,Bergman2006} if every generating set for $S$ has a bounded length function, while $S$ has finite Sierpi\'nski index \cite{MP2012,Sierpinski1935,Banach1935} if there exists a natural number $n$ such that every countable subset of $S$ is contained in a subsemigroup generated by $n$ elements, in which case the least such $n$ is the Sierpi\'nski index.

The current article furthers the above body of work in several directions.  Our main motivating examples are the infinite \emph{partial Brauer monoids} $\PB_X$; these will be defined in Section \ref{sect:PBX}, where we also explain why there are no infinite \emph{full} Brauer monoids.  As well as extending the results of \cite{EF2012,East2014} to $\PB_X$, we introduce new techniques for working with submonoids generated not just by idempotents and two-sided units, but also by idempotents and one-sided units; the latter tend to have much more complicated structures (for one thing, they are not \emph{regular} if there are one-sided units that are not two-sided; see Remark \ref{rem:nonregular} below).  We also develop a general theory of idempotents and one-sided units in arbitrary monoids; we hope this will be useful in subsequent studies.  This general theory is expounded in Section \ref{sect:monoids}, which also gives definitions and background on semigroups and monoids in general.  The partial Brauer monoids $\PB_X$ are introduced in Section \ref{sect:PBX}, as well as a number of parameters (sets and cardinals) associated to the elements of $\PB_X$, and we prove a number of inequalities related to these.  Sections \ref{sect:SX}--\ref{sect:FLFR} study the submonoids of $\PB_X$ generated by all combinations of one- or two-sided units and/or idempotents; in these sections, we characterise the elements of each monoid, calculate the relative ranks of each one modulo any other such monoid it may contain, and classify the minimal-size generating sets modulo any such submonoid.  Section \ref{sect:SB} calculates the Sierpi\'nski rank of each monoid, and determines which of them have the semigroup Bergman property; a centrepiece of this section is a proof (modelled on an ingenious argument of Hyde and P\'eresse \cite{HP2012}) that the Sierpi\'nski rank of $\PB_X$ is equal to $2$.  
The main results, and their locations, are summarised in Table \ref{tab:summary}, which uses the shorthand notation for the various submonoids of $\PB_X$ we consider:~$\EX$ denotes the idempotent-generated submonoid; $\S_X$ is the group of units; $\GLX$ (respectively, $\GRX$) is the monoid of all left (respectively, right) units; $\FX$ is the monoid generated by all idempotents and two-sided units; and $\FLX$ (respectively, $\FRX$) is the monoid generated by all idempotents and left (respectively, right) units.

\begin{table}[ht]
\begin{center}
\begin{tabular}{|l|l|}
\hline
Lemma \ref{lem:PBunits} & Description of $\GLX$ and $\S_X$ \\
Theorem \ref{thm:IGPBX} & Description of $\EX$ \\
Theorem \ref{thm:FPBX} & Description of $\FX$ \\
Theorem \ref{thm:FLFR} & Description of $\FLX$ \\
\hline
Theorem \ref{thm:PBXSX} & $\relrank{\PB_X}{\S_X}=2$ \\
Theorem \ref{thm:PBXGLX} & $\relrank{\PB_X}{\GLX}=1$ \\
Theorem \ref{thm:PBXEX} & $\relrank{\PB_X}{\EX}=2$ \\
Theorem \ref{thm:PBXFX} & $\relrank{\PB_X}{\FX}=2$ \\
Theorem \ref{thm:PBXFLX} & $\relrank{\PB_X}{\FLX}=1$ \\
\hline
\end{tabular}
\quad
\begin{tabular}{|l|l|}
\hline
Theorem \ref{thm:FLXFX} & $\relrank{\FLX}{\FX}=1+\rho$ \\
Theorem \ref{thm:FLXEX} & $\relrank{\FLX}{\EX}=2^{|X|}$ \\
Theorem \ref{thm:FLXGLX} & $\relrank{\FLX}{\GLX}=2+2\rho$ \\
Theorem \ref{thm:FLXSX} & $\relrank{\FLX}{\S_X}=3+3\rho$ \\
\hline
Theorem \ref{thm:FXEX} & $\relrank{\FX}{\EX}=2^{|X|}$ \\
Theorem \ref{thm:FXSX} & $\relrank{\FX}{\S_X}=2+2\rho$ \\
\hline
Theorem \ref{thm:GLXSX} & $\relrank{\GLX}{\S_X}=2+2\rho$ \\
\hline
Theorem \ref{thm:SBPB} & Bergman/Sierpi\'nski in $\PB_X$ \\
Theorem \ref{thm:SBEF} & Bergman/Sierpi\'nski in all other monoids \\
\hline
\end{tabular}
\end{center}
\vspace{-.5cm}
\label{tab:summary}
\caption{Summary and location of the main results.  Any result concerning $\GLX$ or $\FLX$ leads to dual results concerning $\GRX$ or~$\FRX$.  Here, $X$ is an infinite set, and $\rho$ denotes the number of infinite cardinals not exceeding~$|X|$.}
\end{table}

Throughout, we denote the set of natural numbers by $\N=\{0,1,2,\ldots\}$.  We use the $\sqcup$ symbol to denote disjoint union.  When we list the elements of a set as $\{x_1,x_2,\ldots\}$ or $\set{y_i}{i\in I}$, etc., we always assume that different subscripts give rise to different elements of the set.  Functions are generally written to the right of their arguments, and are composed from left to right.  If $a_1\cdots a_k$ denotes a product of elements from some monoid, then this represents the identity element if $k=0$; similar conventions hold for empty sums and lists.  We assume basic results concerning infinite cardinals, such as may be found in \cite[Chapter 5]{Jech2003}, for example.

\section{Monoids}\label{sect:monoids}

In this section, we provide some background on semigroups and monoids, and prove a number of results concerning idempotents and units in arbitrary monoids.  Some of these results are structural (Lemmas~\mbox{\ref{lem:GLGRG}--\ref{lem:FFLFR}}), while some give information concerning relative ranks of various submonoids inside others (Lemmas~\ref{lem:rankMGE} and~\ref{lem:rankME}).  

A semigroup is a set $S$ with an associative binary operation.  If $U$ is a subset of $S$, we write $\la U\ra$ for the subsemigroup of $S$ generated by $U$; so $\la U\ra$ is the smallest subsemigroup of $S$ containing $U$, and consists of all products $u_1\cdots u_k$, where $k\geq1$ and $u_1,\ldots,u_k\in U$.  Following \cite{HM1990}, the \emph{rank} of $S$ is defined by
\[
\rank(S) = \min\bigset{|U|}{U\sub S,\ S=\la U\ra}.
\]
The semigroups we are primarily interested in are all uncountable; for any such semigroup, it is easy to see that $\rank(S)=|S|$.  Thus, a more useful concept for uncountable semigroups is that of \emph{relative} rank.  Following \cite{HRH1998}, if $A\sub S$, the \emph{relative rank of $S$ modulo~$A$} is defined by
\[
\relrank SA = \min\bigset{|U|}{U\sub S,\ S=\la A\cup U\ra}.
\]
It is possible for $S\sm A$ to be uncountable, yet for $\relrank SA$ to be finite; indeed, we provide several examples in the current paper, and many more exist in the literature; see for example \cite{HRH1998,HHMR2003,East2014,MP2011,HP2012,HHR1998,East2012}.

A monoid is a semigroup $M$ with an identity element $1$.  A submonoid of $M$ is a subsemigroup of $M$ that contains $1$.  Following \cite[Section 1.7]{CPbook}, an element $x$ of $M$ is a \emph{left unit} if $ax=1$ for some $a\in M$, in which case we say that $a$ is a \emph{left inverse} of $x$.  Right units and right inverses are defined analogously.  
A \emph{(two-sided) unit} of~$M$ is an element $x$ that is both a left and right unit.  It is a routine exercise to show that a unit $x$ has a unique left inverse and a unique right inverse, and that these are equal, in which case we write $x^{-1}$ for the unique two-sided inverse of $x$.  We denote by $\G_L(M)$ and $\G_R(M)$ the sets of all left and right units of~$M$, respectively, and by $\G(M)=\G_L(M)\cap \G_R(M)$ the set of all units.
\emph{Green's relations} (see \cite[Chapter~2]{Howie}) will not play an explicit role in this paper, but we note that $\G_L(M)$, $\G_R(M)$ and $\G(M)$ are the $\mathscr L$-, $\mathscr R$- and $\mathscr H$-classes of $1$ in $M$, respectively.
%

If $U$ is a subset of a semigroup $S$, we write $E(U)=\set{u\in U}{u=u^2}$ for the set of all idempotents of $U$.  We write
\[
\E(S)=\la E(S)\ra
\]
for the subsemigroup of $S$ generated by all of its idempotents.  A \emph{left ideal} of a semigroup $S$ is a subset $I$ of~$S$ such that $sx\in I$ for all $x\in I$ and $s\in S$.  Right ideals are definied analogously.  An \emph{ideal} is a non-empty subset that is both a left and right ideal.  The proof of the next result is routine, and is omitted; for part (i), see \cite[Theorem 1.10]{CPbook}.

\begin{lemma}\label{lem:GLGRG}
Let $M$ be a monoid, and write $G_L=\G_L(M)$, $G_R=\G_R(M)$ and $G=\G(M)$.  Then
\bit
\itemit{i} $G_L$, $G_R$ and $G$ are all submonoids of $M$, with $G$ a group,
\itemit{ii} $M\sm G_L$ is a left ideal of $M$, and $M\sm G_R$ is a right ideal,
\itemit{iii} $E(G_L)=E(G_R)=E(G)=\{1\}$,
\itemit{iv} $G_L\cap\E(M)=G_R\cap\E(M)=G\cap\E(M)=\{1\}$.  \epfres
\eit
\end{lemma}

%
%
%
%

\begin{rem}
It follows from Lemma \ref{lem:GLGRG}(ii) that $(M\sm G_L)\cap(M\sm G_R)=M\sm(G_L\cup G_R)$ is a subsemigroup of $M$, though it need not be an ideal.
\end{rem}

Recall that a monoid $M$ is \emph{bicyclic} if it is generated by two elements $a,b$ satisfying $ab=1\not=ba$.  All bicyclic monoids are isomorphic to each other, and can be defined by the presentation $\la a,b:ab=1\ra$.  See \cite[pp.~31--32]{Howie} for more details.  Again, the proof of the next result is routine, and is omitted; see Exercise~1(a) of \cite[Section 1.7]{CPbook} and also \cite[Theorem 2.54]{CPbook}.

\begin{lemma}\label{lem:equiv}
Let $M$ be a monoid, and write $G_L=\G_L(M)$, $G_R=\G_R(M)$ and $G=\G(M)$.  Then the following are equivalent:
\bit\bmc2
\itemit{i} $G_L=G$ 
\itemit{ii} $G_R=G$,
\itemit{iii} $M\sm G$ is an ideal of $M$,
\itemit{iv} $M$ has no bicyclic submonoid.  \epfres
\emc\eit
\end{lemma}

%
%
%
%

\begin{rem}
Since bicyclic monoids are infinite, the previous result implies that $\G_L(M)=\G_R(M)=\G(M)$ if $M$ is finite.  
\end{rem}

We have so far considered submonoids consisting of one- and/or two-sided units only.  We now include idempotents.  If $M$ is a monoid, we define
\[
\F(M)=\la E(M)\cup \G(M)\ra \COMMA \F_L(M)=\la E(M)\cup \G_L(M)\ra \COMMA \F_R(M)=\la E(M)\cup \G_R(M)\ra
\]
for the submonoids of $M$ generated by all idempotents and two-sided units, or all idempotents and left units, or all idempotents and right units, respectively.

\begin{lemma}\label{lem:factorisation}
If $M$ is a monoid, then
\bit\bmc2
\itemit{i} $\F_L(M)=\E(M)\G_L(M)$,
\itemit{ii} $\F_R(M)=\G_R(M)\E(M)$,
\itemit{iii} $\F(M)=\E(M)\G(M)=\G(M)\E(M)$.
\emc\eit
\end{lemma}

\pf
Part (iii) is \cite[Lemma 32]{EF2012}.  By duality, it remains to prove (i).  During the proof, we use the abbreviations $E=\E(M)$, $G_L=\G_L(M)$ and $F_L=\F_L(M)$.  Clearly $EG_L\sub\la E\cup G_L\ra=F_L$.  
We can prove the reverse containment by showing that $EG_L$ is a subsemigroup of $M$ containing $E\cup G_L$, since~$F_L$ is the \emph{smallest} such subsemigroup.  As $E\cup G_L\sub EG_L$ is clear, suppose $x,y\in EG_L$, so that $x=eg$ and $y=fh$ for some $e,f\in E$ and $g,h\in G_L$.  Then $1=ag$ for some $a\in M$, and $f=f_1f_2\cdots f_k$ for some $f_1,f_2,\ldots,f_k\in E(M)$.  Then
\[
xy=egf_1f_2\cdots f_kh = egf_1(ag)f_2(ag)\cdots f_k(ag)h = e(gf_1a)(gf_2a)\cdots(gf_ka)gh.
\]
Since $gh\in G_L$ by Lemma \ref{lem:GLGRG}(i), and since $gf_ia\in E(M)$ for each $i$, it follows that $xy\in EG_L$.
\epf

\begin{rem}
The factorisations in Lemma \ref{lem:factorisation} are the reason for the use of the $\F$ symbol.  If $E(M)$ is a submonoid of~$M$ (if~$M$ is inverse, for example), or even if $E(M)^2\sub E(M)\G(M)$, then $\F(M)=E(M)\G(M)$, $\F_L(M)=E(M)\G_L(M)$, and so on; although these simplified factorisations do not hold for arbitrary monoids, we will see in Theorems~\ref{thm:FPBX} and \ref{thm:FLFR} that they do hold when $M$ is a partial Brauer monoid $\PB_X$ (defined in Section \ref{sect:PBX}), even though $E(\PB_X)$ is not a submonoid.
\end{rem}

\begin{rem}\label{rem:factorisation}
Note that we also have $\F_L(M)=\E(M)\G_L(M) = \E(M) [\G(M)\G_L(M)] =\F(M)\G_L(M)$, and similarly $\F_R(M)=\G_R(M)\F(M)$.
\end{rem}

The next two lemmas give some information on what happens when we iterate the above constructions, and consider submonoids of $M$ such as $\G_L(\F_R(M))$.  These will be important when we study $\G_L(M)$, $\F_L(M)$, etc., as monoids in their own right.

\begin{lemma}\label{lem:FFLFR}
If $M$ is a monoid, and if $Q$ is any of $\F_L(M)$, $\F_R(M)$ or $\F(M)$, then
\[
E(Q)=E(M) \COMMA \G_L(Q)=\G_R(Q)=\G(Q)=\G(M) \COMMA \F_L(Q)=\F_R(Q)=\F(Q)=\F(M).
\]
\end{lemma}

\pf
We just prove the statements for $Q=\F_L(M)$, as the others are similar.  During the proof, we also write $E=\E(M)$, $G_L=\G_L(M)$, $G_R=\G_R(M)$, $G=\G(M)$, $F_L=\F_L(M)$, $F_R=\F_R(M)$ and $F=\F(M)$.  

First, $F_L\sub M$ gives $E(F_L)\sub E(M)$ and $\G(F_L)\sub\G(M)=G$.  The reverse containments hold because $E(M)\cup G\sub\la E(M)\cup G\ra\sub\la E(M)\cup G_L\ra=F_L$.  

Next, suppose $x\in \G_L(F_L)$.  Then $1=ax$ for some $a\in F_L$; note that $a\in G_R$.  By Lemma \ref{lem:factorisation}(i), $F_L=EG_L$, so we may write $a=e_1\cdots e_kg$, where $k\geq0$, $e_1,\ldots,e_k\in E(M)$ and $g\in G_L$; we assume that~$k$ is minimal among all such expressions.  If $k\geq1$, then $a=e_1a$, and so $1=ax=e_1ax=e_1$, which gives $a=e_2\cdots e_kg$, contradicting the minimality of $k$.  It follows that $k=0$, and so $a=g\in G_L$.  But then $a\in G_L\cap G_R=G$, and so $1=ax$ gives  $x=a^{-1}\in G$.  This shows that $\G_L(F_L)\sub G=\G(F_L)$.  The reverse containment is obvious, and so $\G_L(F_L)=G$.  By Lemma \ref{lem:equiv}, it also follows that $\G_R(F_L)=G$.

The other statements follow quickly: for example, $\F_R(Q)=\la E(Q)\cup\G_R(Q)\ra=\la E(M)\cup G\ra=F$.
\epf

The proof of the next result is similar; see also Exercise 1(c) of \cite[Section 1.7]{CPbook}.

\begin{lemma}\label{lem:GGLGR}
If $M$ is a monoid, and if $Q$ is any of $\G_L(M)$, $\G_R(M)$ or $\G(M)$, then
\begin{equation}\tag*{\qed}
\G_L(Q)=\G_R(Q)=\G(Q)=\F_L(Q)=\F_R(Q)=\F(Q)=\G(M).
\end{equation}
\end{lemma}


\begin{rem}\label{rem:nonregular}
Recall that a semigroup $S$ is \emph{(von Neumann) regular} if, for each $x\in S$, there exists $u\in S$ such that $x=xux$.  It follows from Lemmas \ref{lem:FFLFR} and \ref{lem:GGLGR} that (using the usual abbreviations) if $G_L\not=G$, then~$G_L$ and $F_L$ are not regular, even if $M$ is itself regular.  Indeed, suppose $x\in G_L$, is such that $x$ is regular in $F_L$.  Then there exist $a\in M$ and $u\in F_L$ such that $1=ax$ and $x=xux$.  But then $1=ax=axux=ux$, so that $x\in \G_L(F_L)=G$.  This shows that no element of $G_L\sm G$ is regular in $F_L$ (even though all elements of $G_L$ are regular in $M$).  By contrast, if $M$ is regular, then $F$ must be regular; this follows quickly from a famous result of FitzGerald \cite{FitzGerald1972}, which says that if $M$ is regular, then so too is $\E(M)$.
\end{rem}

The next two lemmas give some information on relative ranks for certain pairs of (sub)monoids considered so far.
Clearly $\relrank SA=\relrank S{\la A\ra}$ for any subset $A$ of a semigroup $S$.  Thus, for example, $\relrank S{E(S)}=\relrank S{\E(S)}$ for any semigroup $S$, and $\relrank M{E(M)\cup \G(M)}=\relrank M{\F(M)}$ for any monoid $M$.


\begin{lemma}\label{lem:rankMGE}
Let $M$ be a monoid, and write $G_L=\G_L(M)$, $G_R=\G_R(M)$, $G=\G(M)$, $E=\E(M)$ and $F=\F(M)$.  Suppose also that $G_L\not=G$ (or, equivalently by Lemma \ref{lem:equiv}, that $G_R\not=G$).  Then
\bit
\itemit{i} if $M=\la U\ra$, then $U\sm F$ contains at least one element from $G_L\sm G$, and at least one from $G_R\sm G$, 
\bmc3
\itemit{ii} $\relrank MG\geq2$,
\itemit{iii} $\relrank ME\geq2$,
\itemit{iv} $\relrank MF\geq2$.
\emc
\eit
\end{lemma}

\pf (i).  Suppose $M=\la U\ra$.  By duality, it suffices to show that $U\sm F$ contains an element of $G_L\sm G$.  Since $F\sub M=\la U\ra$, certainly
\[
M=\la(U\sm F)\cup F\ra =  \la(U\sm F)\cup\la E\cup G\ra \ra =  \la(U\sm F)\cup E\cup G \ra.
\]
Now let $a\in G_L\sm G$ be arbitrary, and consider an expression $a=u_1\cdots u_k$, where all of the factors belong to $(U\sm F)\cup E\cup G$.  Since $a\not\in G$, the~$u_i$ cannot all belong to $G$.  Let $j=\max\set{i}{u_i\not\in G}$, and put $b=u_{j+1}\cdots u_k\in G$.  Then $ab^{-1}=u_1\cdots u_j$.  By Lemma~\ref{lem:GLGRG}(i), $ab^{-1}\in G_L$.  By Lemma \ref{lem:GLGRG}(ii), $u_j\in G_L$ (or otherwise $ab^{-1}=(u_1\cdots u_{j-1})u_j\in M\sm G_L$, a contradiction).  Since $u_j\not\in G$, it follows that $u_j\in G_L\sm G$.  By Lemma~\ref{lem:GLGRG}(iv), $(G_L\sm G)\cap E=\emptyset$, so it follows that $u_j\in U\sm F$.

\pfitem{iv}  Suppose $M=\la F\cup V\ra$, where $|V|=\relrank MF$.  By part (i), there exist $x,y\in V$ such that $x\in G_L\sm G$ and $y\in G_R\sm G$.  Since $(G_L\sm G)\cap(G_R\sm G)=\emptyset$, it follows that $x\not=y$, and so $\relrank MF=|V|\geq2$.

\pfitem{ii) and (iii}  These follow immediately from (iv), and the fact that $\relrank SA\geq\relrank SB$ for any semigroup $S$ with nested subsets $A\sub B\sub S$.
\epf

\begin{rem}
Lemma \ref{lem:rankMGE} applies to several well-studied monoids, including infinite full and partial transformation monoids, monoids of binary relations on an infinite set, infinite symmetric and dual symmetric inverse monoids, and infinite partition monoids; see for example \cite{HHR1998,HP2012,HRH1998,HHMR2003,East2014}.  
While the lower bounds given in items~(ii)--(iv) may seem crude, they are actually exact values in many of the examples just mentioned; this is also the case when $M$ is an infinite partial Brauer monoid (see Theorems \ref{thm:PBXSX}, \ref{thm:PBXEX} and~\ref{thm:PBXFX}).
\end{rem}

\begin{rem}
If a monoid $M$ satisfies $\G_L(M)\not=\G(M)$ (or equivalently $\G_R(M)\not=\G(M)$, by Lemma \ref{lem:equiv}), then clearly $\G_L(M)\not=M$ and $\G_R(M)\not=M$, and so trivially 
\[
\relrank{\G_L(M)}{\G(M)} \COMMA \relrank{\G_R(M)}{\G(M)} \COMMA \relrank{M}{\G_L(M)} \COMMA \relrank{M}{\G_R(M)}
\]
are all non-zero.  We will see in Theorems \ref{thm:PBXGLX} and \ref{thm:GLXSX} that when $M$ is an infinite partial Brauer monoid, $\relrank{M}{\G_L(M)}$ takes on its minimum possible value of $1$, whereas $\relrank{\G_L(M)}{\G(M)}$ depends on the value of $|X|$.
\end{rem}

Lemma \ref{lem:rankMGE} concerned monoids with one-sided units that are not two-sided.  The next result gives some information about relative ranks in monoids where all one-sided units are two-sided units.  For such a monoid~$M$, we may give a fairly specific formula concerning $\relrank{M}{E(M)}=\relrank{M}{\E(M)}$.  The key property used in the proof is that $M\sm\G(M)$ is an ideal (cf.~Lemma \ref{lem:equiv}).

\begin{lemma}\label{lem:rankME}
Let $M$ be a monoid, and write $G=\G(M)$, $E=\E(M)$ and $F=\F(M)$.  Suppose also that $\G_L(M)=G$ (or, equivalently by Lemma \ref{lem:equiv}, that $\G_R(M)=G$).  
\bit
\itemit{i} If $M=\la U\ra$, then $G=\la G\cap U\ra$.
\itemit{ii} If $G\not=\{1\}$, then $M=\la E\cup U\ra$ if and only if the sets $U_1=G\cap U$ and $U_2=U\sm G$ satisfy $G=\la U_1\ra$ and $M=\la F\cup U_2\ra$.
\itemit{iii} If $G\not=\{1\}$, then $\relrank ME=\rank(G)+\relrank MF$.
\eit
\end{lemma}

\pf
(i).  Suppose $M=\la U\ra$.  Let $g\in G$ be arbitrary, and consider an expression $g=u_1\cdots u_k$, where $u_1,\ldots,u_k\in U$.  By Lemma \ref{lem:equiv}, $M\sm G$ is an ideal of $M$, so it follows that all of the $u_i$ belong to $G$, and so to $G\cap U$: i.e., $g\in\la G\cap U\ra$.  This shows that $G\sub\la G\cap U\ra$; the reverse containment is clear.

\pfitem{ii}  Suppose $G\not=\{1\}$.  If $G=\la U_1\ra$ and $M=\la F\cup U_2\ra$, then 
\[
M=\la F\cup U_2\ra=\la\la E\cup G\ra\cup U_2\ra=\la E\cup G\cup U_2\ra=\la E\cup\la U_1\ra\cup U_2\ra=\la E\cup U_1\cup U_2\ra=\la E\cup U\ra.
\]
Conversely, suppose $M=\la E\cup U\ra$.  By part (i), $G$ is generated by $G\cap(E\cup U)=(G\cap E)\cup(G\cap U)=\{1\}\cup U_1$; the assumption that $G\not=\{1\}$ gives $G=\la U_1\ra$.  Then also
\[
M=\la E\cup U_1\cup U_2\ra=\la E\cup\la U_1\ra\cup U_2\ra=\la E\cup G\cup U_2\ra=\la\la E\cup G\ra\cup U_2\ra=\la F\cup U_2\ra.
\]
(iii).  Suppose $G\not=\{1\}$.  If $M=\la E\cup U\ra$ with $|U|=\relrank ME$, then with $U_1$ and $U_2$ as in part (ii),
\[
\relrank ME=|U|=|U_1|+|U_2|\geq\rank(G)+\relrank MF.
\]
Conversely, if $V_1\sub G$ and $V_2\sub F$ satisfy $G=\la V_1\ra$, $M=\la F\cup V_2\ra$, $|V_1|=\rank(G)$ and $|V_2|=\relrank MF$, then part (ii) gives $M=\la E\cup V_1\cup V_2\ra$, and so $\relrank ME\leq|V_1|+|V_2|=\rank(G)+\relrank MF$.
\epf

\begin{rem}\label{rem:rankME1}
Suppose the monoid $M$ satisfies $\G_L(M)=G$ (equivalently, $\G_R(M)=G$), using the abbreviations of Lemma \ref{lem:rankME}. 
\bit
\itemnit{i} If $G=\{1\}$, then $F=E$, and the conclusion of Lemma \ref{lem:rankME}(iii) says $\relrank ME=1+\relrank ME$, which can only be true if $\relrank ME$ is infinite.  We could get around this by replacing $\rank(G)$ with the smallest size of a \emph{monoid} generating set for $G$ (which coincides with $\rank(G)$ if $G\not=\{1\}$).
\itemnit{ii} If $M=G\not=\{1\}$, then $E=\{1\}$ and $F=M$, so Lemma \ref{lem:rankME}(iii) reduces to $\rank(M)=\rank(G)$.
\itemnit{iii} If $G\not=\{1\}$, and if $U\sub M$ with $|U|=\relrank ME<\aleph_0$, then $M=\la E\cup U\ra$ if and only if the sets $U_1$ and $U_2$ from Lemma \ref{lem:rankME}(ii) additionally satisfy $|U_1|=\rank(G)$ and $|U_2|=\relrank MF$.
\eit
\end{rem}

\begin{rem}\label{rem:rankME2}
For any monoid $M$, Lemma \ref{lem:FFLFR} shows that Lemma \ref{lem:rankME} applies to  $F_L=\F_L(M)$, ${F_R=\F_R(M)}$ and $F=\F(M)$.  Thus, if $G=\G(M)\not=\{1\}$, then 
\[
\relrank{F_L}E=\rank(G)+\relrank{F_L}F \COMMa \relrank{F_R}E=\rank(G)+\relrank{F_R}F \COMMa \relrank FE=\rank(G),
\]
where we have used $\relrank FF=0$ in the last of these.  Moreover, it quickly follows from Lemma \ref{lem:rankME}(iii) that $F=\la E\cup U\ra$ if and only if $G=\la G\cap U\ra$.
\end{rem}

\begin{rem}
It is possible to develop the ideas in \cite{East2012} in order to obtain formulae for ${\relrank{\F(M)}{\G(M)}}$, for an arbitrary monoid $M$, in terms of the minimal size of a subset $U\sub E(M)$ for which $E(M)\sm\{1\}$ is contained in the subsemigroup of $M$ generated by the set $\set{g^{-1}eg}{e\in U,\ g\in \G(M)}$.  However, we will not pursue this idea here.
\end{rem}

In this section, we have considered submonoids of a monoid $M$ generated by various combinations of $E(M)$, $\G(M)$, $\G_L(M)$ and $\G_R(M)$.  The only such submonoids not considered so far are those generated by all one-sided units, or by all idempotents and all one-sided units.  Accordingly, we may define
\[
\G_{LR}(M)=\la \G_L(M)\cup \G_R(M)\ra \AND \F_{LR}(M)=\la E(M)\cup \G_L(M)\cup \G_R(M)\ra.
\]
There does not appear to be a factorisation result akin to Lemma \ref{lem:factorisation} for either of these monoids.  We also cannot establish any positive lower bound on the values of $\relrank M{\G_{LR}(M)}$ or $\relrank M{\F_{LR}(M)}$ in general; for example, Corollary \ref{cor:PB=BB} below shows that when $M$ is an infinite partial Brauer monoid (as defined in Section~\ref{sect:PBX}), $\G_{LR}(M)=\F_{LR}(M)=M$.
This latter property does not hold in general, however; for example, if $M$ is any non-trivial additive monoid of non-negative real numbers, then $\G_{LR}(M)=\F_{LR}(M)=\{0\}\not=M$.
%
%

\begin{figure}[ht]
\begin{center}
\begin{tikzpicture}[scale=.8]
\node[rounded corners,rectangle,draw,fill=blue!20] (M) at (0,10) {\small $M$};
\node[rounded corners,rectangle,draw,fill=blue!20] (FLR) at (0,8) {\small $\F_{LR}(M)$};
\node[rounded corners,rectangle,draw,fill=blue!20] (FL) at (-3,6) {\small $\F_L(M)$};
\node[rounded corners,rectangle,draw,fill=blue!20] (GLR) at (0,6) {\small $\G_{LR}(M)$};
\node[rounded corners,rectangle,draw,fill=blue!20] (FR) at (3,6) {\small $\F_R(M)$};
\node[rounded corners,rectangle,draw,fill=blue!20] (GL) at (-3,4) {\small $\G_L(M)$};
\node[rounded corners,rectangle,draw,fill=blue!20] (F) at (0,4) {\small $\F(M)$};
\node[rounded corners,rectangle,draw,fill=blue!20] (GR) at (3,4) {\small $\G_R(M)$};
\node[rounded corners,rectangle,draw,fill=blue!20] (G) at (0,2) {\small $\G(M)$};
\node[rounded corners,rectangle,draw,fill=blue!20] (E) at (3,2) {\small $\E(M)$};
\node[rounded corners,rectangle,draw,fill=blue!20] (1) at (0,0) {\small $\{1\}$};
\draw (1)--(G)--(F)--(E)--(1);
\fill[white] (1.5,3)circle(.15);
\draw (M)--(FLR)--(FL)--(GL)--(G)--(GR)--(FR)--(FLR);
\draw (GL)--(GLR)--(GR);
\fill[white] (-1.5,5) circle(.15);
\fill[white] (1.5,5) circle(.15);
\draw (FL)--(F)--(FR);
\draw (FLR)--(GLR);
\end{tikzpicture}
\end{center}
\vspace{-.5cm}
\caption{The part of the submonoid lattice of $M$ containing the submonoids considered in Section \ref{sect:monoids}.}
\label{fig:lattice}
\end{figure}
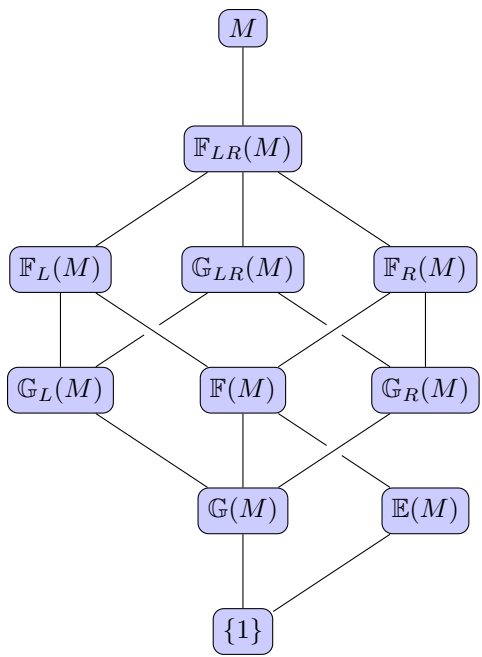

The submonoids of $M$ considered in this section, as well as the inclusion relations satisfied between them, are shown in Figure \ref{fig:lattice}.  

\section{Partial Brauer monoids}\label{sect:PBX}

We now introduce the main objects of our study: the partial Brauer monoids $\PB_X$.  Here we describe the elements and product of $\PB_X$, introduce a number of important parameters, and prove several inequalities that will be used frequently in the remainder of the article.

\begin{figure}[ht]
   \begin{center}
\begin{tikzpicture}[scale=.6]
\begin{scope}[shift={(0,0)}]	
\uvs{1,...,12}
\lvs{1,...,12}
\stlines{2/2,12/11}
\uarcs{1/3,6/8,10/11}
\uarcx59{.7}
\darcs{4/5,8/9}
\darcx36{.7}
\darcx7{10}{.7}
\foreach \x in {1,...,12} {\node () at (\x,2.5) {\tiny $\x$};}
\foreach \x in {1,...,12} {\node () at (\x,-.5) {\tiny $\phantom{'}\x'$};}
\end{scope}
\end{tikzpicture}
\vspace{-4mm}
    \caption{An element of $\PB_X$, where $X=\{1,\ldots,12\}$.}
    \label{fig:alpha}
   \end{center}
 \end{figure}
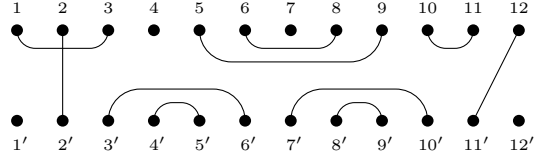

Let $X$ be an arbitrary set, and let $X'=\set{x'}{x\in X}$ be a disjoint copy of $X$.  A \emph{Brauer graph} is a graph with vertex set $X\cup X'$ in which every vertex has degree at most $1$; a Brauer graph is \emph{full} if every vertex has degree equal to $1$.  We write $\PB_X$ for the set of all Brauer graphs, and $\B_X$ for the set of all full Brauer graphs, on vertex set $X\cup X'$.  When drawing Brauer graphs, we draw the vertices from $X$ on an upper row, with those from $X'$ on a lower row directly below.  As an example with $X=\{1,\ldots,12\}$, the Brauer graph with edge set
$
\big\{
\{2,2'\}, \{12,11'\}, 
\{1,3\}, \{5,9\}, \{6,8\}, \{10,11\}, 
\{3',6'\}, \{4',5'\}, \{7',10'\}, \{8',9'\}
\big\}
$
is depicted in Figure \ref{fig:alpha}.

\begin{figure}[ht]
   \begin{center}
\begin{tikzpicture}[scale=.36]
\begin{scope}[shift={(0,0)}]	
\buvs{1,...,7}
\blvs{1,...,7}
\stline11
\uarcs{2/3,4/5,6/7}
\darcs{2/3,4/5,6/7}
\draw[dotted] (8,2)--(12,2);
\draw[dotted] (8,0)--(12,0);
\draw(0.6,1)node[left]{$\al=$};
\draw[->](13.5,-1)--(15.5,-1);
\end{scope}
\begin{scope}[shift={(0,-4)}]	
\buvs{1,...,7}
\blvs{1,...,7}
\uarcs{1/2,3/4,5/6}
\darcs{1/2,3/4,5/6}
\darchalf7{7.5}
\uarchalf7{7.5}
\colv{7.5}1{white}
\draw[dotted] (8,2)--(12,2);
\draw[dotted] (8,0)--(12,0);
\draw(0.6,1)node[left]{$\be=$};
\end{scope}
\begin{scope}[shift={(16,-1)}]	
\buvs{1,...,7}
\blvs{1,...,7}
\stline11
\uarcs{2/3,4/5,6/7}
\darcs{2/3,4/5,6/7}
\draw[dotted] (8,2)--(12,2);
\draw[dotted] (8,0)--(12,0);
\draw[->](13.5,0)--(15.5,0);
\end{scope}
\begin{scope}[shift={(16,-3)}]	
\buvs{1,...,7}
\blvs{1,...,7}
\uarcs{1/2,3/4,5/6}
\darcs{1/2,3/4,5/6}
\darchalf7{7.5}
\uarchalf7{7.5}
\colv{7.5}1{white}
\draw[dotted] (8,2)--(12,2);
\draw[dotted] (8,0)--(12,0);
\end{scope}
\begin{scope}[shift={(32,-2)}]	
\buvs{1,...,7}
\blvs{1,...,7}
\uarcs{2/3,4/5,6/7}
\darcs{1/2,3/4,5/6}
\darchalf7{7.5}
\draw[dotted] (8,2)--(12,2);
\draw[dotted] (8,0)--(12,0);
\draw(12.4,1)node[right]{$=\al\be$};
\end{scope}
\end{tikzpicture}
    \caption{Two Brauer graphs $\al,\be\in\PB_{\N}$ (left), their product $\al\be\in\PB_\N$ (right), and the product graph $\Pi(\al,\be)$ (centre).}
    \label{fig:product2}
   \end{center}
 \end{figure}
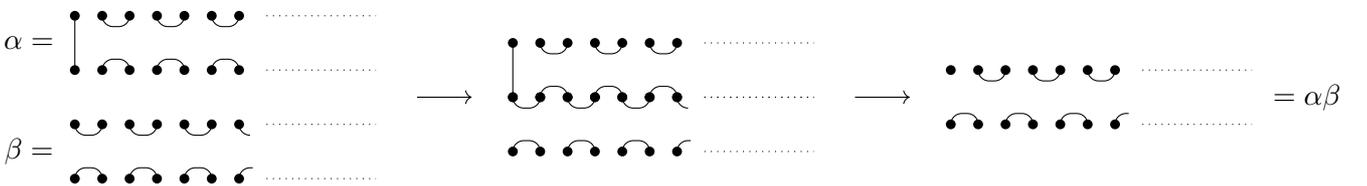

The set $\PB_X$ forms a monoid, called the \emph{partial Brauer monoid}, under a product defined as follows.  Let $\al,\be\in\PB_X$.  First, let $X''=\set{x''}{x\in X}$ be a second disjoint copy of $X$.  Let $\al^\vee$ be the graph obtained by changing each lower vertex $x'$ from $\al$ to $x''$; similarly, let $\be^\wedge$ be the graph obtained by changing each upper vertex $x$ from $\be$ to $x''$.  Now let $\Pi(\al,\be)$ be the graph on vertex set $X\cup X'\cup X''$ with all the edges from both $\al^\vee$ and $\be^\wedge$.  We call $\Pi(\al,\be)$ the \emph{product graph} associated to $\al,\be$, and we note that $\Pi(\al,\be)$ might contain pairs of parallel edges (one coming from $\al$ and one from $\be$).  Finally, $\al\be$ is the graph with vertex set $X\cup X'$, and an edge $\{x,y\}$ whenever $x,y\in X\cup X'$ are distinct and belong to the same connected component of $\Pi(\al,\be)$.  Figures \ref{fig:product1} and \ref{fig:product2} give two example calculations, for finite and (countably) infinite $X$, respectively.

\begin{figure}[ht]
   \begin{center}
\begin{tikzpicture}[scale=.36]
\begin{scope}[shift={(0,0)}]	
\buvs{1,...,12}
\blvs{1,...,12}
\stlines{2/2,12/11}
\uarcs{1/3,6/8,10/11}
\uarcx59{.7}
\darcs{4/5,8/9}
\darcx36{.7}
\darcx7{10}{.7}
\draw(0.6,1)node[left]{$\al=$};
\draw[->](13.5,-1)--(15.5,-1);
\end{scope}
\begin{scope}[shift={(0,-4)}]	
\buvs{1,...,12}
\blvs{1,...,12}
\stlines{4/7,5/4,6/5}
\uarcs{2/3,8/9,10/12}
\darcs{1/3,8/9,11/12}
\draw(0.6,1)node[left]{$\be=$};
\end{scope}
\begin{scope}[shift={(16,-1)}]	
\buvs{1,...,12}
\blvs{1,...,12}
\stlines{2/2,12/11}
\uarcs{1/3,6/8,10/11}
\uarcx59{.7}
\darcs{4/5,8/9}
\darcx36{.7}
\darcx7{10}{.7}
\draw[->](13.5,0)--(15.5,0);
\end{scope}
\begin{scope}[shift={(16,-3)}]	
\buvs{1,...,12}
\blvs{1,...,12}
\stlines{4/7,5/4,6/5}
\uarcs{2/3,8/9,10/12}
\darcs{1/3,8/9,11/12}
\end{scope}
\begin{scope}[shift={(32,-2)}]	
\buvs{1,...,12}
\blvs{1,...,12}
\stlines{2/5}
\uarcs{1/3,6/8,10/11}
\uarcx59{.7}
\darcs{1/3,4/7,8/9,11/12}
\draw(12.4,1)node[right]{$=\al\be$};
\end{scope}
\end{tikzpicture}
    \caption{Two Brauer graphs $\al,\be\in\PB_X$ with $|X|=12$ (left), their product $\al\be\in\PB_X$ (right), and the product graph $\Pi(\al,\be)$ (centre).}
    \label{fig:product1}
   \end{center}
 \end{figure}
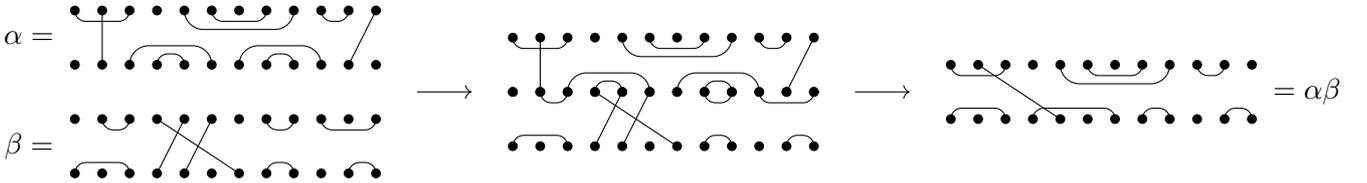

The above product is associative, so $\PB_X$ is a semigroup.  Denote by $1$ the (full) Brauer graph with edge set $\bigset{\{x,x'\}}{x\in X}$.  It is easy to see that $1$ is an identity element, so $\PB_X$ is indeed a monoid.  If~$X$ is finite, then the set $\B_X$ of all full Brauer graphs is a submonoid of $\PB_X$, known as the \emph{Brauer monoid}.  If~$X$ is infinite, then $\B_X$ is not closed under the product.  Figure \ref{fig:product2} exemplifies this last assertion; there, $\al$ and $\be$ are full, but $\al\be$ is not.  In fact, we will see in Corollary \ref{cor:PB=BB} below that \emph{every} element of infinite~$\PB_X$ is the product of two elements from $\B_X$.

A number of parameters associated to Brauer graphs will play a crucial role in all that follows.  First, we note that the connected components of a Brauer graph $\al\in\PB_X$ all have one of the following forms:
\bit
\item $\{x,y'\}$ for distinct $x,y\in X$ --- a \emph{transversal} of $\al$,
\item $\{x,y\}$ for distinct $x,y\in X$ --- an \emph{upper hook} of $\al$,
\item $\{x',y'\}$ for distinct $x,y\in X$ --- a \emph{lower hook} of $\al$,
\item $\{x\}$ for some $x\in X$ --- an \emph{upper singleton} of $\al$,
\item $\{x'\}$ for some $x\in X$ --- a \emph{lower singleton} of $\al$.
\eit
We write $t(\al)$, $h(\al)$, $h^*(\al)$, $s(\al)$ and $s^*(\al)$ for the number of transversals, upper hooks, lower hooks, upper singletons and lower singletons of $\al$, respectively.  Note that $0\leq t(\al),s(\al),s^*(\al)\leq|X|$ and that $0\leq h(\al),h^*(\al)\leq\frac12|X|$, with the ``$\frac12$'' being unnecessary if $X$ is infinite.

We define the \emph{domain} and \emph{codomain} of $\al$ to be the sets
\begin{align*}
\dom(\al) &= \set{x\in X}{x\text{ belongs to a transversal of }\al} , \\
\codom(\al) &= \set{x\in X}{x'\text{ belongs to a transversal of }\al},
\end{align*}
respectively, noting that $|{\dom(\al)}|=|{\codom(\al)}|=t(\al)$; elsewhere in the literature, the cardinal $t(\al)$ is sometimes called the \emph{rank} or \emph{propagating number} of $\al$ and denoted $\rank(\al)$ or $\operatorname{pn}(\al)$; see for example \cite{HR2005,DEG2017}.
It is easy to see that 
\[
\dom(\al\be)\sub\dom(\al) \AND \codom(\al\be)\sub\codom(\be) \qquad\text{for any $\al,\be\in\PB_X$.}
\]
If $x\in\dom(\al)$, we write $x\al$ for the unique element of $\codom(\al)$ for which $\{x,(x\al)'\}$ is a transversal of~$\al$.  If $x\in\codom(\al)$, we write $x\al^{-1}$ for the unique element of $\dom(\al)$ for which $\{x\al^{-1},x'\}$ is a transversal of $\al$.  Note that if $\al,\be\in\PB_X$, and if $x\in\dom(\al)$ is such that $x\al\in\dom(\be)$, then $x\in\dom(\al\be)$ and $x(\al\be)=(x\al)\be$; a dual statement holds for codomains and preimages.  Note, however, that it is not necessary to have $x\al\in\dom(\be)$ in order for $x\in\dom(\al\be)$ to hold; indeed, a transversal of $\al\be$ could arise from a path of length greater than $2$ in the product graph $\Pi(\al,\be)$; see Figure \ref{fig:product1}, for example, where $2\al\not\in\dom(\be)$, even though $2\in\dom(\al\be)$.  If $Y\sub\dom(\al)$ and $Z\sub\codom(\al)$, we will write $Y\al=\set{y\al}{y\in Y}$ and~$Z\al^{-1}=\set{z\al^{-1}}{z\in Z}$.

We also define the \emph{defect} and \emph{codefect} sets and cardinals of $\al$ by
\[
\Defect(\al)=X\sm\dom(\al) \COMMA \Codef(\al)=X\sm\codom(\al) \COMMA \defect(\al)=|{\Defect(\al)}| \COMMA \codef(\al)=|{\Codef(\al)}|.
\]
Note that $\defect(\al) = 2h(\al)+s(\al)$ is the number of points from $X$ that do not belong to a transversal of $\al$, while $\codef(\al) = 2h^*(\al)+s^*(\al)$ is the number of points from $X'$ that do not belong to a transversal.  Since $X=\dom(\al)\sqcup\Defect(\al)=\codom(\al)\sqcup\Codef(\al)$, we have $t(\al)+\defect(\al)=t(\al)+\codef(\al)=|X|$.  Thus, we immediately deduce the following (which does not hold for infinite $X$).

\begin{lemma}\label{lem:finitedefect}
If $X$ is a finite set, then $\defect(\al)=\codef(\al)$ for all $\al\in\PB_X$. \epfres
\end{lemma}

We now describe a convenient \emph{tableau-style} notation for the elements of $\PB_X$.  For $A\sub X$, we write $A'=\set{a'}{a\in A}$.  Let $\al\in\PB_X$, and suppose the transversals, upper hooks and lower hooks of $\al$ are $\bigset{\{a_i,b_i'\}}{i\in I}$, $\bigset{C_j}{j\in J}$ and $\bigset{D_k'}{k\in K}$.  We then write 
\[
\al=\partn{a_i}{C_j}{b_i}{D_k}_{i\in I,\ j\in J,\ k\in K}.
\]
Sometimes we abbreviate this to $\al=\partn{a_i}{C_j}{b_i}{D_k}$, with the indexing sets $I,J,K$ being implied rather than explicitly stated.  Note that with this notation, we have
\[
t(\al)=|I| \COMMA h(\al)=|J| \COMMA h^*(\al)=|K| \COMMA \dom(\al)=\set{a_i}{i\in I} \COMMA \codom(\al)=\set{b_i}{i\in I}.
\]
Note also that the singletons of $\al$ are not listed explicitly in the above notation, although they are implied by it.  We will sometimes use abbreviations of the above notation: we may write $\al=\partn{a_i}{C_j}{b_i}{}$ or $\al=\partn{a_i}{}{b_i}{D_k}$ if $h^*(\al)=0$ or $h(\al)=0$, respectively.  If $h(\al)=h^*(\al)=0$, then we may write $\al=\partnVI{a_i}{b_i}$.  

On a small number of occasions, we will wish to use similar notation, but list \emph{all} of the non-transversals instead of only the hooks.  To do so, if $\al\in\PB_X$, we will write
\[
\al=\sqpartn{a_i}{C_j}{b_i}{D_k}_{i\in I,\ j\in J,\ k\in K},
\]
or just $\sqpartn{a_i}{C_j}{b_i}{D_k}$, which indicates that $\set{C_j}{j\in J}$ and $\set{D_k'}{k\in K}$ are the entire sets of upper and lower non-transversals, respectively, including hooks \emph{and} singletons.

There is also an important anti-involution ${}^*:\PB_X\to\PB_X:\al\mt\al^*$.  With $\al=\partn{a_i}{C_j}{b_i}{D_k}$ as above, we define $\al^*=\partn{b_i}{D_k}{a_i}{C_j}$.  It is easy to check that
\begin{equation}\label{eq:al*}
(\al^*)^*=\al \COMMA \al=\al\al^*\al \COMMA (\al\be)^*=\be^*\al^* \qquad\text{for all $\al,\be\in\PB_X$,}
\end{equation}
so that $\PB_X$ is a \emph{regular $*$-semigroup} in the sense of Nordahl and Scheiblich \cite{NS1978}.  We also have several obvious identities such as
\[
\dom(\al^*)=\codom(\al^*) \COMMA t(\al^*)=t(\al) \COMMA h(\al^*)=h^*(\al) \COMMA s(\al^*)=s^*(\al),
\]
and so on.

In the remainder of this section, we establish a number of inequalities involving the above parameters.  In order to prove them, and for later usage, it will be convenient to list the kinds of connected components that can arise in a product graph $\Pi(\al,\be)$, where $\al,\be\in\PB_X$.  Suppose $\C$ is such a component.  We call $\C$ \emph{trivial} if it is contained in either $X$ or $X'$ or $X''$.  
\bit
\item If $\C\sub X$, then it is an upper non-transversal of $\al$, and remains in the product $\al\be$.
\item If $\C\sub X'$, then it is a lower non-transversal of $\be$, and remains in the product $\al\be$.
\item If $\C\sub X''$, then it is either a loop or a path.  Of course loops involve only finitely many vertices, but paths could be finite or infinite; the latter can extend infinitely in one or two directions.  Such components are essentially ``forgotten'' when we form the product $\al\be$.  (These play an important role, however, in the partial Brauer \emph{algebras}; see for example \cite{MarMaz2014}.)
\eit
We call $\C$ \emph{non-trivial} if it involves at least one vertex from $X''$, and at least one from $X\cup X'$.  A non-trivial component might involve several (even infinitely many) vertices from $X''$, but involves at most two vertices from $X\cup X'$.  There are five types of non-trivial components.
\bit
\item If $\C$ is non-trivial and involves one vertex from $X$ and one from $X'$, then it has the form 
\begin{equation}\label{eq:t}
x\lar{\al} z_1'' \lar{\be} z_2'' \lar{\al} \cdots \lar{\al} z_{2k+1}'' \lar{\be} y' \qquad\text{for some $k\geq0$ and some $z_1,\ldots,z_{2k+1}\in X$.}
\end{equation}
In this case, $\C$ gives rise to the transversal $\{x,y'\}$ in the product $\al\be$.

\item If $\C$ is non-trivial and involves two vertices from $X$, then it has the form 
\begin{equation}\label{eq:h}
x\lar{\al} z_1'' \lar{\be} z_2'' \lar{\al} \cdots \lar{\be} z_{2k}'' \lar{\al} y \qquad\text{for some $k\geq1$ and some $z_1,\ldots,z_{2k}\in X$.}
\end{equation}
In this case, $\C$ gives rise to the upper hook $\{x,y\}$ in the product $\al\be$.

\item If $\C$ is non-trivial and involves one vertex from $X$ and none from $X'$, then it has the form 
\begin{equation}\label{eq:s}
x\lar{\al} z_1'' \lar{\be} z_2'' \lar{\al} \cdots  \qquad\text{for some $z_1,z_2,\ldots\in X$.}
\end{equation}
In this case, $\C$ might be infinite in length, or may terminate at a point corresponding to a lower singleton of $\al$ or an upper singleton of $\be$, but it always gives rise to the upper singleton $\{x\}$ in the product $\al\be$.

\item If $\C$ is non-trivial and involves two vertices from $X'$, then it has the form 
\begin{equation}\label{eq:h*}
x'\lar{\be} z_1'' \lar{\al} z_2'' \lar{\be} \cdots \lar{\al} z_{2k}'' \lar{\be} y' \qquad\text{for some $k\geq1$ and some $z_1,\ldots,z_{2k}\in X$.}
\end{equation}
In this case, $\C$ gives rise to the lower hook $\{x',y'\}$ in the product $\al\be$.

\item If $\C$ is non-trivial and involves one vertex from $X'$ and none from $X$, then it has the form 
\begin{equation}\label{eq:s*}
x'\lar{\be} z_1'' \lar{\al} z_2'' \lar{\be} \cdots  \qquad\text{for some $z_1,z_2,\ldots\in X$.}
\end{equation}
Again, $\C$ might be finite or infinite in this case, but it always gives rise to the lower singleton $\{x'\}$ in the product~$\al\be$.

\eit

\begin{lemma}\label{lem:shdef}
Let $X$ be an arbitrary set, and let $\al,\be\in\PB_X$.  Then
\bit\bmc2
\itemit{i} $s(\al) \leq s(\al\be)$ and $s^*(\be) \leq s^*(\al\be)$,
\itemit{ii} $h(\al) \leq h(\al\be) \leq h(\al)+h(\be)$,
\itemit{iii} $h^*(\be) \leq h^*(\al\be) \leq h^*(\al)+h^*(\be)$,
\itemit{iv} $\defect(\al) \leq \defect(\al\be) \leq \defect(\al)+\defect(\be)$,
\itemit{v} $\codef(\be) \leq \codef(\al\be) \leq \codef(\al)+\codef(\be)$,
\itemit{vi} $t(\al\be)\leq t(\al)$ and $t(\al\be)\leq t(\be)$.
\emc\eit
\end{lemma}

\pf
(i).  Every upper singleton $\{x\}$ of $\al\be$ is either an upper singleton of $\al$, or else arises from some non-trivial component in the product graph $\Pi(\al,\be)$ of the form \eqref{eq:s}.  Thus, if there are $\mu$ of the latter kind of component, then $s(\al\be)=s(\al)+\mu\geq s(\al)$.  The statement concerning $s^*$ is dual.

\pfitem{ii}  Similarly, every upper hook $\{x,y\}$ of $\al\be$ is either an upper hook of $\al$ or else arises from a non-trivial component in $\Pi(\al,\be)$ of the form \eqref{eq:h}.  Thus, if there are $\nu$ of the latter kind of component, then $h(\al\be)=h(\al)+\nu\geq h(\al)$.  Since any component of the form \eqref{eq:h} involves at least one upper hook of $\be$, and since each upper hook of $\be$ is involved in at most one such component, we obtain $\nu\leq h(\be)$.  Thus, $h(\al\be)=h(\al)+\nu\leq h(\al)+h(\be)$.

\pfitem{iv}  With $\mu$ and $\nu$ as above, $\defect(\al\be)=s(\al\be)+2h(\al\be) = s(\al)+\mu+2(h(\al)+\nu) = \defect(\al)+\mu+2\nu \geq \defect(\al)$.  It remains to show that $\mu+2\nu\leq\defect(\be)$.  Since $|{\Defect(\al\be)\sm\Defect(\al)}|=\mu+2\nu$, we may prove the latter by constructing an injective map $\phi:\Defect(\al\be)\sm\Defect(\al)\to\Defect(\be)$.  With this in mind, let $x\in\Defect(\al\be)\sm\Defect(\al)$.  If $\{x\}$ is a singleton of~$\al\be$, then there is a component in $\Pi(\al,\be)$ of the form \eqref{eq:s}, and we define $x\phi=z_1$.  If $x$ belongs to a hook~$\{x,y\}$ of~$\al\be$, then also $y\in\Defect(\al\be)\sm\Defect(\al)$, and there is a component in $\Pi(\al,\be)$ of the form~\eqref{eq:h}; we then define $x\phi=z_1$ and $y\phi=z_{2k}$.

\pfitem{iii) and (v}  These are dual to (ii) and (iv), respectively.

\pfitem{vi}  Any transversal $\{x,y'\}$ of $\al\be$ arises from a non-trivial component in $\Pi(\al,\be)$ of the form \eqref{eq:t}.  Such a component involves the transversals $\{x,z_1'\}$ from $\al$ and $\{z_{2k+1},y'\}$ from $\be$.  The result follows immediately.
\epf
 
\begin{rem}
Lemma \ref{lem:shdef} has no statement of the form $s(\al\be)\leq s(\al)+s(\be)$ or ${s^*(\al\be)\leq s^*(\al)+s^*(\be)}$, because these need not hold.  Examples where $s(\al\be)>s(\al)+s(\be)$ may easily be constructed, even with $|X|=2$.
\end{rem}

The next simple corollary of Lemma \ref{lem:shdef} will be used frequently.  This result, and many more to come, involve cardinals $\mu$ such that $\mu=1$ or $\mu\geq\aleph_0$.  The crucial property of such cardinals is that they cannot be written as a finite sum of smaller cardinals.

\begin{cor}\label{cor:tq}
Suppose $X$ is an arbitrary set, let $\al_1,\ldots,\al_k\in\PB_X$, let $q$ denote any of $h$, $h^*$, $\defect$ or~$\codef$, and suppose $\mu$ is an arbitrary cardinal.  Then 
\bit
\itemit{i}  $t(\al_1\cdots\al_k)\geq\mu \implies t(\al_i)\geq\mu$ for all $i$,
\itemit{ii} if $\mu=1$ or $\mu\geq\aleph_0$, then $q(\al_1\cdots\al_k)\geq\mu \implies q(\al_i)\geq\mu$ for some $i$.
\eit

\end{cor}

\pf
(i).  For any $i$, the two assertions of Lemma \ref{lem:shdef}(vi) give
\[
\mu\leq t(\al_1\cdots\al_{i-1}\al_i\al_{i+1}\cdots\al_k)\leq t(\al_1\cdots\al_{i-1}\al_i) \leq t(\al_i).
\]
(ii).  If $q(\al_i)<\mu$ for all $i$, then, by the relevant part of Lemma \ref{lem:shdef},
$
{q(\al_1\cdots\al_k)\leq q(\al_1)+\cdots+q(\al_k)<\mu},
$
contradicting $q(\al_1\cdots\al_k)\geq\mu$.
\epf

There are dual versions of the next three lemmas, but we will not explicitly state these.
The next result shows how Lemma \ref{lem:shdef} simplifies in the case that $\codef(\al)=0$, which, as we will see in Lemma \ref{lem:PBunits}(ii), is precisely the condition for $\al$ to be a left unit of $\PB_X$.

\begin{lemma}\label{lem:shdef2}
Let $X$ be an arbitrary set, and let $\al,\be\in\PB_X$.  If $\codef(\al)=0$, then
\bit\bmc3
\itemit{i} $s(\al\be)=s(\al)+s(\be)$,
\itemit{ii} $s^*(\al\be)=s^*(\be)$,
\item[]
\itemit{iii} $h(\al\be)=h(\al)+h(\be)$,
\itemit{iv} $h^*(\al\be)=h^*(\be)$,
\item[]
\itemit{v} $\defect(\al\be)=\defect(\al)+\defect(\be)$,
\itemit{vi} $\codef(\al\be)=\codef(\be)$,
\itemit{vii} $t(\al\be)=t(\be)$.
\emc\eit
\end{lemma}

\pf
(i).  As in the proof of Lemma \ref{lem:shdef}, we have $s(\al\be)=s(\al)+\mu$, where $\mu$ denotes the number of non-trivial components of the product graph $\Pi(\al,\be)$ of the form \eqref{eq:s}.  Such a path component either:
\bit
\itemnit{a} is infinite, or
\itemnit{b} terminates at $z_{2k+1}''$ for some $k\geq0$, where $z_{2k+1}$ is an upper singleton of $\be$, or 
\itemnit{c} terminates at $z_{2k}''$ for some $k\geq1$, where $z_{2k}'$ is a lower singleton of $\al$.
\eit
Since $\codef(\al)=0$, there are no components of type (a) or (c), and any component of type (b) must have $k=0$.  Together with the fact that $\codom(\al)=X$, it follows that the path components of the form \eqref{eq:s} are in one-one correspondence with the upper singletons of $\be$.  Thus, $\mu=s(\be)$. 

\pfitem{ii}  Since $\codom(\al)=X$, we have $\al^*\al=1$.  Lemma \ref{lem:shdef}(i) then gives $s^*(\be)\leq s^*(\al\be)\leq s^*(\al^*\al\be)=s^*(\be)$.

\pfitem{iii), (iv) and (vii}  These are proved in similar fashion to (i) and (ii).  

\pfitem{v) and (vi}  These follow from (i)--(iv), with $\defect(\ga)=s(\ga)+2h(\ga)$ and $\codef(\ga)=s^*(\ga)+2h^*(\ga)$.
\epf

Lemma \ref{lem:shdef}(iv) says that $\defect(\al)\leq\defect(\al\be)$ for any $\al,\be\in\PB_X$.  The next result gives a variation on this in the case that $\codef(\al)\leq\defect(\al)$, which, as we will see in Theorem \ref{thm:FLFR}(i), is precisely the condition for $\al$ to be a product of idempotents and left units.

\begin{lemma}\label{lem:shdef4}
Let $X$ be an arbitrary set, and let $\al,\be\in\PB_X$.  If $\codef(\al)\leq\defect(\al)$, then $\defect(\be)\leq\defect(\al\be)$.
\end{lemma}

\pf
Suppose $\codef(\al)\leq\defect(\al)$.  It suffices to demonstrate the existence of an injective map
\[
\phi:\Defect(\be)\to\Defect(\al\be).
\]
By assumption, we may fix an injective map $\psi:\Codef(\al)\to\Defect(\al)$.  Let $\set{\C_i}{i\in I}$ be the set of all connected components in the product graph $\Pi(\al,\be)$ that contain a point $z''$, where $z\in\Defect(\be)$.  
We define $\phi$ by specifying its action on the sets $\set{z\in\Defect(\be)}{z''\in\C_i}$, for each $i\in I$.
\bit
\item If $\C_i$ is a trivial component (i.e., if it is contained wholly in $X''$), then all of its vertices $z''$ are such that~$z$ belongs to both $\Defect(\be)$ and $\Codef(\al)$.  We then define $z\phi=z\psi$ for all such vertices.

\item
If $\C_i$ has the form \eqref{eq:t}, then we must have $k\geq1$ (since $\C_i$ involves at least one point $z''$ with ${z\in\Defect(\be)}$).  In this case, we have
$z_1,\ldots,z_{2k}\in\Defect(\be)$, $z_{2k+1}\in\dom(\be)$, $z_1\in\codom(\al)$ and ${z_2,\ldots,z_{2k+1}\in\Codef(\al)}$.
We then define $z_j\phi=z_{j+1}\psi$ for each $1\leq j\leq 2k$.

\item
If $\C_i$ has the form \eqref{eq:h}, then 
$
z_1,\ldots,z_{2k}\in\Defect(\be)$, $z_1,z_{2k}\in\codom(\al)$ and $z_2,\ldots,z_{2k-1}\in\Codef(\al).
$
We then define $z_1\phi=z_1\al^{-1}$, $z_{2k}\phi=z_{2k}\al^{-1}$ and $z_j\phi=z_j\psi$ for each $2\leq j\leq 2k-1$.

\item
If $\C_i$ has the form \eqref{eq:s}, then (whether this component is finite or infinite)
$
z_1,z_2,\ldots\in\Defect(\be)$, ${z_1\in\codom(\al)}$ and $z_2,z_3,\ldots\in\Codef(\al).
$
We then define $z_1\phi=z_1\al^{-1}$ and $z_j\phi=z_j\psi$ for each~$j\geq2$.

\item
If $\C_i$ has the form \eqref{eq:h*}, then 
$
z_2,\ldots,z_{2k-1}\in\Defect(\be)$, $z_1,z_{2k}\in\dom(\be)$ and $z_1,\ldots,z_{2k}\in\Codef(\al).
$
We then define $z_j\phi=z_j\psi$ for each $2\leq j\leq 2k-1$.

\item
If $\C_i$ has the form \eqref{eq:s*}, then 
$
z_2,z_3,\ldots\in\Defect(\be)$, $z_1\in\dom(\be)$ and $z_1,z_2,\ldots\in\Codef(\al).
$
We then define $z_j\phi=z_j\psi$ for each $j\geq2$.

\eit
We have defined $z\phi$ for each point $z\in\Defect(\be)$, and in each case, one may check that $z\phi\in\Defect(\al\be)$.  The injectivity of $\psi$, and also of $\al^{-1}:\codom(\al)\to\dom(\al)$, ensures that $\phi$ is injective.
\epf

Note that if $\mu$ and $\nu$ are cardinals with $\nu<\mu$, then the difference $\mu-\nu$ is well defined; if $\mu$ is infinite (or if $\nu=0$), then $\mu-\nu=\mu$.


\begin{lemma}\label{lem:shdef3}
Let $X$ be an arbitrary set, and let $\al,\be\in\PB_X$.  Then
\bit
\itemit{i} $s(\be)>\codef(\al) \implies s(\al\be)\geq  s(\al)+ s(\be)-\codef(\al)$,
\itemit{ii} $h(\be)>\codef(\al) \implies h(\al\be)\geq  h(\al)+ h(\be)-\codef(\al)$,
\itemit{iii} if $\defect(\be)=1$ or $\defect(\be)\geq\aleph_0$, then $\defect(\be)>\codef(\al) \implies \defect(\al\be)\geq \defect(\be)$.
\eit
\end{lemma}

\pf
(ii).  Let $\nu$ be as in the proof of Lemma \ref{lem:shdef}(ii).  Also write $\ka=\codef(\al)$, and suppose $h(\be)>\ka$.  Since $h(\al\be)=h(\al)+\nu$, we just need to show that $\nu\geq h(\be)-\ka$.
Now, at most $\ka$ of the upper hooks of $\be$ involve one or more points from $\Codef(\al)$, so at least $h(\be)-\ka$ upper hooks of $\be$ are contained in $\codom(\al)$.  Any such upper hook of $\be$ is involved in a component of type \eqref{eq:h} in the product graph $\Pi(\al,\be)$ with $k=1$, and so uniquely determines an upper hook of $\al\be$ that is not a hook of $\al$.  Thus, $h(\be)-\ka\leq\nu$, as required.

\pfitem{i}  This is almost identical to (ii), but slightly simpler, so we omit the details.

\pfitem{iii}  Suppose $\defect(\be)>\codef(\al)$.  If $\defect(\be)=1$, then $\codef(\al)=0$, and so Lemma \ref{lem:shdef2}(v) gives 
\[
\defect(\al\be)=\defect(\al)+\defect(\be)\geq\defect(\be),
\]
completing the proof in this case.
For the remainder of the proof, we will assume that $\defect(\be)\geq\aleph_0$.  From $s(\be)+2h(\be)=\defect(\be)\geq\aleph_0$, it follows that ${\defect(\be)=\max\{s(\be),h(\be)\}}$.  We assume $\defect(\be)=s(\be)$; the $\defect(\be)=h(\be)$ case is almost identical.  Now, $s(\be)=\defect(\be)>\codef(\al)$, so part (i), above, gives $s(\al\be)\geq s(\al)+s(\be)-\codef(\al)$.  Since $s(\be)>\codef(\al)$ and $s(\be)\geq\aleph_0$, we have $s(\be)-\codef(\al)=s(\be)$.  But then $\defect(\al\be)\geq s(\al\be)\geq s(\al)+s(\be)-\codef(\al) = s(\al)+s(\be) \geq s(\be)=\defect(\be)$.
\epf

\section{Units}\label{sect:SX}

In this section, we study the one- and two-sided units of $\PB_X$.  For simplicity, we will use the abbreviations
\[
\GLX=\G_L(\PB_X) \COMMA \GRX=\G_R(\PB_X) \COMMA \S_X=\G(\PB_X)=\GLX\cap\GRX,
\]
for the monoids of all left units, all right units, or all (two-sided) units of $\PB_X$, respectively.
After characterising the elements of $\GLX$, $\GRX$ and $\S_X$ in Lemma \ref{lem:PBunits}, we calculate the relative ranks
\[
\relrank{\PB_X}{\S_X} \COMMA
\relrank{\PB_X}{\GLX} \COMMA
\relrank{\PB_X}{\GRX} \COMMA
\relrank{\GLX}{\S_X} \COMMA
\relrank{\GRX}{\S_X} ,
\]
in Theorems \ref{thm:PBXSX}, \ref{thm:PBXGLX} and \ref{thm:GLXSX}; these theorems also classify the minimal-size generating sets modulo the stated submonoids.  

We begin with a description of the units.  In what follows, the next result will often be used without explicit reference.

\begin{lemma}\label{lem:PBunits}
If $X$ is an arbitrary set, then
\bit
\itemit{i} $\GRX = \set{\al\in\PB_X}{\dom(\al)=X} = \set{\al\in\PB_X}{\defect(\al)=0}$,
\itemit{ii} $\GLX = \set{\al\in\PB_X}{\codom(\al)=X} = \set{\al\in\PB_X}{\codef(\al)=0}$,
\itemit{iii} $\S_X = \set{\al\in\PB_X}{\dom(\al)=\codom(\al)=X} = \set{\al\in\PB_X}{\defect(\al)=\codef(\al)=0}$,
\itemit{iv} $\GLX=\S_X \iff \GRX=\S_X \iff X$ is finite.
\eit
\end{lemma}

\pf
We just prove (i) and (iv), as (ii) is dual to (i), and (iii) follows from (i) and (ii).

\pfitem{i}  Let $\al\in \GRX$, so that $1=\al\be$ for some $\be\in\PB_X$.  Then $X=\dom(\al\be)\sub\dom(\al)\sub X$, so that $\dom(\al)=X$.  Conversely, if $\dom(\al)=X$, then $1=\al\al^*$, so that $\al\in \GRX$.  

\pfitem{iv}  By Lemma \ref{lem:equiv}, it suffices to show that $\PB_X$ contains a bicyclic submonoid if and only if $X$ is infinite.  If $X$ is infinite, then we take any $\al\in\PB_X$ with $\dom(\al)=X\not=\codom(\al)$, and note that $\al\al^*=1\not=\al^*\al$, so that $\{\al,\al^*\}$ generates a bicyclic submonoid.  If $X$ is finite, then $\PB_X$ cannot contain a bicyclic monoid, since bicyclic monoids are infinite.
\epf

From Lemma \ref{lem:PBunits}(iii) we recover the well-known fact that the group of units $\S_X=\G(\PB_X)$ is isomorphic to the symmetric group on $X$; cf.~\cite[Section 2]{East2014}, \cite[Section 2]{EF2012} and \cite[Lemma 2.3]{ER2018}.  
Note also that $(\GLX)^*=\set{\al^*}{\al\in\GLX}=\GRX$, and similarly $(\GRX)^*=\GLX$.  In fact, if~$M$ is any monoid with an anti-involution $M\to M:x\mt x^*$ (meaning that $(x^*)^*=x$ and $(xy)^*=y^*x^*$ for all $x,y\in M$), then $\G_L(M)^*=\G_R(M)$ and $\G_R(M)^*=\G_L(M)$.
This means that any statement concerning~$\GLX$ has a natural dual statement for $\GRX$, and the latter can be easily deduced from the former.  Thus, we will often only formulate results for one or the other of $\GLX$ or $\GRX$.

The next simple lemma will be used often.

\begin{lemma}\label{lem:action}
Let $X$ be an arbitrary set, and let $\al,\be\in\PB_X$.  Then $\be\in\S_X\al\S_X$ if and only if
\[
t(\al)=t(\be) \COMMA 
h(\al)=h(\be) \COMMA 
h^*(\al)=h^*(\be) \COMMA 
s(\al)=s(\be) \COMMA 
s^*(\al)=s^*(\be) .
\]
\end{lemma}

\pf 
Write $\al=\partn{a_i}{C_j}{b_i}{D_k}$, and let $P$ and $Q'$ be the sets of upper and lower singletons of $\al$, respectively.

\pfitem{$\Rightarrow$}  If $\be=\ga\al\de$ where $\ga,\de\in\S_X$, then $\be=\partn{a_i\ga^{-1}}{C_j\ga^{-1}}{b_i\de}{D_k\de}$, and the upper and lower singleton sets of $\be$ are $P\ga^{-1}$ and $(Q\de)'$, respectively.  Equality of the parameters is immediate.

\pfitem{$\Leftarrow$}  Assuming equality of the parameters, we may write $\be=\partn{e_i}{G_j}{f_i}{H_k}$, using the same indexing sets as for~$\al$.  We also write $R$ and $S'$ for the sets of upper and lower singletons of $\be$, respectively; by assumption, $|P|=|R|$ and $|Q|=|S|$.  We then define $\ga,\de\in\S_X$ so that
\[
e_i\ga=a_i \COMMA b_i\de=f_i \COMMA G_j\ga=C_j \COMMA D_k\de=H_k \qquad\text{for all $i\in I$, $j\in J$, $k\in K$.}
\]
Then $\ga$ must also map $R$ bijectively onto $P$, and we have $\be=\ga\al\de$.
\epf

The next result is key in what follows; it shows that infinite~$\PB_X$ may be generated by $\S_X$ along with two other Brauer graphs of a certain form.

\begin{lemma}\label{lem:alSXbe}
Let $X$ be an infinite set, and let $\al\in\GRX$ and $\be\in\GLX$ with $h^*(\al)=h(\be)=|X|$.  Then $\PB_X=\al\S_X\be$.
\end{lemma}

\pf
Since $\dom(\al)=\codom(\be)=X$, and since $h^*(\al)=h(\be)=|X|$, we may write $\al=\partn{x}{}{a_x}{B_x}$ and $\be=\partn{c_x}{D_x}{x}{}$.  For each $x\in X$, write $B_x=\{b_{x1},b_{x2}\}$ and $D_x=\{d_{x1},d_{x2}\}$.  Fix subsets $Y,Z\sub X$ such that $X=Y\sqcup Z$ and $|X|=|Y|=|Z|$.

Let $\ga\in\PB_X$.  We must show that $\ga=\al\de\be$ for some $\de\in\S_X$.  We give the definition of $\de$ in several stages; see steps (i)--(vi) below.  Write $\ga=\partn{e_i}{G_j}{f_i}{H_k}$, assuming that the indexing sets $I$, $J$ and~$K$ are disjoint (but noting that any or all of them might be empty).  For each $j\in J$ and $k\in K$, write $G_j=\{g_{j1},g_{j2}\}$ and $H_k=\{h_{k1},h_{k2}\}$.  
\bit
\itemnit{i} For every $i\in I$, we define $a_{e_i}\de=c_{f_i}$.  
\eit
Let $Y_J$ and $Y_K$ be subsets of $Y$ such that $|Y_J|=|J|$, $|Y_K|=|K|$ and $|Y\sm Y_J|=|Y\sm Y_K|=|X|$.  Write $Y_J=\set{y_j}{j\in J}$ and $Y_K=\set{y_k}{k\in K}$.  (We do not require that $Y_J$ and $Y_K$ be disjoint.)
\bit
\itemnit{ii} For each $j\in J$, we define $a_{g_{j1}}\de=d_{y_j1}$ and $a_{g_{j2}}\de=d_{y_j2}$.
\itemnit{iii} For each $k\in K$, we define $b_{y_k1}\de=c_{h_{k1}}$ and $b_{y_k2}\de=c_{h_{k2}}$.
\eit
Next, let $V$ be the set of all upper singletons of $\ga$, and $W'$ the set of all lower singletons of $\ga$, where $W\sub X$.  Let $Z_V$ and $Z_W$ be subsets of $Z$ such that $Z_V\cap Z_W=\emptyset$, and $|Z_V|=|V|\aleph_0$ and $|Z_W|=|W|\aleph_0$.  Write $Z_V=\set{z_{vn}}{v\in V,\ n\in\N}$ and $Z_W=\set{z_{wn}}{w\in W,\ n\in\N}$.  
\bit
\itemnit{iv} For each $v\in V$, we define $a_v\de=d_{z_{v0},1}$, \ $b_{z_{vn},1}\de=d_{z_{vn},2}$  and $b_{z_{vn},2}\de=d_{z_{v,n+1},1}$ for each $n\in\N$.
\itemnit{v} For each $w\in W$, we define $b_{z_{w0},1}\de=c_w$, \ $b_{z_{wn},2}\de=d_{z_{wn},1}$ for each $n\in\N$, and $b_{z_{wn},1}\de=d_{z_{w,n-1},2}$ for each $n\in\N\sm\{0\}$.
\eit
So far, $\de$ is defined to be a bijection from $\codom(\al)\cup\bigcup_{x\in Y_K\cup Z_V\cup Z_W}B_x$ to $\dom(\be)\cup\bigcup_{x\in Y_J\cup Z_V\cup Z_W}D_x$.
We denote these sets by $X_1$ and $X_2$, respectively.
Examining steps (i)--(v), note that if the definition of $\de$ is completed arbitrarily (by specifying the edges between the vertices $X_1\cup X_2'$), then each connected component of $\ga$ is a connected component of $\al\de\be$, so that $\ga=\al\de\be$.  Here we wish to show that the definition of $\de$ may be completed in such a way that $\de\in\S_X$.  Now, the complements $X\sm X_1$ and $X\sm X_2$ contain $\bigcup_{x\in Y\sm Y_K}B_x$ and $\bigcup_{x\in Y\sm Y_J}D_x$, respectively, and so $|X\sm X_1|=|X\sm X_2|=|X|$.  Thus, there is a bijection $\ve:X\sm X_1\to X\sm X_2$.
\bit
\itemnit{vi} We complete the definition of $\de$ by defining $x\de=x\ve$ for all $x\in X\sm X_1$.
\eit
Then $\de$ is indeed an element of $\S_X$, and we noted above that $\ga=\al\de\be$.  This completes the proof.
\epf

Lemma \ref{lem:alSXbe} makes no assumption about singletons of $\al,\be$; in particular, it could be the case that $\al,\be$ have no singletons at all: i.e., that $\al,\be\in\B_X$.
Among other things, the next result uses this observation to show that any Brauer graph on an infinite vertex set is a product of two full Brauer graphs.  

\begin{cor}\label{cor:PB=BB}
If $X$ is an infinite set, then
\bit
\itemit{i} $\PB_X$ is generated by its left units and right units; in fact, $\PB_X=\GRX \GLX$,
\itemit{ii} $\PB_X$ is generated by $\B_X$; in fact, $\PB_X=\B_X^2=(\B_X\cap\GRX) (\B_X\cap\GLX)$.
\eit
\end{cor}

\pf
Clearly it suffices to show that $\PB_X\sub (\B_X\cap\GRX) (\B_X\cap\GLX)$.  Let $\al,\be\in\B_X$ be such that ${\dom(\al)=\codom(\be)=X}$ and $h^*(\al)=h(\be)=|X|$.  Then for any $\ga\in\PB_X$, Lemma \ref{lem:alSXbe} gives $\ga=\al\de\be$ for some $\de\in\S_X$.  The proof concludes with the observation that $\al\de\in \B_X\cap\GRX$ and $\be\in \B_X\cap\GLX$.
\epf

\begin{figure}[ht]
\begin{center}
\begin{tikzpicture}[scale=.8]
\node[rounded corners,rectangle,draw,fill=blue!20] (FLR) at (0,8) {\small $\PB_X$};
\node[rounded corners,rectangle,draw,fill=blue!20] (FL) at (-3,6) {\small $\F_L(\PB_X)$};
\node[rounded corners,rectangle,draw,fill=blue!20] (FR) at (3,6) {\small $\F_R(\PB_X)$};
\node[rounded corners,rectangle,draw,fill=blue!20] (GL) at (-3,4) {\small $\G_L(\PB_X)$};
\node[rounded corners,rectangle,draw,fill=blue!20] (F) at (0,4) {\small $\F(\PB_X)$};
\node[rounded corners,rectangle,draw,fill=blue!20] (GR) at (3,4) {\small $\G_R(\PB_X)$};
\node[rounded corners,rectangle,draw,fill=blue!20] (G) at (0,2) {\small $\G(\PB_X)$};
\node[rounded corners,rectangle,draw,fill=blue!20] (E) at (3,2) {\small $\E(\PB_X)$};
\node[rounded corners,rectangle,draw,fill=blue!20] (1) at (0,0) {\small $\{1\}$};
\draw (1)--(G)--(F)--(E)--(1);
\fill[white] (1.5,3)circle(.15);
\draw (FLR)--(FL)--(GL)--(G)--(GR)--(FR)--(FLR) (FL)--(F)--(FR);
\node () at (6.2,4) {$\equiv$};
\begin{scope}[shift={(12,0)}]
\node[rounded corners,rectangle,draw,fill=blue!20] (FLR) at (0,8) {\small $\PB_X$};
\node[rounded corners,rectangle,draw,fill=blue!20] (FL) at (-3,6) {\small $\FLX$};
\node[rounded corners,rectangle,draw,fill=blue!20] (FR) at (3,6) {\small $\FRX$};
\node[rounded corners,rectangle,draw,fill=blue!20] (GL) at (-3,4) {\small $\GLX$};
\node[rounded corners,rectangle,draw,fill=blue!20] (F) at (0,4) {\small $\FX$};
\node[rounded corners,rectangle,draw,fill=blue!20] (GR) at (3,4) {\small $\GRX$};
\node[rounded corners,rectangle,draw,fill=blue!20] (G) at (0,2) {\small $\GX$};
\node[rounded corners,rectangle,draw,fill=blue!20] (E) at (3,2) {\small $\EX$};
\node[rounded corners,rectangle,draw,fill=blue!20] (1) at (0,0) {\small $\{1\}$};
\draw (1)--(G)--(F)--(E)--(1);
\fill[white] (1.5,3)circle(.15);
\draw (FLR)--(FL)--(GL)--(G)--(GR)--(FR)--(FLR) (FL)--(F)--(FR);
\end{scope}
\end{tikzpicture}
\end{center}
\vspace{-5mm}
\caption{The part of the submonoid lattice of $\PB_X$ containing the submonoids studied in this article; the diagram on the right displays the shorthand notation we use for the submonoids.}
\label{fig:lattice2}
\end{figure}
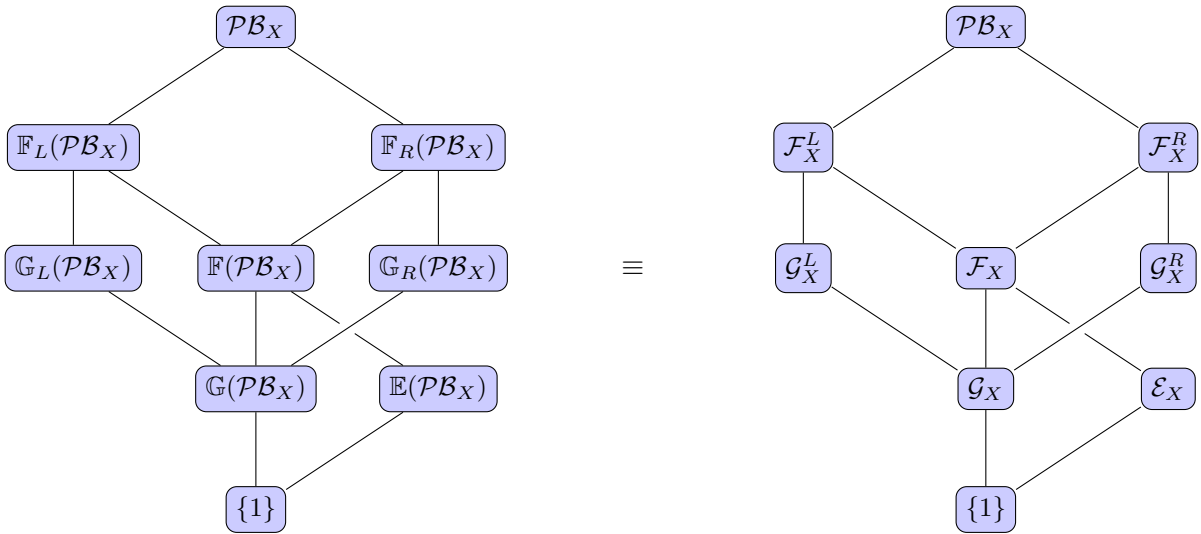

\begin{rem}\label{rem:PB=BB}
It follows from Corollary \ref{cor:PB=BB}(i) that infinite $\PB_X$ is equal to $\G_{LR}(\PB_X)$, in the notation of Section \ref{sect:monoids}.  Since also ${\G_{LR}(M)\sub\F_{LR}(M)\sub M}$ for any monoid $M$, it follows that $\PB_X=\F_{LR}(\PB_X)$ as well.  Thus, the lattice of submonoids given in Figure \ref{fig:lattice} simplifies a little in the case of infinite $\PB_X$.  Figure~\ref{fig:lattice2} pictures this simplified lattice.  All of the submonoids pictured in Figure \ref{fig:lattice2} are distinct, as may be deduced from the descriptions of these in Lemma \ref{lem:PBunits} and Theorems \ref{thm:IGPBX}, \ref{thm:FPBX} and \ref{thm:FLFR}.
\end{rem}

\begin{rem}
Although Corollary \ref{cor:PB=BB}(i) gives $\PB_X=\GRX \GLX$, it is not the case that  $\PB_X=\GLX \GRX$.  Rather, we have $\GLX\GRX=\set{\al\in\PB_X}{t(\al)=|X|}$.  Indeed, if $\al\in\PB_X$ satisfies $t(\al)=|X|$, then we may write $\al=\partn{a_x}{C_j}{b_x}{D_k}$, and it is then easy to see that $\al=\be\ga$, where $\be=\partn{a_x}{C_j}{x}{}\in\GLX$ and $\ga=\partn{x}{}{b_x}{D_k}\in\GRX$.  Conversely, if $\de\in\GLX$ and $\ve\in\GRX$, then from $\codom(\de)=X=\dom(\ve)$, we obtain $\dom(\de\ve)=\dom(\de)$, and so $t(\de\ve)=|X|$.  (As noted in the proof of \cite[Lemma 4.2]{HHR1998}, this also follows from considerations of Green's relations.)
\end{rem}

We are now ready to prove the first main result of this section.

\begin{thm}\label{thm:PBXSX}
Let $X$ be an infinite set.
\bit
\itemit{i} We have $\relrank{\PB_X}{\S_X}=2$.  
\itemit{ii} If $\al,\be\in\PB_X$, then $\PB_X=\la\S_X\cup\{\al,\be\}\ra$ if and only if (renaming if necessary) $\al\in\GRX$, $\be\in\GLX$ and $h^*(\al)=h(\be)=|X|$.
\eit
\end{thm}

\pf
If $\al\in\GRX$ and $\be\in\GLX$ are such that $h^*(\al)=h(\be)=|X|$, then Lemma \ref{lem:alSXbe} gives $\PB_X=\la\S_X\cup\{\al,\be\}\ra$.  This gives the backwards implication in (ii), and also shows that $\relrank{\PB_X}{\S_X}\leq2$; the reverse inequality follows from Lemma \ref{lem:rankMGE}(ii).  

It remains to show the forwards implication in (ii).  With this in mind, suppose $\al,\be\in\PB_X$ are such that $\PB_X=\la\S_X\cup\{\al,\be\}\ra$.  Renaming if necessary, Lemma \ref{lem:rankMGE}(i) gives $\al\in\GRX\sm\S_X$ and $\be\in\GLX\sm\S_X$.  Let $\ga\in\PB_X$ be such that $h^*(\ga)=|X|$, and consider an expression $\ga=\de_1\cdots\de_k$, where $\de_1,\ldots,\de_k\in\S_X\cup\{\al,\be\}$.  Corollary \ref{cor:tq}(ii) gives $h^*(\de_i)=|X|$ for some $i$.  Since $h^*(\be)=0$ (as $\be\in\GLX$) and $h^*(\ve)=0$ for all $\ve\in\S_X$, it follows that $\de_i=\al$, and so $h^*(\al)=|X|$.  A similar argument gives $h(\be)=|X|$.
\epf

\begin{rem}
Note that Theorem \ref{thm:PBXSX}(i) is true for $2\leq|X|<\aleph_0$ as well.  However, if $X$ is finite, then ${\PB_X=\la\S_X\cup\{\al,\be\}\ra}$ if and only if (renaming if necessary) $h(\al)=h^*(\al)=s(\be)=s^*(\be)=1$ and ${s(\al)=s^*(\al)=h(\be)=h^*(\be)=0}$.  This all follows from the proof of \cite[Proposition 3.16]{DEG2017}.
\end{rem}

Now that we have calculated $\relrank{\PB_X}{\S_X}$, it is easy to deduce the values of $\relrank{\PB_X}{\GLX}$ and $\relrank{\PB_X}{\GRX}$.  The next result only gives the statement for $\GLX$; the corresponding result for $\GRX$ is dual.

\begin{thm}\label{thm:PBXGLX}
Let $X$ be an infinite set.
\bit
\itemit{i} We have $\relrank{\PB_X}{\GLX}=1$.  
\itemit{ii} If $\al\in\PB_X$, then $\PB_X=\la \GLX\cup\{\al\}\ra$ if and only if $\al\in \GRX$ and $h^*(\al)=|X|$.
\eit
\end{thm}

\pf
If $\al\in \GRX$ is such that $h^*(\al)=|X|$, then for any $\be\in\GLX$ with $h(\be)=|X|$, Lemma \ref{lem:alSXbe} gives ${\PB_X=\la\S_X\cup\{\al,\be\}\ra\sub\la\GLX\cup\{\al\}\ra}$.  This gives the backwards implication in (ii), and ${\relrank{\PB_X}{\GLX}\leq1}$; the reverse inequality is obvious, since $\PB_X\not=\GLX$.

For the forwards implication in (ii), suppose $\al\in\PB_X$ is such that $\PB_X=\la \GLX\cup\{\al\}\ra$.  By Lemma~\ref{lem:rankMGE}(i), $\GLX\cup\{\al\}$ contains at least one element of $\GRX\sm\S_X$; since $(\GRX\sm\S_X)\cap\GLX=\emptyset$, it follows that this element must be $\al$, and so $\al\in \GRX$.  The proof of Theorem~\ref{thm:PBXSX}(ii) works virtually unmodified to show that $h^*(\al)=|X|$, noting that $h^*(\ve)=0$ for all $\ve\in \GLX$.
\epf

Next, we wish to calculate the relative ranks of $\GLX$ and $\GRX$ modulo $\S_X$.  In contrast to the previous situations (Theorems \ref{thm:PBXSX} and \ref{thm:PBXGLX}), we will see that $\relrank{\GLX}{\S_X}$ and $\relrank{\GRX}{\S_X}$ depend on the value of~$|X|$: more specifically, they depend on the number of infinite cardinals not exceeding $|X|$.  Again, we just treat the $\GLX$ case.

\begin{lemma}\label{lem:GLXSX1}
Let $X$ be an infinite set, and let $\Om=\set{\al_\mu,\be_\mu}{\mu=1\text{ or }\aleph_0\leq\mu\leq|X|}\sub\GLX$, where
\[
h(\al_\mu)=s(\be_\mu)=\mu \AND s(\al_\mu)=h(\be_\mu)=0.
\]
Then $\GLX=\la\S_X\cup\Om\ra$.
\end{lemma}

\pf
For $n\in\N\sm\{0\}$, define $\al_n=\al_1^n$ and $\be_n=\be_1^n$.  Then Lemma \ref{lem:shdef2}(i) and (iii) gives
\[
h(\al_n)=s(\be_n)=n \AND s(\al_n)=h(\be_n)=0.
\]
We also let $\al_0,\be_0$ be arbitrary elements of $\S_X$.  Now let $\ga\in\GLX$ be arbitrary, and write $\mu=h(\ga)$ and $\nu=s(\ga)$.  Then, again by Lemma \ref{lem:shdef2}(i) and (iii), $\al_\mu\be_\nu\in\la\S_X\cup\Om\ra$ satisfies
\[
h(\al_\mu\be_\nu)=h(\al_\mu)+h(\be_\nu)=\mu+0=\mu \ANDSIM s(\al_\mu\be_\nu)=\nu.
\]
Since also $\codom(\al_\mu\be_\nu)=X$, as $\al_\mu\be_\nu\in\GLX$, Lemma \ref{lem:action} gives $\ga\in\S_X\al_\mu\be_\nu\S_X\sub\la\S_X\cup\Om\ra$.
\epf

\begin{lemma}\label{lem:GLXSX2}
Let $X$ be an infinite set, and suppose $\Om\sub\GLX$ is such that $\GLX=\la\S_X\cup\Om\ra$.  Then $\Om$ contains a subset of the form described in Lemma \ref{lem:GLXSX1}.
\end{lemma}

\pf
Let $\mu$ be any cardinal such that either $\mu=1$ or $\aleph_0\leq\mu\leq|X|$.  We must show that there exist elements $\al,\be\in\Om$ such that 
\[
h(\al)=s(\be)=\mu \AND s(\al)=h(\be)=0.
\]
We just prove the existence of $\al$, as the argument for $\be$ is almost identical.
Let $\si\in\GLX$ be such that $h(\si)=\mu$ and $s(\si)=0$, and consider an expression $\si=\al_1\cdots\al_k$, where $\al_1,\ldots,\al_k\in\Om$.  Then Lemma~\ref{lem:shdef2}(i) gives
\[
\mu=h(\si)=h(\al_1\cdots\al_k)=h(\al_1)+\cdots+h(\al_k) \ANDSIM 0=s(\al_1)+\cdots+s(\al_k).
\]
The latter gives $s(\al_i)=0$ for all $i$, and the former gives $h(\al_i)=\mu$ for some $i$; we take $\al=\al_i$.
\epf

Here is the final main result of this section; it follows quickly from Lemmas \ref{lem:GLXSX1} and \ref{lem:GLXSX2}, after checking that the set $\Om$ from Lemma \ref{lem:GLXSX1} has the appropriate size.

\begin{thm}\label{thm:GLXSX}
Let $X$ be an infinite set, and let $\rho$ be the number of cardinals $\mu$ satisfying $\aleph_0\leq\mu\leq|X|$.
\bit
\itemit{i} We have $\relrank{\GLX}{\S_X}=2+2\rho$.  
\itemit{ii} If $\rho<\aleph_0$, and if $\Om\sub \GLX$ with $|\Om|=2+2\rho$, then $\GLX=\la\S_X\cup \Om\ra$ if and only if $\Om$ has the form described in Lemma \ref{lem:GLXSX1}. \epfres
\eit
\end{thm}

\begin{rem}\label{rem:GLXSX}
The assumption $\rho<\aleph_0$ is essential in Theorem \ref{thm:GLXSX}(ii); indeed, if $\rho\geq\aleph_0$, then $\Om\sub\GLX$ could contain a \emph{proper} subset of the form described in Lemma \ref{lem:GLXSX1}, yet still have $|\Om|=2+2\rho$.
\end{rem}

\begin{rem}
If we write $|X|=\aleph_\al$, where $\al$ is an ordinal, then $\rho=1+|\al|$.  Thus, 
\begin{align*}
\relrank{\GLX}{\S_X}<\aleph_0 &\iff \rho<\aleph_0 \iff |X|=\aleph_n \text{ for some } n\in\N.
\intertext{If $\rho\geq\aleph_0$, then $\relrank{\GLX}{\S_X}=\rho$.  Thus, writing $\om$ and $\om_1$ for the first countable and uncountable ordinals, respectively,}
\relrank{\GLX}{\S_X}=\aleph_0 &\iff \rho=\aleph_0 \iff \aleph_\om\leq|X|<\aleph_{\omega_1}.
\end{align*}
In particular, $\relrank{\GLX}{\S_X}$ is countable for uncountably many values of $|X|$.  Similar comments may be made for other relative ranks whose values involve the parameter~$\rho$; see Theorems \ref{thm:FXSX}, \ref{thm:FLXFX}, \ref{thm:FLXGLX} and \ref{thm:FLXSX}.
\end{rem}

\section{Idempotents}\label{sect:EX}

All other submonoids of $\PB_X$ we consider will include the set $E(\PB_X)$ of all idempotents among their generators.  
Accordingly, in this section, we investigate the submonoid $\E(\PB_X) = \la E(\PB_X)\ra$ generated by all such idempotents.  For simplicity, we will write $\EX$ for $\E(\PB_X)$ from this point on.  The main results of this section include a characterisation of the elements of $\EX$ in Theorem~\ref{thm:IGPBX}, and the calculation of the relative rank of $\PB_X$ modulo $\EX$ (equivalently, modulo $E(\PB_X)$) in Theorem~\ref{thm:PBXEX}, where we also classify the minimal generating sets modulo $\EX$.  

The idempotents of $\PB_X$ were described (and enumerated) in \cite{DEEFHHL1}; however, we do not need the full classification here.  Rather, we just need to know that certain simple Brauer graphs are idempotents.  The next result follows from \cite[Theorem 5]{DEEFHHL1}, but we include a simple proof for convenience; we will often use this result without explicit reference.  

\begin{lemma}\label{lem:idempotents}
If $X$ is an arbitrary set, and if $\al\in\PB_X$ is such that $x\al=x$ for all $x\in\dom(\al)$, then $\al$ is an idempotent.
\end{lemma}

\pf
For any $\al\in\PB_X$, all of the non-transversals of $\al$ remain in the product $\al^2$.  The stated assumption ensures that this is the case for the transversals of $\al$ as well.
\epf

To describe $\EX$, we must first define some more parameters associated to Brauer graphs.  We define the \emph{fix}, \emph{support} and \emph{shift} sets and cardinals of $\al\in\PB_X$ to be
\begin{align*}
\Fix(\al)&=\set{x\in\dom(\al)}{x\al=x}, & \Supp(\al)&=X\sm\Fix(\al), & \Sh(\al)&=\dom(\al)\sm\Fix(\al), \\
\fix(\al)&=|\Fix(\al)|, & \supp(\al)&=|\Supp(\al)|, & \sh(\al)&=|\Sh(\al)|.
\end{align*}
Note that the condition ``$x\al=x$ for all $x\in\dom(\al)$'' in Lemma \ref{lem:idempotents} could be restated as ``${\Fix(\al)=\dom(\al)}$'' or, equivalently, ``$\sh(\al)=0$''.  It is easy to construct idempotents of $\PB_X$ where these conditions do not hold.  Note also that $\Supp(\al)=\Defect(\al)\sqcup\Sh(\al)$.

Two important steps in the proof of Theorem \ref{thm:IGPBX} (which describes the elements of infinite $\EX$) have been completed elsewhere in the literature.  Namely, the monoid $\EX$ was described in the case of \emph{finite}~$X$ in \cite{DEG2017}, and the idempotent-generated subsemigroup of the larger \emph{partition monoid} $\P_X$ was described in \cite{EF2012}.  We will postpone a discussion of the latter (see Lemma \ref{lem:EF} and the preceding paragraphs).  The next result is part of \cite[Theorem~3.18]{DEG2017}.  

\begin{thm}\label{thm:DEG}
If $X$ is a finite set, then
\[
\epfreseq
\EX = \set{\al\in\PB_X}{\defect(\al)\leq1\text{ and }\sh(\al)=0} \cup \set{\al\in\PB_X}{\defect(\al)\geq2}.
\]
\end{thm}

A key role in the proof of Theorem \ref{thm:IGPBX} is played by another important submonoid, which has been useful in a number of other contexts \cite{EMRT2018,EF2012,JEgrpm,East2011}.  By parts (ii) and (iii) of Lemma~\ref{lem:shdef}, the set
\[
\I_X = \set{\al\in\PB_X}{h(\al)=h^*(\al)=0}
\]
is a submonoid of $\PB_X$.  It was noted in \cite[Section 2]{EF2012} that $\I_X$ is isomorphic to the symmetric inverse monoid on the set $X$: i.e., the set of all injective partial transformations of~$X$ under the operation of relational composition.
Note that $\I_X$ is closed under the $\al\mt\al^*$ map discussed in Section \ref{sect:PBX}.  Indeed, if~$\al\in\I_X$, then $\al^*\equiv\al^{-1}$ is the inverse mapping of $\al$.  
The main remaining step in establishing Theorem~\ref{thm:IGPBX} is to describe the elements of $\I_X$ that are products of idempotents from $\PB_X$; this is accomplished in Lemma~\ref{lem:IX}, the proof of which requires the next three preliminary lemmas.

\begin{lemma}\label{lem:mcycle}
Let $W$ be a finite set of size $3m$, where $m\geq2$, let $x_1,\ldots,x_m$ be distinct elements of $W$, and let $\al=\partnII{x_1}\cdots{x_{m-1}}{x_m}{x_2}\cdots{x_m}{x_1}\in\I_W$.  Then $\al=\be\ga\de$ for some $\be,\ga,\de\in E(\PB_W)$.
\end{lemma}

\begin{figure}[ht]
\begin{center}
\begin{tikzpicture}[scale=0.5]
\begin{scope}
\node[above] at (1,2) {\tiny $x_1$};
\node[above] at (4,2) {\tiny $x_2$};
\node[above] at (9,2) {\tiny $x_{m-1}$};
\node[above] at (12,2) {\tiny $x_m$};
\stline14
\stline47
\stline{9}{12}
\stline{12}1
\foreach \x in {1,4,9,12} {\uv{\x}} ;
\foreach \x in {2,3,10,11,13,14} {\uvg{\x}} ;
\foreach \x in {1,4,7,12} {\lv{\x}} ;
\foreach \x in {2,3,5,6,13,14} {\lvg{\x}} ;
\draw[dotted] (5,2)--(8,2);
\draw[dotted] (8,0)--(11,0);
\node at (16.5,1) {$=$};
\draw[|-|] (0,2)--(0,0); \node[left] at (0,1) {$\al$};
\end{scope}
\begin{scope}[shift={(18,2)}]
\node[above] at (1,2) {\tiny $x_1$};
\node[above] at (4,2) {\tiny $x_2$};
\node[above] at (9,2) {\tiny $x_{m-1}$};
\node[above] at (12,2) {\tiny $x_m$};
\stline11
\stline44
\stline{9}{9}
\stline{12}{12}
\darc23
\darc56
\darc{10}{11}
\darc{13}{14}
\foreach \x in {1,4,9,12} {\uv{\x}} ;
\foreach \x in {2,3,10,11,13,14} {\uvg{\x}} ;
\foreach \x in {1,4,12} {\lv{\x}} ;
\foreach \x in {2,3,5,6,13,14} {\lvg{\x}} ;
\draw[dotted] (5,2)--(8,2);
\draw[dotted] (7,0)--(8,0);
\draw[|-] (15,2)--(15,0); \node[right] at (15,1) {$\be$};
\end{scope}
\begin{scope}[shift={(18,0)}]
\stline33
\stline66
\stline{11}{11}
\stline{14}{14}
\uarc12
\uarc45
\uarc{9}{10}
\uarc{12}{13}
\darc24
\darc57
\darc{10}{12}
\darcx1{13}{.9}
\foreach \x in {1,4,9,12} {\uv{\x}} ;
\foreach \x in {2,3,5,6,10,11,13,14} {\uvg{\x}} ;
\foreach \x in {1,4,7,12} {\lv{\x}} ;
\foreach \x in {2,3,5,6,13,14} {\lvg{\x}} ;
\draw[|-|] (15,2)--(15,0); \node[right] at (15,1) {$\ga$};
\end{scope}
\begin{scope}[shift={(18,-2)}]
\stline11
\stline44
\stline77
\stline{12}{12}
\uarc23
\uarc56
\uarc{10}{11}
\uarc{13}{14}
\foreach \x in {1,4,7,12} {\uv{\x}} ;
\foreach \x in {2,3,5,6,10,11,13,14} {\uvg{\x}} ;
\foreach \x in {1,4,7,12} {\lv{\x}} ;
\foreach \x in {2,3,5,6,13,14} {\lvg{\x}} ;
\draw[dotted] (8,2)--(9,2);
\draw[dotted] (8,0)--(11,0);
\draw[-|] (15,2)--(15,0); \node[right] at (15,1) {$\de$};
\end{scope}
\end{tikzpicture}
\end{center}
\vspace{-5mm}
\caption{Verification of the equation $\al=\be\ga\de$ from the proof of Lemma \ref{lem:mcycle}; see the text for more details.}
\label{fig:mcycle}
\end{figure}
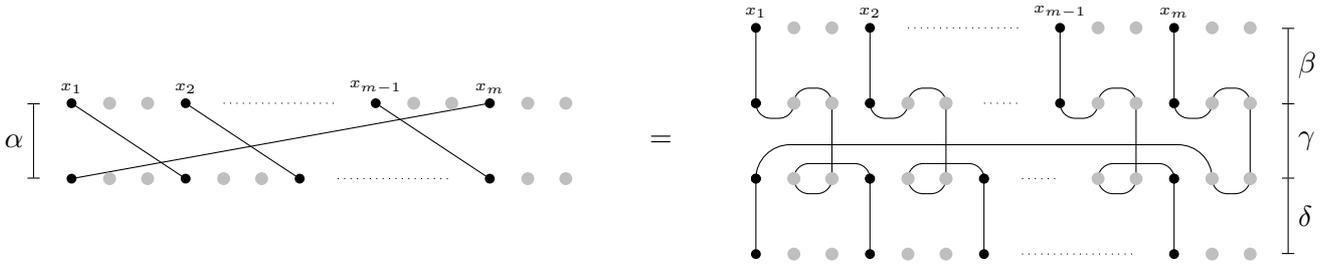

\pf
In Figure \ref{fig:mcycle}, we define the idempotents $\be,\ga,\de\in E(\PB_W)$ and show that $\al=\be\ga\de$.  In the figure, the~$2m$ elements of $W\sm\{x_1,\ldots,x_m\}$ are shaded gray. 
\epf

\begin{lemma}\label{lem:mtrail}
Let $W$ be a finite set of size $3m-2$, where $m\geq2$, let $x_1,\ldots,x_m$ be distinct elements of $W$, and let $\al=\partnIII{x_1}\cdots{x_{m-1}}{x_2}\cdots{x_m}\in\I_W$.  Then $\al=\be\ga\de$ for some $\be,\ga,\de\in E(\PB_W)$.
\end{lemma}

\pf
The proof is almost identical to that of Lemma \ref{lem:mcycle}.  In fact, Figure \ref{fig:mcycle} may easily be modified to work here as well.  We simply remove the last two gray vertices from each row as well as any blocks from~$\be$,~$\ga$ and $\de$ that involve any of these vertices, and also the transversals $\{x_m,x_1'\}$, $\{x_m,x_m'\}$ and $\{x_1,x_1'\}$ from~$\al$,~$\be$ and~$\de$, respectively.
\epf

The proof of Lemma \ref{lem:mcycle} may also be easily modified to prove the following.

\begin{lemma}\label{lem:infinitecycletrail}
Let $W$ be a countably infinite set, let $\{\ldots,x_1,x_2,x_3,\ldots\}$ be a subset of $W$ with infinite complement, and let $\al$ be any of the three elements $\partnII{x_1}{x_2}{x_3}\cdots{x_2}{x_3}{x_4}\cdots$, $\partnII{x_2}{x_3}{x_4}\cdots{x_1}{x_2}{x_3}\cdots$ or $\partnIV{x_1}{x_2}{x_3}{x_2}{x_3}{x_4}$ of $\I_W$.  Then $\al=\be\ga\de$ for some $\be,\ga,\de\in E(\PB_W)$. \epfres
\end{lemma}

If $\al\in\PB_X$, and if $W\sub X$ is such that any edge $\{x,y\}$ of $\al$ satisfies either $x,y\in W\cup W'$ or $x,y\in (X\sm W)\cup(X\sm W)'$, then we define the \emph{restriction} of $\al$ to $W$ to be the induced subgraph of $\al$ on vertex set $W\cup W'$; note that this restriction belongs to $\PB_W$.

If $\set{W_i}{i\in I}$ is some collection of pairwise disjoint sets, and if $\al_i\in\PB_{W_i}$ for all $i$, then we denote by $\bigcup_{i\in I}\al_i$ the Brauer graph with vertex set $\bigcup_{i\in I}W_i$ and with edge set equal to the union of the edge sets of the $\al_i$. 
Sometimes this operation is denoted $\oplus$ or $\otimes$ (see for example \cite{DEEFHHL1,LZ2015}), but since we view the elements of partial Brauer monoids as graphs, $\cup$ seems more appropriate for our purposes.

The proof of the next lemma uses \emph{cycle-trail} notation for elements of $\I_X$, which we now describe.
\bit
\itemnit{i} A \emph{finite cycle} is a permutation $(x_1,\ldots,x_m)$ of a set $\{x_1,\ldots,x_m\}$ that maps $x_1\mt x_2\mt\cdots\mt x_m\mt x_1$.
\itemnit{ii} An \emph{infinite cycle} is a permutation $(\ldots,x_{-1},x_0,x_1,x_2,\ldots)$ of a set $\{\ldots,x_{-1},x_0,x_1,x_2,\ldots\}$ that maps $\cdots\mt x_{-1}\mt x_0\mt x_1\mt x_2\mt\cdots$.
\itemnit{iii} A \emph{finite trail} is a partial bijection $[x_1,\ldots,x_m]$ of a set $\{x_1,\ldots,x_m\}$ that maps $x_1\mt x_2\mt\cdots\mt x_m$.  This trail has domain $\{x_1,\ldots,x_{m-1}\}$ and codomain $\{x_2,\ldots,x_m\}$.
\itemnit{iv} A \emph{right-infinite} trail is a partial bijection $[x_1,x_2,\ldots]$ of a set $\{x_1,x_2,\ldots\}$ that maps ${x_1\mt x_2\mt\cdots}$.  This trail has domain $\{x_1,x_2,\ldots\}$ and codomain $\{x_2,x_3,\ldots\}$.
\itemnit{v} A \emph{left-infinite} trail is a partial bijection $[\ldots,x_2,x_1]$ of a set $\{x_1,x_2,\ldots\}$ that maps ${\cdots\mt x_2\mt x_1}$.  This trail has domain $\{x_2,x_3,\ldots\}$ and codomain $\{x_1,x_2,\ldots\}$.
\eit
The cycle in (i) is called an \emph{$m$-cycle}, and the trail in (iii) an \emph{$m$-trail}; these are called \emph{trivial} if $m=1$, or \emph{non-trivial} if $m\geq2$.  Note that a trivial cycle is the identity map on a one-element set, while a trivial trail is the empty map on a one-element set.  Cycles and trails may be regarded as elements of suitable partial Brauer monoids, using the identification of $\I_X$ with a submonoid of $\PB_X$ described above.  It is easy to see that any element of $\I_X$ may be uniquely decomposed as a (disjoint) union of cycles and trails.

\begin{lemma}\label{lem:IX}
If $X$ is an infinite set, and if $\al\in\I_X$ is such that $\defect(\al)=\codef(\al)\geq\max(\aleph_0,\sh(\al))$, then~${\al\in\EX}$.
\end{lemma}

\pf
For the proof, we define the \emph{fail} set and cardinal of $\al$ by
\[
\Fail(\al)=X\sm(\dom(\al)\cup\codom(\al))=\Defect(\al)\cap\Codef(\al) \AND \fail(\al)=|{\Fail(\al)}|.
\]
We consider two cases, according to whether $\sh(\al)\leq\fail(\al)$ or $\sh(\al)>\fail(\al)$.

\pfcase{1}  Suppose first that $\sh(\al)\leq\fail(\al)$.  We first claim that $\fail(\al)\geq\aleph_0$.  To prove this, suppose to the contrary that $\fail(\al)<\aleph_0$.  Then
\[
\codef(\al) = |{\Codef(\al)}| = |{\Fail(\al)\sqcup(\Codef(\al)\sm\Fail(\al))}| = \fail(\al) +|{\Codef(\al)\sm\Fail(\al)}| .
\]
Since $\codef(\al)\geq\aleph_0$ and $\fail(\al)<\aleph_0$, it follows that $|{\Codef(\al)\sm\Fail(\al)}|=\codef(\al)\geq\aleph_0$.  Now, 
\[
\Codef(\al)\sm\Fail(\al) = \big(X\sm\codom(\al)\big) \sm \big(X\sm(\dom(\al)\cup\codom(\al))\big) = \dom(\al)\sm\codom(\al) \sub \Sh(\al).
\]
It follows that $\aleph_0\leq|{\Codef(\al)\sm\Fail(\al)}| \leq \sh(\al) \leq \fail(\al) <\aleph_0$, a contradiction.  This completes the proof of the claim that $\fail(\al)\geq\aleph_0$.  Consequently, and using $\sh(\al)\leq\fail(\al)$, we may fix two subsets $Y,Z\sub\Fail(\al)$ such that $Y\cap Z=\emptyset$ and $|Y|=|Z|=\sh(\al)$.  We also fix bijections $\phi:\Sh(\al)\to Y$ and $\psi:\Sh(\al)\to Z$.

Now suppose $\al$ has
\bit\bmc2
\itemnit{i} non-trivial finite cycles $\set{\al_i}{i\in I}$,
\itemnit{ii} infinite cycles $\set{\al_j}{j\in J}$,
\item[] ~
\itemnit{iii} non-trivial finite trails $\set{\al_k}{k\in K}$,
\itemnit{iv} right-infinite trails $\set{\al_l}{l\in L}$,
\itemnit{v} left-infinite trails $\set{\al_m}{m\in M}$,
\emc\eit
where the indexing sets $I,J,K,L,M$ are assumed to be pairwise disjoint.  Note that we have not listed the trivial cycles and trails.  Write $Q=I\cup J\cup K\cup L\cup M$.  Note that $\Sh(\al) = \bigsqcup_{q\in Q}\dom(\al_q)$.  For each $q\in Q$, let $X_q=\dom(\al_q)\cup\codom(\al_q)$, and put
\[
W_q=X_q\cup \dom(\al_q)\phi\cup \dom(\al_q)\psi.
\]
If $i\in I$ and $\al_i$ is an $m$-cycle, then $|W_i|=3m$; if $k\in K$ and $\al_k$ is an $m$-trail, then $|W_k|=3m-2$; and if $q\in J\cup L\cup M$, then $|X_q|=|W_q|=|W_q\sm X_q|=\aleph_0$.  For each $q\in Q$, let $\be_q$ be the unique element of $\I_{W_q}$ with the same transversals as $\al_q$.  Then by Lemma \ref{lem:mcycle}, \ref{lem:mtrail} or \ref{lem:infinitecycletrail}, as appropriate, there exist $\ga_q,\de_q,\ve_q\in E(\PB_{W_i})$ such that $\be_q=\ga_q\de_q\ve_q$.  
Note that $\Fix(\al) \sub X\sm\bigcup_{q\in Q}W_q \sub \Fix(\al)\cup\Fail(\al)$.  
Put $V=X\sm\bigcup_{q\in Q}W_q$, and let $\zeta\in\PB_V$ be the Brauer graph with edge set $\bigset{\{x,x'\}}{x\in\Fix(\al)}$.  With $\xi$ denoting any of $\ga,\de,\ve$, we define $\xi = \zeta\cup\bigcup_{q\in Q} \xi_q$.  Then by construction, $\ga,\de,\ve\in E(\PB_X)$ and $\al=\ga\de\ve$.

\pfcase{2}  Now suppose $\sh(\al)>\fail(\al)$.  For simplicity, write $A=\dom(\al)$ and $B=\codom(\al)$, and put $B_1=(A\cap B)\al$, noting that $B\sm B_1=(A\sm B)\al$.  In the proof of \cite[Lemma 27]{EF2012}, it was shown that $|A\sm B|=|B\sm A|=|B\sm B_1|$.  Fix a bijection $\phi:A\sm B\to B\sm A$, and define $\be,\ga\in\I_X$ by
\[
\be = \partnV aba{b\phi}_{a\in A\cap B,\ b\in A\sm B} \AND \ga = \partnV a{b\phi}{a\al}{b\al}_{a\in A\cap B,\ b\in A\sm B}.
\]
Evidently, we have $\al=\be\ga$, so the proof will be complete if we can show that $\be,\ga\in\EX$.  Now, $\ga$ maps~$B$ bijectively onto itself, so $\defect(\ga)=\codef(\ga)=\fail(\ga)=|X\sm B|=\codef(\al)\geq\aleph_0$.  
Because of Case~1, we will be able to conclude that $\ga\in\EX$ if we can show that $\defect(\ga)\geq\sh(\ga)$.  
Define $C=\set{x\in A\cap B}{x\al\not=x}$, and note that $\Sh(\al)=(A\sm B)\cup C$.  But also $\Sh(\ga)\sub(B\sm A)\cup C$, and so
\[
\sh(\ga)\leq|B\sm A|+|C|=|A\sm B|+|C|=\sh(\al)\leq\codef(\al)=\defect(\ga).
\]
As noted above, this completes the proof that $\ga\in\EX$.  It remains to show that $\be\in\EX$.

Now, $\be$ maps $A$ bijectively onto $B$, mapping $A\cap B$ identically onto itself, and $A\sm B$ onto $B\sm A$.  Since $A\sm B$ and $B\sm A$ are disjoint, it follows that the cycle-trail decomposition of $\be$ consists of 1-cycles (one for each element of $A\cap B$), 1-trails (one for each element of $\Fail(\be)=X\sm(A\cup B)$), and 2-trails (one for each element of $A\sm B$).    Note also that $\sh(\be)=|A\sm B|=|B\sm A|$.  We must consider two subcases.

\pfcase{2.1}  First consider the case in which $\sh(\be)<\aleph_0$.  Now, 
\[
\defect(\be) = |X\sm A| = \defect(\al) \AND \codef(\be) = |X\sm B| = \codef(\al).
\]
Thus, $\defect(\be)=\codef(\be)\geq\aleph_0=\max(\aleph_0,\sh(\be))$.  Also, 
\[
\aleph_0 \leq \defect(\be) = |X\sm A| = |B\sm A| + |X\sm(A\cup B)| = \sh(\be) + \fail(\be).
\]
From $\sh(\be)<\aleph_0$, it then follows that $\fail(\be)\geq\aleph_0$.  Consequently, $\sh(\be)<\aleph_0\leq\fail(\be)$, and so $\be\in\EX$, by Case 1.

\begin{figure}[ht]
\begin{center}
\begin{tikzpicture}[scale=0.5]
\begin{scope}
\foreach \x in {1,2,3,4} {\uv{\x}\lv{\x}} ;
\node[above] at (1,2) {\tiny $a_h$};
\node[above] at (2,2) {\tiny $b_h$};
\node[above] at (3,2) {\tiny $c_h$};
\node[above] at (4,2) {\tiny $d_h$};
\stline13
\stline24
\node at (6.5,1) {$=$};
\draw[|-|] (0,2)--(0,0); \node[left] at (0,1) {$\be_h$};
\end{scope}
\begin{scope}[shift={(8,2)}]
\foreach \x in {1,2,3,4} {\uv{\x}\lv{\x}} ;
\node[above] at (1,2) {\tiny $a_h$};
\node[above] at (2,2) {\tiny $b_h$};
\node[above] at (3,2) {\tiny $c_h$};
\node[above] at (4,2) {\tiny $d_h$};
\stline11
\stline22
\darc34
\draw[|-] (5,2)--(5,0); \node[right] at (5,1) {$\eta_h$};
\end{scope}
\begin{scope}[shift={(8,0)}]
\foreach \x in {1,2,3,4} {\uv{\x}\lv{\x}} ;
\stline11
\stline44
\uarc23
\darc23
\draw[|-|] (5,2)--(5,0); \node[right] at (5,1) {$\si_h$};
\end{scope}
\begin{scope}[shift={(8,-2)}]
\foreach \x in {1,2,3,4} {\uv{\x}\lv{\x}} ;
\stline33
\stline44
\uarc12
\draw[-|] (5,2)--(5,0); \node[right] at (5,1) {$\pi_h$};
\end{scope}
\end{tikzpicture}
\end{center}
\vspace{-5mm}
\caption{Verification of the equation $\be_h=\eta_h\si_h\pi_h$ from the proof of Lemma \ref{lem:IX}; see the text for more details.}
\label{fig:disjoint}
\end{figure}
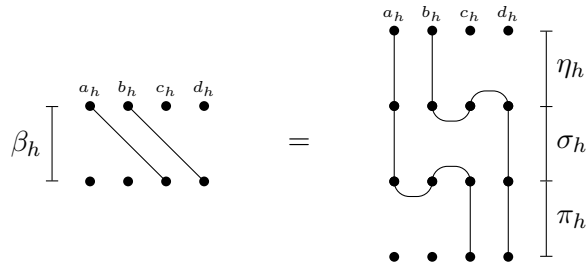

\pfcase{2.2}  Finally, suppose $\sh(\be)\geq\aleph_0$.  Choose some indexing set $H$ with $|H|=\sh(\be)=|A\sm B|$.  Since $|H|\geq\aleph_0$, we may write
\[
A\sm B=\set{a_h}{h\in H} \sqcup \set{b_h}{h\in H}.
\]
For $h\in H$, put $c_h=a_h\be$ and $d_h=b_h\be$.  Then the $2$-trails of $\be$ are $\bigset{[a_h,c_h],[b_h,d_h]}{h\in H}$.  For each $h\in H$, let $W_h=\{a_h,b_h,c_h,d_h\}$, and let $\be_h=[a_h,c_h]\cup[b_h,d_h]\in\PB_{W_h}$ be the restriction of $\be$ to $W_h$.  In Figure~\ref{fig:disjoint}, we show that $\be_h=\eta_h\si_h\pi_h$ for some $\eta_h,\si_h,\pi_h\in E(\PB_{W_h})$ for each $h$.  
Put $V=(A\cap B)\cup(X\sm(A\cup B))$, and let $\zeta\in\PB_V$ be the Brauer graph with edge set $\bigset{\{x,x'\}}{x\in A\cap B}$.
With $\xi$ denoting any of $\eta,\si,\pi$, we define $\xi = \zeta\cup\bigcup_{h\in H} \xi_h$.  Then $\eta,\si,\pi\in E(\PB_X)$ and $\be=\eta\si\pi$, completing the proof.
\epf

We have already noted that $\PB_X$ is a submonoid of the larger \emph{partition monoid} $\P_X$.  The idempotent-generated subsemigroup of $\P_X$ was described in \cite[Theorem 30]{EF2012}.  We do not need to give the full details of this result, or even fully define $\P_X$ itself, but we will make some comments that are relevant to the current situation.  To an element~$\al$ of $\P_X$, one may associate the the \emph{singularity} and \emph{cosingularity} parameters, denoted $\sing(\al)$ and $\cosing(\al)$, respectively; see \cite[p.~115]{EF2012}.  Of crucial importance here is that when $\al\in\PB_X$, we have $\sing(\al)=\defect(\al)$ and $\cosing(\al)=\codef(\al)$ in our current terminology.  There is also a notion of the \emph{shift} and \emph{support}, $\sh(\al)$ and $\supp(\al)$, of an element $\al$ of $\P_X$, and these coincide with our current definitions in the case that $\al\in\PB_X$; but note that $\supp(\al)$ was denoted $\operatorname{warp}(\al)$ in \cite{EF2012}.

Recall from \cite{NS1978,EF2012} that an element $\al\in\PB_X$ is a \emph{projection} if $\al^2=\al=\al^*$.  It is easy to see, using~\eqref{eq:al*}, that for any $\al\in\PB_X$, both $\al\al^*$ and $\al^*\al$ are projections.  Hence, if $\al$ is any idempotent from $\PB_X$, then again using \eqref{eq:al*}, $\al=\al\al^*\al=\al(\al\al)^*\al=(\al\al^*)(\al^*\al)$ is the product of two projections.  Thus, $\EX$ is also equal to the subsemigroup of $\PB_X$ generated by all projections.
The next result follows from the first paragraph of the proof of \cite[Theorem 30]{EF2012}.

\begin{lemma}\label{lem:EF}
Let $X$ be an infinite set, let $\be_1,\ldots,\be_k\in\PB_X$ be projections, and put $\al=\be_1\cdots\be_k$.  If $\supp(\be_i)\geq\aleph_0$ for some $i$, then $\defect(\al)=\codef(\al)\geq\max(\aleph_0,\sh(\al))$. \epfres
\end{lemma}

We are now ready to state and prove the first main result of this section, which characterises the elements of the idempotent-generated subsemigroup $\EX=\E(\PB_X)$ of $\PB_X$.

\begin{thm}\label{thm:IGPBX}
If $X$ is an infinite set, then
\begin{align*}
\EX = \set{\al\in\PB_X}{\defect(\al)\leq1 \text{ and } \sh(\al)=0} 
&\cup \set{\al\in\PB_X}{\defect(\al)\geq2 \text{ and } \supp(\al)<\aleph_0} \\
&\cup \set{\al\in\PB_X}{\defect(\al)=\codef(\al)\geq\max(\aleph_0,\sh(\al))}.
\end{align*}
\end{thm}

\pf
During the proof, we will write
\bit
\item $\Om_1 = \set{\al\in\PB_X}{\defect(\al)\leq1 \text{ and } \sh(\al)=0}$, 
\item ${\Om_2 = \set{\al\in\PB_X}{\defect(\al)\geq2 \text{ and } \supp(\al)<\aleph_0}}$,
\item $\Om_3 = \set{\al\in\PB_X}{\defect(\al)=\codef(\al)\geq\max(\aleph_0,\sh(\al))}$.
\eit
First suppose $\al\in\EX$.  As discussed above, we may write $\al=\be_1\cdots\be_k$, where $\be_1,\ldots,\be_k\in\PB_X$ are projections.  
If $\supp(\be_i)\geq\aleph_0$ for some $i$, then $\al\in\Om_3$, by Lemma \ref{lem:EF}.
Next, suppose $\supp(\be_i)<\aleph_0$ for all~$i$.  Put $W_i=\Supp(\be_i)$ for each $i$, and let $W=\bigcup_{i=1}^k W_i$, noting that $|W|<\aleph_0$ and $\Supp(\al)\sub W$.  For each~$i$, let $\ga_i\in E(\PB_W)$ be the restriction of $\be_i$ to $W$. Then $\ga_1\cdots\ga_k\in\E(\PB_W)$, and it quickly follows from Theorem \ref{thm:DEG} that $\al=\be_1\cdots\be_k\in\Om_1\cup\Om_2$.  This completes the proof that $\EX\sub\Om_1\cup\Om_2\cup\Om_3$.

To prove the reverse inclusion, first note that $\Om_1\cup\Om_2\sub\EX$ also follows quickly from Theorem \ref{thm:DEG}.  Now suppose $\al\in\Om_3$, and write $\al=\partn{a_i}{C_j}{b_i}{D_k}$.  Then $\al=\be\ga\de$, where $\be=\partn {a_i}{C_j}{a_i}{}$, $\ga=\partnVI{a_i}{b_i}$ and ${\de=\partn{b_i}{}{b_i}{D_k}}$.  By Lemma \ref{lem:idempotents}, $\be,\de\in E(\PB_X)$.  Also, $\ga\in\I_X$ satisfies $\defect(\ga)=\defect(\al)$, $\codef(\ga)=\codef(\al)$ and ${\sh(\ga)=\sh(\al)}$, so that $\defect(\ga)=\codef(\ga)\geq\max(\aleph_0,\sh(\ga))$.  Lemma \ref{lem:IX} then gives $\ga\in\EX$, and the proof is complete.
\epf

\begin{rem}\label{rem:EXinFX}
We note for later reference that any element $\al$ of $\EX$ satisfies $\defect(\al)=\codef(\al)$.  Indeed, this is obvious if $\al\in\Om_3$ (in the notation of the above proof), and follows quickly from Lemma \ref{lem:finitedefect} if $\al\in\Om_1\cup\Om_2$.
\end{rem}

\begin{rem}
Even though $\PB_X$ is a submonoid of $\P_X$, the idempotent-generated subsemigroup $\E(\PB_X)$ is not simply the intersection of $\E(\P_X)$ with $\PB_X$.  Indeed, any $\al\in\PB_X$ with $\defect(\al)=1$ and $\supp(\al)<\aleph_0$ is a product of idempotents from $\P_X$, as follows from \cite[Theorem 30]{EF2012} or \cite[Proposition 16]{East2011}; however, such an $\al$ is only a product of idempotents from $\PB_X$ if $\al$ is itself an idempotent, as follows from Theorem \ref{thm:IGPBX}.
\end{rem}

Now that we have characterised the elements of $\EX=\E(\PB_X)$, we wish to calculate the relative rank of~$\PB_X$ modulo $\EX$.  In Lemma \ref{lem:alSXbe} above, we proved that $\PB_X=\al\S_X\be$ for suitably chosen one-sided units $\al\in \GRX$ and $\be\in \GLX$.  The next lemma gives the analogous result for $\EX$; instead of requiring that ${h^*(\al)=h(\be)=|X|}$, we make the weaker assumption that ${\codef(\al)=\defect(\be)=|X|}$.

\begin{lemma}\label{lem:alEXbe}
Let $X$ be an infinite set, and let $\al\in\GRX$ and $\be\in\GLX$ with ${\codef(\al)=\defect(\be)=|X|}$.  Then $\PB_X=\al\EX\be$. 
\end{lemma}

\pf
Choose any partitions $\set{A_x}{x\in X}$ and $\set{B_x}{x\in X}$ of $\Codef(\al)$ and $\Defect(\be)$, respectively, with $|A_x|=|B_x|=2$ for all $x\in X$.  Define (full) Brauer graphs $\al_1 = \partn {x\al}{A_x}{x\al}{A_x}$, $\be_1 = \partn {x\be^{-1}}{B_x}{x\be^{-1}}{B_x}$, $\al_2=\al\al_1$ and $\be_2=\be_1\be$.  Then 
\[
\al_1,\be_1\in E(\PB_X) \COMMA \al_2=\partn {x}{}{x\al}{A_x}\in \GRX \COMMA \be_2=\partn {x\be^{-1}}{B_x}{x}{}\in \GLX \COMMA h^*(\al_2)=h(\be)=|X|.
\]
Now let $\ga\in\PB_X$ be arbitrary.  We follow steps (i)--(v) in the proof of Lemma \ref{lem:alSXbe} to define an element ${\de\in\PB_X}$ such that $\ga=\al_2\de\be_2$; however, after step (v), we instead define all the elements of ${(X\sm X_1)\cup(X\sm X_2)'}$ to be singletons of $\de$.  Since $\defect(\de)=\codef(\de)=|X|$, Theorem \ref{thm:IGPBX} gives $\de\in\EX$.  But then ${\ga=\al_2\de\be_2 = \al (\al_1\de\be_1)\be \in\al\EX\be}$, as required.
\epf

We may now calculate the relative rank of $\PB_X$ modulo $\EX=\E(\PB_X)$.

\begin{thm}\label{thm:PBXEX}
Let $X$ be an infinite set.
\bit
\itemit{i} We have $\relrank{\PB_X}{\EX}=2$.  
\itemit{ii} If $\al,\be\in\PB_X$, then $\PB_X=\la \EX\cup\{\al,\be\}\ra$ if and only if (renaming if necessary) $\al\in\GRX$, $\be\in\GLX$ and $\codef(\al)=\defect(\be)=|X|$.
\eit
\end{thm}

\pf Lemma \ref{lem:alEXbe} gives the backwards implication in (ii), and also $\relrank{\PB_X}{\EX}\leq2$; the reverse inequality follows from Lemma \ref{lem:rankMGE}(iii).  

For the forwards implication in (ii), suppose $\al,\be\in\PB_X$ are such that $\PB_X=\la\EX\cup\{\al,\be\}\ra$.  By Lemma~\ref{lem:rankMGE}(i), we may assume without loss of generality that $\al\in \GRX$ and $\be\in \GLX$.  We must show that $\codef(\al)=\defect(\be)=|X|$.  By duality, it suffices to prove the statement concerning $\al$.  To do so, let $\ga\in\PB_X$ be such that $\defect(\ga)<|X|=\codef(\ga)$, and consider an expression $\ga=\de_1\cdots\de_k$, where all of the factors belong to $\EX\cup\{\al,\be\}$.  For $0\leq i\leq k$, let $\ga_i=\de_1\cdots\de_i$.  Then for any $i$, Lemma \ref{lem:shdef}(iv) gives $\defect(\ga_i)\leq\defect(\ga_i\de_{i+1}\cdots\de_k)=\defect(\ga)<|X|$.  Since $\codef(\ga_0)=0$ and $\codef(\ga_k)=\codef(\ga)=|X|$, we may define $j=\min\set{i}{\codef(\ga_i)=|X|}$, noting that $1\leq j\leq k$ and $\codef(\ga_{j-1})<|X|$.  If also $\codef(\de_j)<|X|$, then Lemma \ref{lem:shdef}(iv) would give $|X|=\codef(\ga_j)=\codef(\ga_{j-1}\de_j)\leq\codef(\ga_{j-1})+\codef(\de_j)<|X|$, a contradiction.  So we must have $\codef(\de_j)=|X|$; thus, the proof will be complete if we can show that $\al=\de_j$.  Now, if also $\defect(\de_j)=|X|$, then since $|X|\geq\aleph_0$ and $|X|>\codef(\ga_{j-1})$, Lemma \ref{lem:shdef3}(iii) would give ${\defect(\ga_{j-1}\de_j)\geq\defect(\de_j)=|X|}$, contradicting $\defect(\ga_{j-1}\de_j)=\defect(\ga_j)<|X|$.  So we must in fact have ${\defect(\de_j)<|X|=\codef(\de_j)}$.  As noted in Remark \ref{rem:EXinFX}, we $\defect(\ve)=\codef(\ve)$ for all $\ve\in\EX$, so it follows that $\de_j\in\{\al,\be\}$.  Also, since $\be\in\GLX$, we have $\codef(\be)=0\not=|X|=\codef(\de_j)$, and so $\de_j=\al$.
\epf

\begin{rem}
Note that Theorem \ref{thm:PBXEX}(i) is true for $2\leq|X|<\aleph_0$ as well.  However, if $X$ is finite, then $\PB_X=\la\EX\cup\{\al,\be\}\ra$ if and only if $\S_X=\la\al,\be\ra$.  This all follows from the proof of \cite[Proposition 3.16]{DEG2017}.
\end{rem}

\section{Idempotents and two-sided units}\label{sect:FX}

We now turn our attention to the submonoid $\F(\PB_X) = \la E(\PB_X)\cup\G(\PB_X)\ra$ of $\PB_X$ generated by its idempotents and (two-sided) units.  We will continue to write $\S_X=\G(\PB_X)$, $\EX=\E(\PB_X)$, and so on, and from now on, we will also write $\FX=\F(\PB_X)$.  By Lemma \ref{lem:factorisation}(iii), we have $\FX=\EX\S_X=\S_X\EX$.  In Theorem \ref{thm:FPBX}, we characterise the elements of $\FX$.  Theorems \ref{thm:PBXFX}, \ref{thm:FXEX} and \ref{thm:FXSX} calculate $\relrank{\PB_X}{\FX}$, $\relrank{\FX}{\EX}$ and $\relrank{\FX}{\S_X}$, respectively; these theorems also characterise the minimal-size generating sets modulo the stated submonoids.

In order to prove the main results, we will need several preparatory lemmas.  Although the main focus of the current section is idempotents and \emph{two-sided} units, some of these lemmas hold in the larger submonoids
\[
\F_L(\PB_X)=\la E(\PB_X)\cup \G_L(\PB_X)\ra \AND \F_R(\PB_X)=\la E(\PB_X)\cup \G_R(\PB_X)\ra,
\]
and will also be of use when we study these submonoids in Section \ref{sect:FLFR}.  For simplicity, we will denote these submonoids by $\FLX$ and $\FRX$, respectively.  The elements of these monoids are described in Theorem \ref{thm:FLFR}.

We begin with a characterisation of the elements of $\FX=\F(\PB_X)$.  In what follows, we will often use the next result without explicit reference; since its statement and proof hold regardless of whether $X$ is finite or infinite, we make no restrictions on the size of $X$.

\begin{thm}\label{thm:FPBX}
If $X$ is an arbitrary set, then
\[
\FX = E(\PB_X)\S_X = \S_X E(\PB_X) = \set{\al\in\PB_X}{\defect(\al)=\codef(\al)}.
\]
\end{thm}

\pf First note that $E(M)\G(M)=\G(M)E(M)\sub\F(M)$ for any monoid $M$.  It remains to show that
\bit\bmc2
\itemnit{i} $\FX\sub\set{\al\in\PB_X}{\defect(\al)=\codef(\al)}$, 
\itemnit{ii} $\set{\al\in\PB_X}{\defect(\al)=\codef(\al)}\sub E(\PB_X)\S_X$.
\emc\eit
(i).  By Lemma \ref{lem:finitedefect}, $\defect(\al)=\codef(\al)$ for all $\al\in\PB_X$ if $X$ is finite.  So suppose $X$ is infinite, and let~$\al\in\FX$.  Then $\al=\be\ga$ for some $\be\in\EX$ and $\ga\in \S_X$.  Then $\defect(\be)=\codef(\be)$, by Theorem \ref{thm:IGPBX} (cf.~Remark \ref{rem:EXinFX}).  Lemma~\ref{lem:action} then gives $\defect(\al)=\defect(\be)=\codef(\be)=\codef(\al)$.

\pfitem{ii}  Suppose $\al\in\PB_X$ is such that $\defect(\al)=\codef(\al)$, and write $\al=\partn{a_i}{C_j}{b_i}{D_k}$.  Since $\defect(\al)=\codef(\al)$, there is a permutation $\be\in\S_X$ such that $a_i\be=b_i$ for all $i\in I$.  But then $\al=(\al\be^{-1})\be$, with $\al\be^{-1}=\partn{a_i}{C_j}{a_i}{D_k\be^{-1}}$ an idempotent, by Lemma \ref{lem:idempotents}.
\epf

\begin{rem}\label{rem:UR}
An element $x$ of a monoid $M$ is \emph{unit regular} if $x=xax$ for some unit $a\in\G(M)$.  As noted in \cite[Section 3]{FL1993}, $x$ is unit regular if and only if $x=eg$ for some idempotent $e\in E(M)$ and unit~$g\in\G(M)$.  Thus, Theorem \ref{thm:FPBX} shows that $\FX=E(\PB_X)\G(\PB_X)$ is the set of all unit regular elements of~$\FX$, and that the unit regular elements form a submonoid of $\PB_X$; cf.~\cite[Corollary 3.7]{FL1993}.  The unit regular elements of an arbitrary monoid do not necessarily form a submonoid; for example, this is not the case for finite partition monoids, as may easily be shown using GAP \cite{GAP}.  Unit regularity also plays an important role in ring theory; see for example \cite{Ehrlich1976,Handelman1977}.
\end{rem}

Now that we have characterised the elements of $\FX$, we can calculate its relative rank in $\PB_X$.

\begin{thm}\label{thm:PBXFX}
Let $X$ be an infinite set.
\bit
\itemit{i} We have $\relrank{\PB_X}{\FX}=2$.  
\itemit{ii} If $\al,\be\in\PB_X$, then $\PB_X=\la \FX\cup\{\al,\be\}\ra$ if and only if (renaming if necessary) $\al\in\GRX$, $\be\in\GLX$ and $\codef(\al)=\defect(\be)=|X|$.
\eit
\end{thm}

\pf
(i).  Since $\EX\sub\FX$, Theorem~\ref{thm:PBXEX} gives $\relrank{\PB_X}{\FX} \leq \relrank{\PB_X}{\EX} = 2$.  Lemma \ref{lem:rankMGE}(iv) gives the reverse inequality.

\pfitem{ii}  
If the stated conditions on $\al,\be$ hold, then Lemma \ref{lem:alEXbe} gives $\PB_X=\la\EX\cup\{\al,\be\}\ra$, so that certainly $\PB_X=\la\FX\cup\{\al,\be\}\ra$.  Conversely, suppose $\PB_X=\la\FX\cup\{\al,\be\}\ra$.  The proof of Theorem \ref{thm:PBXEX} works almost unmodified to show that $\al,\be$ satisfy the stated conditions.  The only difference is that the elements $\de_1,\ldots,\de_k$ used during the proof belong now to $\FX\cup\{\al,\be\}$, rather than to $\EX\cup\{\al,\be\}$.  The key property of elements~$\ve\in\EX$ used in the proof of Theorem \ref{thm:PBXEX} was that $\defect(\ve)=\codef(\ve)$; but this is also true if instead $\ve\in\FX$, by Theorem~\ref{thm:FPBX}.
\epf

\begin{rem}
Comparing Theorems \ref{thm:PBXEX} and \ref{thm:PBXFX}, we see that for any $\al,\be\in\PB_X$,
\[
\PB_X=\la \EX\cup\{\al,\be\}\ra \iff \PB_X=\la \FX\cup\{\al,\be\}\ra.
\]
\end{rem}

Since $\FX$ contains both $\S_X$ and $\EX$ as submonoids, we would naturally like to calculate the relative rank of~$\FX$ modulo these two submonoids.  The case of $\EX$ is easily dealt with; the following is an immediate consequence of Lemma~\ref{lem:rankME} (cf.~Remark \ref{rem:rankME2}), and the fact that $\rank(\S_X)=|\S_X|=2^{|X|}$ (as $\S_X$ is uncountable).

\begin{thm}\label{thm:FXEX}
Let $X$ be an infinite set.
\bit
\itemit{i} We have $\relrank{\FX}{\EX}=2^{|X|}$.  
\itemit{ii} If $\Om\sub \FX$, then $\FX=\la\EX\cup \Om\ra$ if and only if $\S_X=\la\S_X\cap\Om\ra$. \epfres
\eit
\end{thm}

The value of $\relrank{\FX}{\S_X}$ is harder to determine; again, it involves the number of infinite cardinals not exceeding $|X|$.  As noted above, some of the preliminary results we require will be formulated so as to be of use when we study the larger monoids $\FLX=\F_L(\PB_X)$ and $\FRX=\F_R(\PB_X)$ in Section \ref{sect:FLFR}.  We begin by characterising the elements of $\FLX$ and $\FRX$.  Recall that we write $\GLX=\G_L(\PB_X)$ and $\GRX=\G_R(\PB_X)$.  Also recall that $\FLX=\EX\GLX=\FX\GLX$ (cf.~Remark \ref{rem:factorisation}).  Again, we will often use the next result without explicit reference.

\begin{thm}\label{thm:FLFR}
If $X$ is an arbitrary set, then
\bit
\itemit{i} $\FLX = E(\PB_X)\GLX = \set{\al\in\PB_X}{\codef(\al)\leq\defect(\al)}$,
\itemit{ii} $\FRX = \GRX E(\PB_X) = \set{\al\in\PB_X}{\defect(\al)\leq\codef(\al)}$.
\eit
\end{thm}

\pf
We just prove (i), as (ii) is dual.  Since $E(M)\G_L(M)\sub\F_L(M)$ for any monoid $M$, it suffices to show that
\bit\bmc2
\itemnit{a} $\FLX \sub \set{\al\in\PB_X}{\codef(\al)\leq\defect(\al)}$,
\itemnit{b} $\bigset{\al\in\PB_X}{\codef(\al)\leq\defect(\al)}\sub E(\PB_X)\GLX$.
\emc\eit
(a).  Suppose $\al\in\FLX$, and write $\al=\be\ga$, where $\be\in\FX$ and $\ga\in \GLX$.  By Theorem \ref{thm:FPBX} and Lemma \ref{lem:PBunits}(ii), we have $\defect(\be)=\codef(\be)$ and $\codef(\ga)=0$.  Combined with parts (iv) and (v) of Lemma~\ref{lem:shdef}, it follows that
$\codef(\al) = \codef(\be\ga) \leq \codef(\be)+\codef(\ga) = \defect(\be) \leq \defect(\be\ga) =\defect(\al)$.

\pfitem{b}  Suppose $\al\in\PB_X$ is such that $\codef(\al)\leq\defect(\al)$.  Write $\al=\sqpartn{a_i}{C_j}{b_i}{D_k}$, recalling that this notation lists \emph{all} of the non-transversals of $\al$, not only the hooks.  Choose an injective map $\phi:\Codef(\al)\to\Defect(\al)$, and let $U=X\sm(\dom(\al)\cup\Codef(\al)\phi)$.  Let $\be\in\PB_X$ have transversals, and upper and lower non-transverals 
\[
\bigset{\{a_i,a_i'\}}{i\in I} \COMMA \set{C_j}{j\in J} \COMMA {\set{D_k\phi}{k\in K}\cup\bigset{\{u'\}}{u\in U}},
\]
respectively.  Let $\ga\in\PB_X$ have transversals and upper non-transversals
\[
\bigset{\{a_i,b_i'\}}{i\in I}\cup\bigset{\{x\phi,x'\}}{x\in\Codef(\al)} \AND \bigset{\{u\}}{u\in U},
\]
respectively.  Then $\al=\be\ga$, with $\be\in E(\PB_X)$ by Lemma \ref{lem:idempotents}, and $\ga\in \GLX$ by Lemma~\ref{lem:PBunits}(ii).
\epf

\begin{rem}[cf.~Remark \ref{rem:UR}]
\label{rem:LUR}
An element $x$ of a monoid $M$ is \emph{right-unit regular} if $x=xax$ for some right unit $a\in\G_R(M)$; left-unit regularity is defined analogously.  (See \cite{Chen1997} for the corresponding concept in ring theory.)  If we write $U_R(M)$ for the set of all right-unit regular elements of $M$, then one may show that
\begin{equation}\label{eq:URM}
\G_L(M)E(M) \sub U_R(M) \sub E(M)\G_L(M).
\end{equation}
We do not have $E(M)\G_L(M)\sub \G_L(M)E(M)$ in general, however.  Indeed, consider an element $\al\in\PB_X$ with $h(\al)=s^*(\al)=1$ and $s(\al)=h^*(\al)=0$.  Then $\al\in\FLX=E(\PB_X)\GLX$, by Theorem \ref{thm:FLFR}(i).  However, we claim that $\al\not\in\GLX E(\PB_X)$.  To see why this is the case, suppose to the contrary that $\al=\be\ga$, where $\be\in\GLX$ and $\ga\in E(\PB_X)$.  
Now, $\codef(\ga)\leq\codef(\be\ga)=\codef(\al)=1$.  If $\codef(\ga)=1$, then since also $\defect(\ga)=\codef(\ga)$, by Remark \ref{rem:EXinFX}, we must have $s(\ga)=1$, in which case Lemma \ref{lem:shdef2}(i) would give $0=s(\al)=s(\be\ga)=s(\be)+s(\ga)\geq s(\ga)=1$, a contradiction.  Thus, $\codef(\ga)=0$.  But then Theorem \ref{thm:IGPBX} gives $\sh(\ga)=0$ (since certainly $\ga\in\EX$), and so $\ga=1$, and $\al=\be\ga=\be\in\GLX$, contradicting $\codef(\al)\not=0$.

On the other hand, it is not hard to show that $U_R(\PB_X)=E(\PB_X)\GLX=\FLX$.  Indeed, by \eqref{eq:URM}, it suffices to show that $\FLX\sub U_R(\PB_X)$.  To do so, let $\al\in\FLX$, write $U=\codom(\al)$ and $V=\Codef(\al)$, and fix an injective map $\phi:V\to\Defect(\al)$.  Then it is easy to see that $\al=\al\be\al$, where $\be=\partnV uv{u\al^{-1}}{v\phi}\in\GRX$.
\end{rem}

The next result will be used often, and highlights an important property of the elements of $\FLX$.

\begin{lemma}\label{lem:defFLX}
If $X$ is an arbitrary set, and if $\al_1,\ldots,\al_k\in\FLX$, then $\defect(\al_i\cdots\al_j)\leq\defect(\al_1\cdots\al_k)$ for all $1\leq i\leq j\leq k$.
\end{lemma}

\pf
We have $\defect(\al_i\cdots\al_j)\leq\defect(\al_i\cdots\al_j(\al_{j+1}\cdots\al_k))\leq\defect((\al_1\cdots\al_{i-1})\al_i\cdots\al_j(\al_{j+1}\cdots\al_k))$, by Lemmas \ref{lem:shdef}(iv) and~\ref{lem:shdef4}, respectively.
\epf

The next lemma introduces a certain special kind of subset $\Om$ of $\FX$ that will play an important role in this section and the next.  The lemma immediately following will show that for any such subset $\Om$, we have $\FX=\la\S_X\cup\Om\ra$, and we will see later that $\Om$ is of minimal possible size with respect to this property.

\begin{lemma}\label{lem:U1}
Let $X$ be an infinite set, and let $\Om = \set{\al_\mu,\be_\mu}{\mu=1\text{ or }\aleph_0\leq\mu\leq|X|}\sub\FX$, where for each~$\mu$,
\[
t(\al_\mu)=t(\be_\mu)=|X| \COMMA s(\al_\mu)<\mu=h(\al_\mu) \COMMA h(\be_\mu)<\mu=s(\be_\mu),
\]
and either
\bit
\itemit{i} $s^*(\al_\mu)<\mu=h^*(\al_\mu)$ and $h^*(\be_\mu)<\mu=s^*(\be_\mu)$, or 
\itemit{ii} $h^*(\al_\mu)<\mu=s^*(\al_\mu)$ and $s^*(\be_\mu)<\mu=h^*(\be_\mu)$.
\eit
Then for any cardinal $0\leq\nu\leq|X|$, there exists $\si_\nu,\tau_\nu\in\la\S_X\cup \Om\ra$ such that
\[
t(\si_\nu)=t(\tau_\nu)=|X|
\COMMA
h(\si_\nu)=h^*(\si_\nu) = s(\tau_\nu)=s^*(\tau_\nu) = \nu
\COMMA
s(\si_\nu)=s^*(\si_\nu) = h(\tau_\nu)=h^*(\tau_\nu) = 0.
\]
\end{lemma}

\pf
We use transfinite induction.  First, let $\si_0,\tau_0$ be any elements of $\S_X$.  Now suppose $1\leq\nu\leq|X|$ is such that elements $\si_\ka,\tau_\ka$ of the desired form exist for all cardinals $\ka<\nu$.

\pfcase{1}  Suppose first that $1\leq\nu<\aleph_0$.  Consider the elements $\si_{\nu-1},\tau_{\nu-1}\in\la\S_X\cup\Om\ra$, guaranteed to exist by the above induction hypothesis.  
Note that when $\mu=1$, we must be in case (i), as $\codef(\al_1)=\defect(\al_1)=s(\al_1)+2h(\al_1)=2$, and similarly $\codef(\be_1)=1$.  Write $\Codef(\al_1)=\{a,b\}$ and $\Codef(\be_1)=\{c\}$, where $a,b,c\in X$, and where $a,b$ are distinct (but note that possibly $c\in\{a,b\}$). 
Choose distinct $u,v\in\dom(\si_{\nu-1})$ and any $w\in\dom(\tau_{\nu-1})$, and let $\ga,\de\in\S_X$ be such that $a\ga=u$, $b\ga=v$ and $c\de=w$.  Then $\si_\nu=\al_1\ga\si_{\nu-1}$ and $\tau_\nu=\be_1\de\tau_{\nu-1}$ have the desired properties.  (For future reference, we note that in fact $\si_\nu,\tau_\nu\in\la\S_X\cup\{\al_1,\be_1\}\ra$.)

\pfcase{2}  Suppose now that $\aleph_0\leq\nu\leq|X|$.  We begin by proving the existence of $\si_\nu$.  We first claim that:
\begin{equation}\label{eq:ve}
\text{there exists $\ga\in\la\S_X\cup \Om\ra$ such that $t(\ga)=|X|$, $h(\ga)=\nu$ and $s(\ga)=0$.  }
\end{equation}
By assumption, $\al_\nu\in \Om$ satisfies $t(\al_\nu)=|X|$ and $s(\al_\nu)<\nu=h(\al_\nu)$.  Clearly \eqref{eq:ve} holds if $s(\al_\nu)=0$ (we take $\ga=\al_\nu$), so suppose $s(\al_\nu)\geq1$.  For simplicity, we write $\ka=s(\al_\nu)$.  (Note that $\ka$ might be finite.)  By the induction hypothesis, since $\ka<\nu$, there exists $\si_\ka\in\la\S_X\cup \Om\ra$ with $t(\si_\ka)=|X|$, $h(\si_\ka)=h^*(\si_\ka)=\ka$ and $s(\si_\ka)=s^*(\si_\ka)=0$.
Let $V$ be the set of upper singletons of $\al_\nu$.  Since $\codef(\si_\ka)=2\ka$ and $|V|=\ka$, with $\ka<|X|$, there is a permutation $\pi\in \S_X$ that maps $\Codef(\si_\ka)$ bijectively onto a subset $W\sub V\cup\dom(\al_\nu)$ with $V\sub W$ (we may take $W=V$ if $\ka$ is infinite).  Then $\ga=\si_\ka\pi\al_\nu$ satisfies the conditions of \eqref{eq:ve}.  A dual argument (using~$\al_\nu$ in case (i) or $\be_\nu$ in case (ii)) shows that:
\begin{equation}\label{eq:eta}
\text{there exists $\de\in\la\S_X\cup \Om\ra$ such that $t(\de)=|X|$, $h^*(\de)=\nu$ and $s^*(\de)=0$.}
\end{equation}
With $\ga$ and $\de$ as in \eqref{eq:ve} and \eqref{eq:eta}, and since $\ga,\de\in\FX$, Theorem \ref{thm:FPBX} gives
\[
\codef(\ga) = \defect(\ga) = 2\nu = \codef(\de) = \defect(\de).
\]
Since also $t(\ga)=|X|=t(\de)$, there exists a permutation $\ve\in\S_X$ that maps $\codom(\ga)$ bijectively onto $\dom(\de)$, and it follows that $\si_\nu = \ga\ve\de$ has the desired properties.

The existence of $\tau_\nu$ is demonstrated in almost identical fashion, with the symbols $s$ and $h$ swapped, and using $\be_\nu$ in place of $\al_\nu$.  The only place where special care is required is as follows.  In order to prove the analogue of claim \eqref{eq:ve}---i.e., to prove that there exists $\ga\in\la\S_X\cup \Om\ra$ such that $t(\ga)=|X|$, $s(\ga)=\nu$ and $h(\ga)=0$---we write $\ka=h(\be_\nu)$, but we then utilise the element $\tau_{2\ka}$ (rather than $\si_\ka$) to ensure that a permutation $\pi\in\S_X$ exists so that $\ga=\tau_{2\ka}\pi\be_\nu$ has the desired properties.  Note that $2\ka<\nu$ follows from $\ka<\nu$, since $\nu\geq\aleph_0$.
\epf

\begin{lemma}\label{lem:U2}
Let $X$ be an infinite set, and let $\Om\sub\FX$ be as in Lemma \ref{lem:U1}.  Then $\FX=\la\S_X\cup \Om\ra$.
\end{lemma}

\pf
Let $\ga\in\FX$ be arbitrary, and write $\lam=t(\ga)$, $\mu_1=h(\ga)$, $\mu_2=h^*(\ga)$, $\nu_1=s(\ga)$ and $\nu_2=s^*(\ga)$.  By Lemma \ref{lem:action}, it suffices to demonstrate the existence of any $\de\in\la\S_X\cup \Om\ra$ with
\[
t(\de)=\lam \COMMA h(\de)=\mu_1 \COMMA h^*(\de)=\mu_2 \COMMA s(\de)=\nu_1 \COMMA s^*(\de)=\nu_2,
\]
since then $\ga\in\S_X\de\S_X$.  
Theorem \ref{thm:FPBX} gives $2\mu_1+\nu_1 = \defect(\ga) = \codef(\ga) = 2\mu_2+\nu_2$.

\pfcase{1} Suppose first that $\lam=|X|$.  Consider the elements $\si_{\mu_1},\si_{\mu_2},\tau_{\nu_1},\tau_{\nu_2}\in\la\S_X\cup \Om\ra$ as given by Lemma~\ref{lem:U1}.  Post-multiplying $\si_{\mu_1}$ by a suitable permutation if necessary, and keeping in mind that ${t(\si_{\mu_1})=t(\tau_{\nu_1})=|X|}$, we may assume without loss of generality that $\Defect(\tau_{\nu_1})\sub\codom(\si_{\mu_1})$ and also ${|{\codom(\si_{\mu_1})\sm\Defect(\tau_{\nu_1})}|=|X|}$.  It then follows that
\[
t(\si_{\mu_1}\tau_{\nu_1}) = |X| \COMMA h(\si_{\mu_1}\tau_{\nu_1}) = \mu_1 \COMMA s(\si_{\mu_1}\tau_{\nu_1}) = \nu_1.
\]
Similarly, we may assume that
\[
t(\si_{\mu_2}\tau_{\nu_2}) = |X| \COMMA h^*(\si_{\mu_2}\tau_{\nu_2}) = \mu_2 \COMMA s^*(\si_{\mu_2}\tau_{\nu_2}) = \nu_2.
\]
Since $\codef(\si_{\mu_1}\tau_{\nu_1})=\defect(\si_{\mu_1}\tau_{\nu_1})=2\mu_1+\nu_1=2\mu_2+\nu_2=\codef(\si_{\mu_2}\tau_{\nu_2}) =\defect(\si_{\mu_2}\tau_{\nu_2})$, and since $t(\si_{\mu_1}\tau_{\nu_1})=t(\si_{\mu_2}\tau_{\nu_2})=|X|$, there is a permutation $\pi\in\S_X$ that maps $\codom(\si_{\mu_1}\tau_{\nu_1})$ bijectively onto $\dom(\si_{\mu_2}\tau_{\nu_2})$.  Then $\de=(\si_{\mu_1}\tau_{\nu_1})\pi(\si_{\mu_2}\tau_{\nu_2})$ has the desired parameter values.

\pfcase{2}  Suppose now that $\lam<|X|$.  Since $|X|=\lam+\defect(\ga)=\lam+\codef(\ga)$, we have ${\defect(\ga)=\codef(\ga)=|X|}$.  Write $\dom(\ga)=A$.  By Lemma \ref{lem:action}, post-multiplying by a permutation if necessary, we may assume that~${\ga=\sqpartn a{B_x}a{C_x}}$, again recalling that this notation lists \emph{all} of the non-transversals.  Clearly $\ga=\de\ve$, where $\de=\sqpartn a{B_x}a{B_x}$ and $\ve=\sqpartn a{C_x}a{C_x}$; by symmetry, it suffices to show that $\de\in\la\S_X\cup\Om\ra$.  Fix a decomposition $X=Y\sqcup Z$ where $|X|=|Y|=|Z|$, and put $U=\bigcup_{y\in Y}B_y$ and $V=\bigcup_{z\in Z}B_z$.  Then $\de=\de_1\de_2$, where $\de_1=\sqpartn x{B_y}x{B_y}_{x\in A\cup V,\ y\in Y}$ and $\de_2=\sqpartn x{B_z}x{B_z}_{x\in A\cup U,\ z\in Z}$.  By Case 1, we have $\de_1,\de_2\in\la\S_X\cup\Om\ra$.
\epf

The previous lemma will be used to give an upper bound on the size of generating sets for $\FX$ modulo~$\S_X$.  The next two lemmas work towards establishing that this is also a lower bound. 

\begin{lemma}\label{lem:F2}
Let $X$ be an infinite set, and let $\mu$ be a cardinal such that either $\mu=1$ or $\aleph_0\leq\mu\leq|X|$.  Let $\al_1,\ldots,\al_k\in\FLX$, write $\be=\al_1\cdots\al_k$, and suppose $\be\in\FX$ and $t(\be)=|X|$.
\bit
\itemit{i} If $s(\be)<\mu= h(\be)$, then $\al_i\in\FX$, $s(\al_i)<\mu= h(\al_i)$ and $t(\al_i)=|X|$ for some $i\in\{1,\ldots,k\}$.
\itemit{ii} If $h(\be)<\mu= s(\be)$, then $\al_i\in\FX$, $h(\al_i)<\mu= s(\al_i)$ and $t(\al_i)=|X|$ for some $i\in\{1,\ldots,k\}$.
\eit
\end{lemma}

\pf
Throughout the proof, we will frequently use Lemma \ref{lem:shdef} and Theorems \ref{thm:FPBX} and \ref{thm:FLFR} without explicit mention.

\pfitem{ii}  Suppose $h(\be)<\mu= s(\be)$.  Note then that $\codef(\be)=\defect(\be)=s(\be)+2h(\be)=\mu$.  For $0\leq i\leq k$, put $\be_i=\al_1\cdots\al_i$.  Then for any such $i$, Lemma \ref{lem:defFLX} gives
$\codef(\be_i)\leq\defect(\be_i)\leq\defect(\al_1\cdots\al_k)=\defect(\be)=\mu$.
Since $\codef(\be_0)=0$, we may define $j=\max\set{i}{\codef(\be_i)<\mu}$.  Since $\codef(\be_k)=\codef(\be)=\mu$, we have $0\leq j<k$.  Since Corollary \ref{cor:tq}(i) gives $t(\al_{j+1})=|X|$, the proof will be complete if we can show that
\bit\bmc3
\itemnit{a} $\al_{j+1}\in\FX$,
\itemnit{b} $h(\al_{j+1})<\mu$,
\itemnit{c} $s(\al_{j+1})=\mu$.
\emc\eit
From the definition of $j$, we have 
\[
\codef(\be_j)<\mu \AND \mu\leq\codef(\be_{j+1})=\codef(\be_j\al_{j+1})\leq\codef(\be_j)+\codef(\al_{j+1}) . 
\]
Thus, by the form of $\mu$, $\codef(\al_{j+1})\geq\mu$.  Together with Lemma \ref{lem:defFLX}, it follows that
\[
\mu\leq\codef(\al_{j+1})\leq\defect(\al_{j+1})\leq\defect(\al_1\cdots\al_k)=\defect(\be)=\mu.
\]
Thus, $\defect(\al_{j+1})=\codef(\al_{j+1})=\mu$.  In particular, (a) holds.

To show that (b) holds, suppose to the contrary that $h(\al_{j+1})\geq\mu$.  Then $h(\al_{j+1})\geq\mu>\codef(\be_j)$, so Lemma \ref{lem:shdef3}(ii) gives $h(\be_j\al_{j+1}) \geq h(\be_j)+h(\al_{j+1})-\codef(\be_j) \geq h(\al_{j+1})-\codef(\be_j)$.  If $\mu=1$, then $\codef(\be_j)=0$, while if $\mu\geq\aleph_0$, then $h(\al_{j+1})\geq\mu\geq\aleph_0$.  Thus, in either case, ${h(\al_{j+1})-\codef(\be_j)=h(\al_{j+1})}$, and so $\mu>h(\be)=h(\be_j\al_{j+1}\cdots \al_k)\geq h(\be_j\al_{j+1}) \geq h(\al_{j+1})-\codef(\be_j)=h(\al_{j+1})\geq\mu$, a contradiction.  This completes the proof of (b).

In light of $\mu=\defect(\al_{j+1})=s(\al_{j+1})+2h(\al_{j+1})$ and $h(\al_{j+1})<\mu$, and by the form of $\mu$, it follows that $s(\al_{j+1})=\mu$, giving (c).  

\pfitem{i}  Suppose $s(\be)<\mu=h(\be)$.  This time, $\codef(\be)=\defect(\be)=2\mu$.  If $\mu\geq\aleph_0$, then $2\mu=\mu$, and the proof carries on in essentially the same way as in (ii), above.  For the $\mu=1$ case, we define $j=\min\set{i}{\defect(\al_i)\geq1}$, noting that $\codef(\al_i)\leq\defect(\al_i)=0$ for all $i<j$, so that ${\ga=\al_1\cdots\al_{j-1}}$ belongs to $\S_X$.  Then $\ga^{-1}\be=\al_j\cdots\al_k$, with $s(\ga^{-1}\be)=s(\be)=0$, by Lemma \ref{lem:action}.  Lemma \ref{lem:defFLX} gives ${1\leq\defect(\al_j)\leq\defect(\al_1\cdots\al_k)=\defect(\be)=2}$.  Also, $s(\al_j)\leq s(\al_j\cdots\al_k)=s(\ga^{-1}\be)=0$, so that $s(\al_j)=0$.  Together with $1\leq\defect(\al_j)\leq2$, it then follows that $\defect(\al_j)=2$, and $h(\al_j)=1$.
\epf

\begin{lemma}\label{lem:F3}
Let $X$ be an infinite set, and suppose $\Om\sub\FLX$ is such that $\FX\sub\la\S_X\cup \Om\ra$.  For any cardinal~$\mu$ such that $\mu=1$ or $\aleph_0\leq\mu\leq|X|$,
\bit
\itemit{i} there exists $\al_\mu\in \Om$ with $\al_\mu\in\FX$, $s(\al_\mu)<\mu=h(\al_\mu)$ and $t(\al_\mu)=|X|$,
\itemit{ii} there exists $\be_\mu\in \Om$ with $\be_\mu\in\FX$, $h(\be_\mu)<\mu=s(\be_\mu)$ and $t(\be_\mu)=|X|$. 
\eit
\end{lemma}

\pf
The proofs being essentially identical, we just prove (i).  Let $\ga\in\FX$ be such that $s(\ga)<\mu=h(\ga)$ and $t(\ga)=|X|$, where $\mu=1$ or $\aleph_0\leq\mu\leq|X|$, and consider an expression $\ga=\de_1\cdots\de_k$, where each of the factors belong to $\S_X\cup\Om$.  Lemma \ref{lem:F2}(i) says that one of the $\de_i$ satisfies $\de_i\in\FX$, $s(\de_i)<\mu=h(\de_i)$ and $t(\de_i)=|X|$.  So we may take $\al_\mu=\de_i$, noting that $\de_i\not\in\S_X$ (as $h(\de_i)>0$).
\epf

We have now gathered all the facts needed to prove the final main result of this section.

\begin{thm}\label{thm:FXSX}
Let $X$ be an infinite set, and let $\rho$ be the number of cardinals $\mu$ satisfying $\aleph_0\leq\mu\leq|X|$.
\bit
\itemit{i} We have $\relrank{\FX}{\S_X}=2+2\rho$.  
\itemit{ii} If $\rho<\aleph_0$, and if $\Om\sub \FX$ with $|\Om|=2+2\rho$, then $\FX=\la\S_X\cup \Om\ra$ if and only if $\Om$ has the form described in Lemma \ref{lem:U1}.
\eit
\end{thm}

\pf
First, if $\Om\sub\FX$ is of the form given in Lemma \ref{lem:U1}, then Lemma \ref{lem:U2} gives $\FX=\la\S_X\cup\Om\ra$.  This gives $\relrank{\FX}{\S_X}\leq|\Om|=2+2\rho$, and also the backwards implication in (ii).

Next, suppose $\Om\sub\FX$ is such that $\FX=\la\S_X\cup\Om\ra$ and $|\Om|=\relrank{\FX}{\S_X}$.  
By Lemma \ref{lem:F3} (noting that $\Om\sub\FX\sub\FLX$), for any cardinal $\mu$ such that $\mu=1$ or $\aleph_0\leq\mu\leq|X|$,
\bit
\itemnit{a} there exists $\al_\mu\in \Om$ with $s(\al_\mu)<\mu=h(\al_\mu)$ and $t(\al_\mu)=|X|$,
\itemnit{b} there exists $\be_\mu\in \Om$ with $h(\be_\mu)<\mu=s(\be_\mu)$ and $t(\be_\mu)=|X|$.
\eit
The elements from (a) and (b) are distinct, and there are $2+2\rho$ of them.  Thus, $\relrank{\FX}{\S_X}=|\Om|\geq2+2\rho$.  This completes the proof of (i).  

To complete the proof of (ii), suppose from now on that $\rho<\aleph_0$.
By the dual version of Lemma \ref{lem:F3} (noting also that $\Om\sub\FX\sub\FRX$), for any cardinal $\mu$ such that $\mu=1$ or $\aleph_0\leq\mu\leq|X|$,
\bit
\itemnit{c} there exists $\ga_\mu\in \Om$ with $s^*(\ga_\mu)<\mu=h^*(\ga_\mu)$ and $t(\ga_\mu)=|X|$,
\itemnit{d} there exists $\de_\mu\in \Om$ with $h^*(\de_\mu)<\mu=s^*(\de_\mu)$ and $t(\de_\mu)=|X|$.
\eit
The elements from (c) and (d) are also distinct, and there are $2+2\rho$ of them.  Since $|\Om|=2+2\rho<\aleph_0$, it follows that
\[
\Om = \set{\al_\mu,\be_\mu}{\mu=1\text{ or }\aleph_0\leq\mu\leq|X|} = \set{\ga_\mu,\de_\mu}{\mu=1\text{ or }\aleph_0\leq\mu\leq|X|}.
\]
Now, $\defect(\al_1)=\codef(\ga_1)=2$, $\defect(\be_1)=\codef(\de_1)=1$ and  $\defect(\al_\mu)=\defect(\be_\mu)=\codef(\ga_\mu)=\codef(\de_\mu)=\mu$ for all $\aleph_0\leq\mu\leq|X|$.  Since $\Om\sub\FX$, all these elements have equal defect and codefect, so ${\al_1=\ga_1}$, $\be_1=\de_1$, and $\{\al_\mu,\be_\mu\}=\{\ga_\mu,\de_\mu\}$ for $\mu\geq\aleph_0$.  It quickly follows that $\Om=\set{\al_\mu,\be_\mu}{\mu=1\text{ or }\aleph_0\leq\mu\leq|X|}$ satisfies the conditions of Lemma \ref{lem:U1}.
\epf

\begin{rem}
Again, the assumption $\rho<\aleph_0$ is essential in Theorem \ref{thm:FXSX}(ii).  Indeed, if $\rho\geq\aleph_0$, then not only could we add a superfluous element $\ga$ to a generating set $\Om$ of the form given in Lemma \ref{lem:U1} without increasing its size (cf.~Remark \ref{rem:GLXSX}), but the elements given in (a)--(b) in the above proof might have little overlap with the elements given in (c)--(d).  For example, we might have $s^*(\al_1)=s(\ga_1)=2$, so that~$\al_1\not=\ga_1$.
\end{rem}

\section{Idempotents and one-sided units}\label{sect:FLFR}

This section concerns the submonoids 
\[
\F_L(\PB_X)=\la E(\PB_X)\cup \G_L(\PB_X)\ra \AND \F_R(\PB_X)=\la E(\PB_X)\cup \G_R(\PB_X)\ra
\]
of $\PB_X$ generated by its idempotents and left units, or idempotents and right units, respectively.
We will continue to use the abbreviations $\EX=\E(\PB_X)$, $\GLX=\G_L(\PB_X)$, $\FLX=\F_L(\PB_X)$, and so on.  The elements of $\FLX$ and $\FRX$ were characterised in Theorem \ref{thm:FLFR}.  
The main results of this section calculate the relative rank of $\PB_X$ modulo $\FLX$ (Theorem \ref{thm:PBXFLX}), and the relative ranks of $\FLX$ modulo each of the submonoids $\FX$, $\EX$, $\GLX$ and $\S_X$ (Theorems \ref{thm:FLXFX}, \ref{thm:FLXEX}, \ref{thm:FLXGLX} and \ref{thm:FLXSX}, respectively); we also classify the minimal-size generating sets modulo the stated submonoids.
The corresponding statements for~$\FRX$ are dual, and are easily deduced.

\begin{thm}\label{thm:PBXFLX}
Let $X$ be an infinite set.  
\bit
\itemit{i} We have $\relrank{\PB_X}{\FLX}=1$.  
\itemit{ii} If $\al\in\PB_X$, then $\PB_X=\la \FLX\cup\{\al\}\ra$ if and only if $\al\in\GRX$ and $\codef(\al)=|X|$.  
\eit
\end{thm}

\pf
If $\al\in\GRX$ and $\codef(\al)=|X|$, then for any $\be\in\GLX$ with $\defect(\be)=|X|$, Theorem \ref{thm:PBXFX}(ii) gives  $\PB_X=\la\FX\cup\{\al,\be\}\ra\sub\la\FLX\cup\{\al\}\ra\sub\PB_X$.  This gives the backwards implication in (ii), and also the inequality $\relrank{\PB_X}{\FLX}\leq1$; since $\PB_X\not=\FLX$, the reverse inequality is obvious.  

For the forwards implication in (ii), suppose $\al\in\PB_X$ is such that $\PB_X=\la\FLX\cup\{\al\}\ra$.  Lemma \ref{lem:rankMGE}(i) says that $(\FLX\sm\FX)\cup\{\al\}$ contains an element of $\GRX$.  Every element $\ga$ of $\FLX\sm\FX$ satisfies $\codef(\ga)<\defect(\ga)$, but every element $\de$ of $\GRX$ satisfies $\defect(\de)=0$.  It follows that $\al\in\GRX$.  It remains to show that $\codef(\al)=|X|$.  To do so, let $\si\in\PB_X$ be such that $\defect(\si)<|X|=\codef(\si)$, and consider an expression $\si=\ga_1\cdots\ga_k$, where $\ga_1,\ldots,\ga_k\in\FLX\cup\{\al\}$.  Corollary \ref{cor:tq}(ii) gives $\codef(\ga_i)=|X|$ for some $i$.  It suffices to show that $\ga_i=\al$.  To do so, suppose to the contrary that $\ga_i\in\FLX$, and note then that $|X|=\codef(\ga_i)\leq\defect(\ga_i)$, giving $\defect(\ga_i)=|X|$.  Then Lemma \ref{lem:defFLX} gives $|X|=\defect(\ga_i)\leq\defect(\ga_1\cdots\ga_k)=\defect(\si)<|X|$, a contradiction.
\epf

We now begin the task of calculating the relative rank of $\FLX$ modulo the submonoids mentioned above.

\begin{lemma}\label{lem:Om15}
Let $X$ be an infinite set, and let $\Om=\set{\al_\mu}{\mu=1\text{ or }\aleph_0\leq\mu\leq|X|}\sub\FLX$, where for each~$\mu$,
\[
\codef(\al_\mu)<\mu=\defect(\al_\mu).
\]
Then for any $0\leq\nu\leq|X|$, there exists $\si_\nu\in\la\S_X\cup\Om\ra$ with $\defect(\si_\nu)=\nu$ and $\codef(\si_\nu)=0$.
\end{lemma}

\pf
We use transfinite induction.  The result is clearly true for $\nu=0$ (take $\si_0$ to be any element of~$\S_X$).  Next, suppose $1\leq\nu\leq\aleph_0$ is such that the lemma holds for all cardinals $\ka<\nu$.  If $\nu<\aleph_0$, then $\si_\nu=\si_{\nu-1}\al_1$ has the desired properties; indeed, Lemma \ref{lem:PBunits}(ii) gives $\si_{\nu-1},\al_1\in \GLX$, and so $\si_\nu\in \GLX$, which gives $\codef(\si_\nu)=0$, while Lemma \ref{lem:shdef2}(v) gives $\defect(\si_\nu)=\defect(\si_{\nu-1})+\defect(\al_1)=\nu$.
Next suppose $\nu\geq\aleph_0$.  If $\codef(\al_\nu)=0$, then we just take $\si_\nu=\al_\nu$, so suppose $\codef(\al_\nu)\geq1$, and write $\ka=\codef(\al_\nu)$.  Since $\ka<\nu$, $\si_\ka\in\la\S_X\cup\Om\ra$ exists, by the induction hypothesis, and we have $t(\al_\nu)=t(\si_\ka)=|X|$.  Since also $\codef(\al_\nu)=\defect(\si_\ka)=\ka$, it follows that there is a permutation $\pi\in\S_X$ such that $\codom(\al_\nu)\pi=\dom(\si_\ka)$.  Then $\si_\nu=\al_\nu\pi\si_\ka$ has the desired properties.  
\epf

\begin{rem}\label{rem:Om15}
Since every element $\al$ of $\Om$ (as in Lemma \ref{lem:Om15}) has $\codef(\al)<|X|$, it follows that $t(\al)=|X|$ for all $\al\in\Om$.
\end{rem}

\begin{lemma}\label{lem:Om2}
If $X$ is an infinite set, and if $\Om\sub\FLX$ is as in Lemma \ref{lem:Om15}, then ${\FLX=\la \FX\cup\Om\ra}$.
\end{lemma}

\pf
Since $\FLX=\FX\GLX$ (cf.~Remark \ref{rem:factorisation}), it suffices to show that $\GLX\sub\la\FX\cup\Om\ra$.  So let $\ga=\partn{a_x}{B_i}x{}\in \GLX$ be arbitrary.  Put $\de=\partn{a_x}{B_i}{a_x}{}\in E(\PB_X)$, and write $\mu=\defect(\ga)$.  Then $\codef(\de)=\defect(\de)=\mu$ also.  Since $t(\de)=t(\ga)=|X|=t(\si_\mu)$, where $\si_\mu\in\la\S_X\cup\Om\ra$ is as in Lemma \ref{lem:Om15} (cf.~Remark \ref{rem:Om15}), there is a permutation $\pi\in\S_X$ such that $a_x\pi=x\si_\mu^{-1}$ for all $x\in X$.  It follows that $\ga=\de\pi\si_\mu\in\la\FX\cup\Om\ra$.
\epf

\begin{lemma}\label{lem:Om1}
Let $X$ be an infinite set, and suppose $\Om\sub\FLX$ is such that $\FLX=\la\FX\cup\Om\ra$.  Then $\Om$ contains a subset of the form described in Lemma \ref{lem:Om15}.
\end{lemma}

\pf
Let $\mu$ be a cardinal such that $\mu=1$ or $\aleph_0\leq\mu\leq|X|$.  We must show that there exists $\al\in\Om$ such that $\codef(\al)<\mu=\defect(\al)$.  Let $\si\in\FLX$ be such that ${\codef(\si)<\mu=\defect(\si)}$, and consider an expression $\si=\al_1\cdots\al_k$, where each factor belongs to $\FX\cup\Om$.  By Corollary~\ref{cor:tq}(ii), $\defect(\al_i)\geq\mu$ for some $i$.  Let $j=\max\set{i}{\defect(\al_i)\geq\mu}$.  Lemma \ref{lem:defFLX} gives ${\defect(\al_j)\leq\defect(\al_1\cdots\al_k)=\defect(\si)=\mu}$, and so $\defect(\al_j)=\mu$.  The proof will be complete if we can show that $\codef(\al_j)<\mu$, since then also $\al_j\not\in\FX$, which would give $\al_j\in\Om$.  Suppose to the contrary that $\codef(\al_j)\geq\mu$.  
Combined with $\codef(\al_j)\leq\defect(\al_j)=\mu$, it follows that $\codef(\al_j)=\mu$.  
Put $\be=\al_{j+1}\cdots\al_k$, and note that by Lemma \ref{lem:shdef}(iv) and the definition of $j$, $\defect(\be)\leq\defect(\al_{j+1})+\cdots+\defect(\al_k)<\mu$.  Thus, we have
\[
[\codef(\al_j)=1 \text{ or } \codef(\al_j)\geq\aleph_0] \AND \codef(\al_j)>\defect(\be),
\]
so we obtain $\codef(\al_j\be)\geq\codef(\al_j)=\mu$ from the dual of Lemma \ref{lem:shdef3}(iii).  Lemma \ref{lem:shdef}(v) then gives $\mu>\codef(\si)=\codef(\al_1\cdots\al_j\be)\geq\codef(\al_j\be)\geq\mu$, a contradiction.
\epf

The next result follows quickly from Lemmas \ref{lem:Om2} and \ref{lem:Om1}.

\begin{thm}\label{thm:FLXFX}
Let $X$ be an infinite set, and let $\rho$ be the number of cardinals $\mu$ satisfying $\aleph_0\leq\mu\leq|X|$.
\bit
\itemit{i} We have $\relrank{\FLX}{\FX}=1+\rho$.  
\itemit{ii} If $\rho<\aleph_0$, and if $\Om\sub \FLX$ with $|\Om|=1+\rho$, then $\FLX=\la \FX\cup \Om\ra$ if and only if $\Om$ has the form described in Lemma \ref{lem:Om15}. \epfres
\eit
\end{thm}

Since $\G(\FLX)=\S_X$ (cf.~Lemma \ref{lem:FFLFR}), and since $\FLX\sm\S_X$ is an ideal of $\FLX$ (cf.~Lemmas \ref{lem:equiv} and \ref{lem:FFLFR}), we may also quickly deal with the situation modulo $\EX$.  

\begin{thm}\label{thm:FLXEX}
Let $X$ be an infinite set.
\bit
\itemit{i} We have $\relrank{\FLX}{\EX}=2^{|X|}$.  
\itemit{ii} If $\Om\sub \FLX$, then $\FLX=\la \EX\cup \Om\ra$ if and only if $\S_X=\la\S_X\cap\Om\ra$ and $\Om$ contains a subset of the form described in Lemma \ref{lem:Om15}. 
\eit
\end{thm}

\pf
(i).  This follows from Lemma \ref{lem:rankME}(iii), and the fact that $\rank(\S_X)=|\S_X|=2^{|X|}=|\FLX|$.

\pfitem{ii}  If $\Om\sub\FLX$, then Lemma \ref{lem:rankME}(ii) says that $\FLX=\la\EX\cup\Om\ra$ if and only if $\S_X=\la\S_X\cap\Om\ra$ and $\FLX=\la\FX\cup(\Om\sm\S_X)\ra$.  By Lemmas \ref{lem:Om2} and \ref{lem:Om1}, this latter condition is equivalent to $\Om\sm\S_X$ (and hence $\Om$) having a subset of the form described in Lemma \ref{lem:Om15}.
\epf

Now we move on to the task of calculating $\relrank{\FLX}{\GLX}$.

\begin{lemma}\label{lem:FLXSX}
Let $X$ be an infinite set, and let $\Om=\set{\al_\mu,\be_\mu}{\mu=1\text{ or }\aleph_0\leq\mu\leq|X|}\sub\FX$, where for all~$\mu$,
\bit\bmc2
\itemit{i} $s(\al_\mu)<\mu=h(\al_\mu)$ and $t(\al_\mu)=|X|$,
\itemit{ii} $h(\be_\mu)<\mu=s(\be_\mu)$ and $t(\be_\mu)=|X|$,
\itemit{iii} $h^*(\al_\mu)=s^*(\be_\mu)=\mu$ or $s^*(\al_\mu)=h^*(\be_\mu)=\mu$.
\emc
\eit
Then $\FLX=\la\GLX\cup\Om\ra$.
\end{lemma}

\pf
Since $\FLX=\FX\GLX$, by Remark \ref{rem:factorisation}, it suffices to show that $\FX\sub\la\GLX\cup\Om\ra$.  Let $\mu$ be a cardinal such that $\mu=1$ or $\aleph_0\leq\mu\leq|X|$.  We claim that there exist elements $\ga_\mu,\de_\mu\in\la\GLX\cup\Om\ra$ such that 
\[
\ga_\mu,\de_\mu\in\FX \COMMA t(\ga_\mu)=t(\de_\mu)=|X| \COMMA s(\ga_\mu)<\mu=h(\ga_\mu) \COMMA h(\de_\mu)<\mu=s(\de_\mu),
\]
and either
\bit
\itemnit{a} $s^*(\ga_\mu)<\mu=h^*(\ga_\mu)$ and $h^*(\de_\mu)<\mu=s^*(\de_\mu)$, or
\itemnit{b} $h^*(\ga_\mu)<\mu=s^*(\ga_\mu)$ and $s^*(\de_\mu)<\mu=h^*(\de_\mu)$.
\eit
We prove the claim only in the case in which $h^*(\al_\mu)=s^*(\be_\mu)=\mu$ (i.e., the first option in assumption~(iii)) holds, with the other case being virtually identical.
Suppose the set of lower singletons of $\al_\mu$ is~$V'$, where~${V\sub X}$.  Since $\codom(\al_\mu)\sub X\sm V$ and $t(\al_\mu)=|X|$, we have $|X\sm V|=|X|$.  Let $\ve\in\GLX$ be any element such that $\dom(\ve)=X\sm V$.  Then $\ga_\mu=\al_\mu\ve$ satisfies
\[
s(\ga_\mu)=s(\al_\mu)<\mu=h(\al_\mu)=h(\ga_\mu) \COMMA s^*(\ga_\mu)=0 \COMMA h^*(\ga_\mu)=\mu \COMMA t(\ga_\mu)=|X|.
\]
Similarly, there exists $\eta\in\GLX$ such that $\de_\mu=\be_\mu\eta$ satisfies
\[
h(\de_\mu)<\mu=s(\de_\mu) \COMMA h^*(\de_\mu)=0 \COMMA s^*(\de_\mu)=\mu \COMMA t(\de_\mu)=|X|.
\]
This completes the proof of the claim.
We now note that the set $\Ga=\set{\ga_\mu,\de_\mu}{\mu=1\text{ or }\aleph_0\leq\mu\leq|X|}$ is of the form described in Lemma \ref{lem:U1}.  Lemma \ref{lem:U2} then gives $\FX=\la\S_X\cup\Ga\ra\sub\la\GLX\cup\Om\ra$, as required.
\epf

\begin{rem}\label{rem:FLXSX}
If $\Om$ is of the form described in Lemma \ref{lem:U1}, then $\Om$ is also of the form described in Lemma~\ref{lem:FLXSX}, but the converse is not necessarily true.
\end{rem}

\begin{lemma}\label{lem:FLXSX2}
Let $X$ be an infinite set, and suppose $\Om\sub\FLX$ is such that $\FX\sub\la\S_X\cup\Om\ra$.  For any cardinal~$\mu$ such that $\mu=1$ or $\aleph_0\leq\mu\leq|X|$,
\bit
\itemit{i} there exists $\al_\mu\in \Om$ with $\al_\mu\in\FX$, $h^*(\al_\mu)=\mu$, $\defect(\al_\mu)=2\mu$ and $t(\al_\mu)=|X|$,
\itemit{ii} there exists $\be_\mu\in \Om$ with $\be_\mu\in\FX$, $s^*(\be_\mu)=\mu$, $\defect(\be_\mu)=\mu$ and $t(\be_\mu)=|X|$.
\eit
\end{lemma}

\pf
The proofs being almost identical, we just prove (i).  Let $\si\in\FX$ be such that $s^*(\si)<\mu=h^*(\si)$, noting that $\codef(\si)=2\mu$, and consider an expression $\si=\ga_1\cdots\ga_k$, where the factors all belong to $\S_X\cup\Om$.  
By Corollary \ref{cor:tq}(ii), $h^*(\ga_i)\geq\mu$ for some $i$.  Combined with Lemma \ref{lem:defFLX}, we obtain
\[
2\mu\leq2h^*(\ga_i)\leq\codef(\ga_i)\leq\defect(\ga_i)\leq\defect(\ga_1\cdots\ga_k)= \defect(\si)=2\mu,
\]
so we have equality throught.  In particular, it follows that $h^*(\ga_i)=\mu$, and that $\defect(\ga_i)=\codef(\ga_i)=2\mu$; the latter also gives $\ga_i\in\FX$.  
Corollary \ref{cor:tq}(i) gives $t(\ga_i)=|X|$.  We put $\al_\mu=\ga_i$ (note that $\ga_i\not\in\S_X$, because~$\defect(\ga_i)\not=0$).
\epf

\begin{rem}\label{rem:FLXSX2}
Note that there could be some overlap between the elements from (i) and (ii) in Lemma~\ref{lem:FLXSX2}: namely, if $\mu\geq\aleph_0$, then it is possible to have $\al_\mu=\be_\mu$.  However, if $\mu\not=\nu$, then $\al_\mu\not=\al_\nu$ and $\be_\mu\not=\be_\nu$.
\end{rem}

\begin{lemma}\label{lem:FLXSX3}
Let $X$ be an infinite set, let $\rho$ be the number of cardinals $\mu$ satisfying $\aleph_0\leq\mu\leq|X|$, and suppose $\rho<\aleph_0$.  If $\Om\sub\FLX$ is such that $|\Om\cap\FX|\leq2+2\rho$ and $\FX\sub\la\S_X\cup\Om\ra$, then $\Om\cap\FX$ has the form described in Lemma \ref{lem:FLXSX}, in which case $|\Om\cap\FX|=2+2\rho$.
\end{lemma}

\pf
By Lemma \ref{lem:F3}, $\Om$ contains a subset $\Ga_1=\set{\al_\mu,\be_\mu}{\mu=1\text{ or }\aleph_0\leq\mu\leq|X|}$, where for each $\mu$,
\bit
\itemnit{a} $\al_\mu\in\FX$, $s(\al_\mu)<\mu=h(\al_\mu)$ and $t(\al_\mu)=|X|$,
\itemnit{b} $\be_\mu\in\FX$, $h(\be_\mu)<\mu=s(\be_\mu)$ and $t(\be_\mu)=|X|$. 
\eit
Since $\Ga_1\sub\FX$, it follows that $\Ga_1\sub\Om\cap\FX$.  Since $2+2\rho=|\Ga_1|\leq|\Om\cap\FX|\leq2+2\rho$, it follows that $|\Om\cap\FX|=2+2\rho=|\Ga_1|$.  Since $\rho<\aleph_0$, it also follows that $\Om\cap\FX=\Ga_1$.  
It remains to show that for each~$\mu$, either
\begin{equation}\label{eq:alorbe}
h^*(\al_\mu)=s^*(\be_\mu)=\mu \OR s^*(\al_\mu)=h^*(\be_\mu)=\mu.
\end{equation}
By Lemma \ref{lem:FLXSX2}, $\Om$ contains a subset $\Ga_2=\set{\ga_\mu,\de_\mu}{\mu=1\text{ or }\aleph_0\leq\mu\leq|X|}$, where for each $\mu$,
\bit
\itemnit{c} $\ga_\mu\in\FX$, $h^*(\ga_\mu)=\mu$, $\defect(\ga_\mu)=2\mu$ and $t(\ga_\mu)=|X|$,
\itemnit{d} $\de_\mu\in\FX$, $s^*(\de_\mu)=\mu$, $\defect(\de_\mu)=\mu$ and $t(\de_\mu)=|X|$.
\eit
Again, $\Ga_2\sub\Om\cap\FX=\Ga_1$ (but, as in Remark \ref{rem:FLXSX2}, we might have $|\Ga_2|<|\Ga_1|$).  Now, $\al_1,\be_1,\ga_1,\de_1$ are the only elements of finite (co)defect in (a)--(d).  Since
$\defect(\al_1)=2$, $\defect(\be_1)=1$, $\defect(\ga_1)=2$ and $\defect(\de_1)=1$, it follows that $\al_1=\ga_1$ and $\be_1=\de_1$.  It follows that \eqref{eq:alorbe} holds when $\mu=1$.
Next, suppose $\mu\geq\aleph_0$.  Then $\defect(\al_\mu)=\defect(\be_\mu)=\defect(\ga_\mu)=\defect(\de_\mu)=\mu$, and so $\ga_\mu,\de_\mu\in\{\al_\mu,\be_\mu\}$.  Thus, one of the following must hold:
\bit\bmc2
\itemnit{e} $h^*(\al_\mu)=s^*(\be_\mu)=\mu$, or
\itemnit{f} $h^*(\be_\mu)=s^*(\al_\mu)=\mu$, or
\itemnit{g} $h^*(\al_\mu)=s^*(\al_\mu)=\mu$, or
\itemnit{h} $h^*(\be_\mu)=s^*(\be_\mu)=\mu$.
\emc\eit
If (e) or (f) holds, then \eqref{eq:alorbe} holds.  If (g) holds, then because $\aleph_0\leq\mu=\codef(\be_\mu)=s^*(\be_\mu)+2h^*(\be_\mu)$, we must have $s^*(\be_\mu)=\mu$ or $h^*(\be_\mu)=\mu$, so that \eqref{eq:alorbe} still holds.  Case (h) is treated similarly.
\epf


\begin{thm}\label{thm:FLXGLX}
Let $X$ be an infinite set, and let $\rho$ be the number of cardinals $\mu$ satisfying $\aleph_0\leq\mu\leq|X|$.
\bit
\itemit{i} We have $\relrank{\FLX}{\GLX}=2+2\rho$.  
\itemit{ii} If $\rho<\aleph_0$, and if $\Om\sub \FLX$ with $|\Om|=2+2\rho$, then $\FLX=\la \GLX\cup\Om\ra$ if and only if $\Om$ is of the form described in Lemma \ref{lem:FLXSX}. 
\eit
\end{thm}

\pf
Lemma \ref{lem:FLXSX} gives $\relrank{\FLX}{\GLX}\leq2+2\rho$ and the backwards implication in (ii).  

Next, suppose $\Om\sub\FLX$ is such that $|\Om|=\relrank{\FLX}{\GLX}$ and $\FLX=\la\GLX\cup\Om\ra=\la\S_X\cup(\GLX\sm\S_X)\cup\Om\ra$.  For simplicity, we will write $\Ga=(\GLX\sm\S_X)\cup\Om$, so that $\FLX=\la\S_X\cup\Ga\ra$.  By Lemma \ref{lem:FLXSX2} (cf.~Remark~\ref{rem:FLXSX2}), $\Ga$ contains a subset $\Ga_1$ such that $\Ga_1\sub\FX$ and $|\Ga_1|\geq1+\rho$.  By parts (ii) and (iii) of Lemma \ref{lem:PBunits}, ${\GLX\sm\S_X = \set{\al\in\PB_X}{\codef(\al)=0\not=\defect(\al)}}$.  It follows that $(\GLX\sm\S_X)\cap\FX=\emptyset$, and so ${\Ga_1\sub\Ga\cap\FX=\Om\cap\FX}$.  Hence, $\relrank{\FLX}{\GLX}=|\Om|\geq|\Om\cap\FX|\geq|\Ga_1|\geq1+\rho$.  If $\rho\geq\aleph_0$, then $1+\rho=\rho=2+2\rho$, and so the proof of (i) is complete in this case.  For the remainder of the proof, we assume that $\rho<\aleph_0$.

Now, $\Ga\sub\FLX$ is such that $|\Ga\cap\FX|=|\Om\cap\FX|\leq|\Om|=\relrank{\FLX}{\GLX}\leq2+2\rho$ and $\FX\sub\la\S_X\cup\Ga\ra$.  Lemma~\ref{lem:FLXSX3} then says that $\Om\cap\FX=\Ga\cap\FX$ has the form described in Lemma \ref{lem:FLXSX}, and has size~$2+2\rho$.  But then $2+2\rho\geq\relrank{\FLX}{\GLX}=|\Om|\geq|\Om\cap\FX|=2+2\rho$.  It follows that $\relrank{\FLX}{\GLX}=2+2\rho$, completing the proof of (i).  We also have $|\Om|=|\Om\cap\FX|$; thus, since $\Om$ is finite, $\Om=\Om\cap\FX$ has the specified form.
\epf

The last task of this section is to calculate $\relrank{\FLX}{\S_X}$.  We have already done a lot of the preliminary work for this, but we require one more lemma.

\begin{lemma}\label{lem:FLXSX4}
Let $X$ be an infinite set, and let $\Om=\Om_1\cup\Om_2$, where
\[
\Om_1=\set{\al_\mu,\be_\mu}{\mu=1\text{ or }\aleph_0\leq\mu\leq|X|}\sub\FX \AND {\Om_2=\set{\ga_\mu}{\mu=1\text{ or }\aleph_0\leq\mu\leq|X|}\sub\FLX}
\]
satisfy
\bit\bmc2
\itemit{i} $s(\al_\mu)<\mu=h(\al_\mu)$ and $t(\al_\mu)=|X|$,
\itemit{ii} $h(\be_\mu)<\mu=s(\be_\mu)$ and $t(\be_\mu)=|X|$,
\itemit{iii} $h^*(\al_\mu)=s^*(\be_\mu)=\mu$ or $s^*(\al_\mu)=h^*(\be_\mu)=\mu$,
\itemit{iv} $\codef(\ga_\mu)<\mu=\defect(\ga_\mu)$. 
\emc
\eit
Then $\FLX=\la\S_X\cup\Om\ra$.
\end{lemma}

\pf
We first claim that $\la\S_X\cup\Om\ra$ contains a subset $\Ga=\set{\si_\mu,\tau_\mu}{\mu=1\text{ or }\aleph_0\leq\mu\leq|X|}$ such that for all~$\mu$,
\[
t(\si_\mu)=t(\tau_\mu)=|X| \COMMa
s(\si_\mu)<\mu=h(\si_\mu) \COMMa
s^*(\si_\mu)<\mu=h^*(\si_\mu) \COMMa
h(\tau_\mu)<\mu=s(\tau_\mu) \COMMa
h^*(\tau_\mu)<\mu=s^*(\tau_\mu).
\]
Before we prove the claim, we note that the lemma will then follow.  Indeed, $\Ga$ has the form described in Lemma~\ref{lem:U1}, and $\Om_2$ has the form described in Lemma \ref{lem:Om15}, so Lemmas \ref{lem:U2} and \ref{lem:Om2}, respectively, give $\FX=\la\S_X\cup\Ga\ra$ and $\FLX=\la\FX\cup\Om_2\ra$.  But then
\[
\FLX \supseteq \la\S_X\cup\Om\ra = \la\S_X\cup\Ga\cup\Om\ra \supseteq \la\S_X\cup\Ga\cup\Om_2\ra = \la\la\S_X\cup\Ga\ra\cup\Om_2\ra = \la\FX\cup\Om_2\ra = \FLX.
\]
To establish the claim, we consider two cases.

\pfcase{1}  Suppose first that $\mu=1$.  From assumption (ii), we have $s(\be_1)=1$ and $h(\be_1)=0$, and so $\defect(\be_1)=1$.  Since $\be_1\in\FX$, we have $\codef(\be_1)=1$ as well.  Consequently, we must have $h^*(\be_1)=0$ and $s^*(\be_1)=1$.  Because of $h^*(\be_1)=0$, assumption (iii) gives $h^*(\al_1)=1$.  From assumption (i), we have $h(\al_1)=1$ and $s(\al_1)=0$, so that $\codef(\al_1)=\defect(\al_1)=2$; together with $h^*(\al_1)=1$, it follows that $s^*(\al_1)=0$.  Thus, we may take $\si_1=\al_1$ and $\tau_1=\be_1$.

Thus, we have established the claim in the case $\mu=1$.  Before we consider the infinite case, we note that the argument in Case 1 of the proof of Lemma~\ref{lem:U1} (see in particular the final, parenthesised, sentence) shows that for any $1\leq\nu<\aleph_0$, there exist $\si_\nu,\tau_\nu\in\la\S_X\cup\{\al_1,\be_1\}\ra$ such that
\[
h(\si_\nu)=h^*(\si_\nu) = s(\tau_\nu)=s^*(\tau_\nu) = \nu
\AND
s(\si_\nu)=s^*(\si_\nu) = h(\tau_\nu)=h^*(\tau_\nu) = 0.
\]
Since $\nu<\aleph_0$, these trivially satisfy $t(\si_\nu)=t(\tau_\nu)=|X|$.

\pfcase{2}  We prove the claim for infinite $\mu$ by transfinite induction.  Suppose $\aleph_0\leq\mu\leq|X|$ is such that appropriate elements $\si_\ka,\tau_\ka$ exist for all cardinals $\ka<\mu$.  We just prove the existence of $\si_\mu$, as the existence of $\tau_\mu$ is similar.  Let~$\de\in\{\al_\mu,\be_\mu\}$ be such that $h^*(\de)=\mu$.  Since $t(\al_\mu)=t(\de)=|X|$, and since ${\codef(\al_\mu)=\defect(\de)=\mu}$, there is a permutation $\pi\in\S_X$ that maps $\codom(\al_\mu)$ bijectively onto $\dom(\de)$.  Let $\ve=\al_\mu\pi\de$, so that 
\[
\ve\in\la\S_X\cup\Om\ra \COMMA t(\ve)=|X| \COMMA s(\ve)<\mu=h(\ve) \COMMA h^*(\ve)=\mu.
\]
If we also had $s^*(\ve)<\mu$, then we could take $\si_\mu=\ve$.  So suppose instead that $s^*(\ve)=\mu$.  Let $V'$ be the set of lower singletons of $\ve$, where $V\sub X$, noting that $|V|=\mu$.  We observed above that $\Om_2$ has the form described in Lemma \ref{lem:Om15}, so by that lemma, there exists $\eta\in\la\S_X\cup\Om_2\ra\sub\la\S_X\cup\Om\ra$ such that $\defect(\eta)=\mu$ and $\codef(\eta)=0$.  Let $\xi\in\S_X$ be any permutation that maps $V$ bijectively onto $\Defect(\eta)$.  Then $\si_\mu=\ve\xi\eta\in\la\S_X\cup\Om\ra$ has the desired properties.
\epf

\begin{rem}\label{rem:FLXSX4}
The sets $\Om_1$ and $\Om_2$ in the statement of Lemma \ref{lem:FLXSX4} have the forms described in Lemmas~\ref{lem:FLXSX} and \ref{lem:Om15}, respectively.
\end{rem}

Here is the final main result of this section.

\begin{thm}\label{thm:FLXSX}
Let $X$ be an infinite set, and let $\rho$ be the number of cardinals $\mu$ satisfying $\aleph_0\leq\mu\leq|X|$.
\bit
\itemit{i} We have $\relrank{\FLX}{\S_X}=3+3\rho$.  
\itemit{ii} If $\rho<\aleph_0$, and if $\Om\sub \FLX$ with $|\Om|=3+3\rho$, then $\FLX=\la \S_X\cup\Om\ra$ if and only if $\Om$ is of the form described in Lemma \ref{lem:FLXSX4}. 
\eit
\end{thm}

\pf
Lemma \ref{lem:FLXSX4} gives $\relrank{\FLX}{\S_X}\leq3+3\rho$ and the backwards implication in (ii).  

Next, suppose $\Om\sub\FLX$ is such that $|\Om|=\relrank{\FLX}{\S_X}$ and $\FLX=\la\S_X\cup\Om\ra$.  
Since $\FX\sub\FLX=\la\S_X\cup\Om\ra$, Lemma~\ref{lem:F3} shows that $\Om$ contains a subset $\Om_1$ such that $\Om_1\sub\FX$ and $|\Om_1|=2+2\rho$.  
Also, since certainly $\FLX=\la\FX\cup\Om\ra$, Lemma \ref{lem:Om1} shows that $\Om$ contains a subset $\Om_2$ of the form described in Lemma \ref{lem:Om15}.  Note that $|\Om_2|=1+\rho$, and that $\Om_2\sub\FLX\sm\FX$.  
In particular, since $\Om_1\cap\Om_2=\emptyset$, it follows that $\relrank{\FLX}{\S_X}=|\Om|\geq|\Om_1|+|\Om_2|=3+3\rho$, completing the proof of (i).  

To complete the proof of (ii), suppose $\rho<\aleph_0$, and let $\Om,\Om_1,\Om_2$ be as in the previous paragraph.  Then by finiteness of all three sets, and since $|\Om|=|\Om_1|+|\Om_2|$, we must have $\Om=\Om_1\cup\Om_2$.  We have already noted that $\Om_2$ has the form described in Lemma \ref{lem:Om15}, and that $\Om_2\sub\FLX\sm\FX$ and $\Om_1\sub\FX$.  Thus, $|\Om\cap\FX|=|\Om_1|=2+2\rho$, and so Lemma \ref{lem:FLXSX3} says that $\Om_1=\Om\cap\FX$ has the form described in Lemma~\ref{lem:FLXSX}.  Thus, $\Om=\Om_1\cup\Om_2$ has the form described in Lemma \ref{lem:FLXSX4} (cf.~Remark \ref{rem:FLXSX4}).
\epf

\section{Sierpi\'nski rank and the semigroup Bergman property}\label{sect:SB}

Recall from \cite{MP2012} that the \emph{Sierpi\'nski rank} of a semigroup $S$, denoted $\SR(S)$, is the least integer $n$ such that every countable subset of $S$ is contained in an $n$-generator subsemigroup of $S$, if such an integer exists; otherwise, we say $S$ has infinite Sierpi\'nski rank and write $\SR(S)=\infty$.  Every finitely generated semigroup trivially has finite Sierpi\'nski rank, and this then coincides with the \emph{rank} of the semigroup, as defined in Section \ref{sect:monoids}.
Recall from \cite{MMR2009} that a semigroup~$S$ has the \emph{semigroup Bergman property} if every generating set for $S$ has a bounded length function.  Finite semigroups trivially have the Bergman property, but this is not true of arbitrary finitely generated semigroups (consider a free semigroup of finite rank).
The main results of this section (Theorems \ref{thm:SBPB} and \ref{thm:SBEF}) use results of previous sections to calculate the Sierpi\'nski rank for each of the monoids $\PB_X,\EX,\GLX,\GRX,\FX,\FLX,\FRX$, and also determine which of them have the semigroup Bergman property.

For the proof of the first lemma, we recall again that $\PB_X$ is a submonoid of the larger \emph{partition monoid}~$\P_X$.  As before, we will not recall the full definition of $\P_X$ here; the reader may refer to \cite{EF2012,East2014}, where the focus was on the infinite case.  
Recall from \cite{MMR2009} that a semigroup $S$ is \emph{strongly distorted} if there exists a sequence $(a_1,a_2,a_3,\ldots)$ of natural numbers, and a natural number $N$ such that, for all sequences $(s_1,s_2,s_3,\ldots)$ of elements from~$S$, there exists a subset $T$ of $S$ with $|T|=N$ such that each $s_n$ can be factorised as a product of length at most~$a_n$ over $T$.  It follows from \cite[Lemma 2.4 and Proposition 2.2(i)]{MMR2009} that a strongly distorted semigroup that is not finitely generated has the semigroup Bergman property.

\begin{lemma}\label{lem:SD}
If $X$ is an infinite set, then $\PB_X$ is strongly distorted.
\end{lemma}

\pf
Let $(\ga_1,\ga_2,\ga_3,\ldots)$ be a sequence of elements of the partition monoid $\P_X$.  It was shown in \cite[Theorem~37]{East2014} that there exist elements $\al,\be\in\P_X$ such that $\ga_n=\al\be\al^n\be^2(\al^*)^n\be^*\al^*$ for each $n$.  We will not repeat the construction here, but it is easy to check that if all the elements $\ga_n$ belong to $\PB_X$, then the constructed elements $\al,\be$ (and hence also $\al^*,\be^*$) belong to $\PB_X$ as well.  It follows that $\PB_X$ is strongly distorted; we take $N=4$ and $a_n=2n+6$ for all $n$.
\epf

It follows from the previous proof that $\SR(\PB_X)\leq4$.  In Theorem \ref{thm:SBPB}, we will show that $\SR(\PB_X)=2$; the proof we give is an adaptation of an ingenious argument of Hyde and P\'eresse \cite{HP2012} originally purposed for the symmetric inverse monoid $\I_X$.  
Recall that a permutation $\al\in\S_X$ is an \emph{involution} if $\al^2=1$ (we consider the identity element to be an involution).  Part~(i) of the following lemma was proved in \cite[Lemma~2.2]{Galvin1995}, and part (ii) in \cite[Lemma 2.4]{HP2012}.

\begin{lemma}\label{lem:SX}
Let $X$ be an arbitrary set, and let $\al\in\S_X$.
\bit
\itemit{i} There exist two involutions $\be,\ga\in\S_X$ such that $\al=\be\ga$.
\itemit{ii} There exists an involution $\de\in\S_X$ such that $\al^{-1}\in\la\al,\al\de\ra$. \epfres
\eit
\end{lemma}

Recall that if $\al\in\PB_X$, and if $Y\sub\dom(\al)$ and $Z\sub\codom(\al)$, then $Y\al=\set{y\al}{y\in Y}$ and $Z\al^{-1}=\set{z\al^{-1}}{z\in Z}$.

\begin{thm}\label{thm:SBPB}
If $X$ is an infinite set, then $\PB_X$ has the semigroup Bergman property, and $\SR(\PB_X)=2$.
\end{thm}

\pf
Since $\PB_X$ is uncountable, it is not finitely generated.  The Bergman property then follows immediately from Lemma \ref{lem:SD} and the above-mentioned results from~\cite{MMR2009}.  We noted above that $\SR(\PB_X)\leq4$.  Since also $\SR(\PB_X)\geq2$ (as $\PB_X$ is not commutative), it suffices to show that any three elements of $\PB_X$ belong to a subsemigroup generated by two elements.  With this in mind, let $\ga_1,\ga_2,\ga_3\in\PB_X$ be arbitrary.  We will construct elements $\al,\be\in\PB_X$ such that $\ga_1,\ga_2,\ga_3\in\la\al,\be\ra$.

Since $X$ is infinite, we may fix a decomposition $X=\bigsqcup_{i=0}^\infty X_i=X_0\sqcup X_1\sqcup X_2\sqcup\cdots$, where $|X_i|=|X|$ for each $i$.  We let $\al\in\PB_X$ be any element with $h^*(\al)=|X|$ and $X_i\al=X_{i+1}$ for all $i\geq0$.  Note that $\dom(\al)=X$ and that $\codom(\al)=X\sm X_0=\bigcup_{i=1}^\infty X_i$.  The definition of $\be$ is far more involved, and is achieved in a number of stages.  First, for each $i\in\{11,\ldots,18\}$, let $\si_i\in\S_{X_i}$ be an involution of $X_i$; the exact definition of the $\si_i$ will be given later.  We then let $\be$ be any element of $\PB_X$ such that
\bit
\bmc2
\item $h(\be)=|X|$,
\item $\dom(\be)=\bigcup_{i=11}^\infty X_i$,
\emc
\item $X_i\be=X_i$ for all $i\in\{11,\ldots,18\}$, and the restriction of $\be$ to $X_i$ is $\si_i$ for all such $i$,
\bmc4
\item $X_{19}\be=X_{19}\cup X_{20}$,
\item $X_{20}\be=\bigcup_{i=21}^\infty X_i$,
\item $X_{21}\be=\bigcup_{i=1}^{10}X_i$,
\item $\big(\bigcup_{i=22}^\infty X_i\big)\be=X_0$.
\emc
\eit
It is easy to check that $\codom(\be)=X$.  Also, since the $\si_i$ are involutions, $\be^2$ does not depend on the choices of~$\si_i$.  This means that we may use $\be^2$ to define the involutions $\si_i$, with no fear of circularity.  Note also that $\dom(\be^2)=\bigcup_{i=11}^{20}X_i$, $\codom(\be^2)=X$ and $h(\be^2)=|X|$, with the last of these following from Lemma \ref{lem:shdef}(ii).  It follows from Lemma \ref{lem:alSXbe} that $\PB_X=\al\S_X\be^2$.  Thus, there exist $\de_1,\de_2,\de_3\in\S_X$ such that $\ga_i=\al\de_i\be^2$, for $i=1,2,3$.  By Lemma \ref{lem:SX}(i), there exist involutions $\ve_1,\ldots,\ve_6\in\S_X$ such that $\de_1,\de_2,\de_3\in\la\ve_1,\ldots,\ve_6\ra$.  Note that the $\de_i$ and $\ve_i$ do not depend on the involutions $\si_i$, since their definitions involve only $\al$ and $\be^2$.  (This is why we used $\PB_X=\al\S_X\be^2$ instead of $\PB_X=\al\S_X\be$, which is also true, in order to define the $\de_i$ and $\ve_i$.)  The proof will be complete if we can show that $\ve_1,\ldots,\ve_6\in\la\al,\be\ra$.

Now we define $\pi=\al^{22}\be$ and $\tau=\al^{21}\be\al^{10}\be^2$.  It is routine to check that
\begin{equation}\label{eq:pitau}
\dom(\pi)=\codom(\tau)=X \AND \codom(\pi)=\dom(\tau)=X_0,
\end{equation}
and that $\pi,\tau$ do not depend on the choices of the involutions $\si_i$.  (Note also that $h^*(\pi)=h(\tau)=|X|$, but this will not concern us.)  In light of \eqref{eq:pitau}, we have $\tau^*\tau=1$, and $\pi\tau\in\S_X$.  By Lemma \ref{lem:SX}(ii), there is an involution $\ve_7\in\S_X$ such that $(\pi\tau)^*=(\pi\tau)^{-1}\in\la\pi\tau,(\pi\tau)\ve_7\ra$.  For reasons that will become clear later, we also let $\ve_8=1$ be the identity of $\S_X$.  We now use the involutions $\ve_1,\ldots,\ve_8\in\S_X$ to define the involutions~${\si_i\in\S_{X_i}}$, for $i=11,\ldots,18$.  First, it is easy to check that
\[
\dom(\tau^*\al^n)=X \AND \codom(\tau^*\al^n)=X_n \qquad\text{for any $n\in\N$,}
\]
and that $\tau^*\al^n$ does not depend on the choices of the $\si_i$.
It follows that for $i\in\{1,\ldots,8\}$, the domain and codomain of $(\tau^*\al^{10+i})^* \ve_i(\tau^*\al^{10+i})$ are both equal to $X_{10+i}$; for any such $i$, we let $\si_{10+i}$ be the restriction of $(\tau^*\al^{10+i})^* \ve_i(\tau^*\al^{10+i})$ to $X_{10+i}$.  
So $\si_i\in\S_{X_i}$ for each such $i$, and each $\si_i$ is an involution because the $\ve_i\in\S_X$ are involutions.  We have now completed the definition of $\be$.  

For $i\in\{1,\ldots,8\}$, define $\eta_i=\al^{22}\be\al^{10+i}\be\al^{11-i}\be\al^{10}\be^2$.  One may check that $\dom(\eta_i)=\codom(\eta_i)=X$, so that $\eta_i\in\S_X$.  Furthermore, if $i\in\{1,\ldots,8\}$, then for any $x\in X$, we have $x\al^{22}\be\al^{10+i}=x\pi\al^{10+i}\in X_{10+i}$, so that $x\al^{22}\be\al^{10+i}\be=(x\pi\al^{10+i})\si_{10+i}$.  Using this, and the fact that $\al\al^*=\tau^*\tau=1$, we then calculate
\begin{align*}
x\eta_i = (x \al^{22}\be\al^{10+i}\be)\al^{11-i}\be\al^{10}\be^2 &= (x\pi\al^{10+i})\si_{10+i}\al^{11-i}\be\al^{10}\be^2 \\
&= x\pi\al^{10+i}[(\tau^*\al^{10+i})^* \ve_i(\tau^*\al^{10+i})]\al^{11-i}\be\al^{10}\be^2 \\
&= x\pi\al^{10+i}(\al^{10+i})^*\tau \ve_i\tau^*\al^{21}\be\al^{10}\be^2 \\
&= x\pi[\al^{10+i}(\al^{10+i})^*]\tau \ve_i[\tau^*\tau] 
=  x\pi\tau \ve_i.
\end{align*}
Thus, $\eta_i$ and $(\pi\tau)\ve_i$ contain the same transversals.  Since $\eta_i$ and $(\pi\tau)\ve_i$ both belong to $\S_X$, it follows that they are equal.  In particular, $(\pi\tau)\ve_i\in\la\al,\be\ra$ for all $i\in\{1,\ldots,8\}$.  Taking $i=8$, and recalling that $\ve_8=1$, we obtain $\pi\tau=(\pi\tau)\ve_8\in\la\al,\be\ra$.  Then also $(\pi\tau)^{-1}\in\la\pi\tau,(\pi\tau)\ve_7\ra\sub\la\al,\be\ra$.  It follows that for any $i\in\{1,\ldots,8\}$, $\ve_i = (\pi\tau)^{-1}(\pi\tau)\ve_i\in\la\al,\be\ra$.  As noted above, this completes the proof.
\epf

\begin{rem}
In the definition of $\al$ and $\be$ in the above proof, we specified that $h^*(\al)=h(\be)=|X|$, but said nothing about singletons.  Thus, as with Lemma \ref{lem:alSXbe} (cf.~Corollary \ref{cor:PB=BB}), $\al,\be$ could be chosen to have \emph{no} singletons: i.e., to belong to $\B_X$, the set of all \emph{full} Brauer graphs.  This means that any countable subset of~$\PB_X$ belongs to a subsemigroup of $\PB_X$ generated by two elements of $\B_X$.
\end{rem}

We now move on to consider the monoids $\EX,\GLX,\GRX,\FX,\FLX,\FRX$.  We require the following two results; the first is \cite[Proposition 5 and Remark 7]{East2012}, and the second follows from \cite[Theorem~3.5]{Galvin1995}.

\begin{lemma}\label{lem:SRM}
Let $M$ be a monoid, write $G=\G(M)$, and suppose $M\sm G$ is an ideal of $M$.  If $\SR(G)$ and $\relrank MG$ are both finite, then $\SR(M)=\SR(G)+\relrank MG$.  \epfres
\end{lemma}

\begin{thm}\label{thm:Galvin}
If $X$ is an infinite set, then the symmetric group $\S_X$ has Sierpi\'nski rank $2$. \epfres
\end{thm}

We are now ready to prove the second main result of this section.

\begin{thm}\label{thm:SBEF}
Let $X$ be an infinite set.  Then
\bit
\itemit{i} $\SR(\EX)=\infty$,\\[-6mm]
\itemit{ii} $\SR(\GLX)=\SR(\GRX)=\SR(\FX)=\begin{cases}
2n+6 &\text{if $|X|=\aleph_n$, where $n\in\N$}\\
\infty &\text{otherwise,}
\end{cases}$
\itemit{iii} $\SR(\FLX)=\SR(\FRX)=\begin{cases}
3n+8 &\text{if $|X|=\aleph_n$, where $n\in\N$}\\
\infty &\text{otherwise,}
\end{cases}$
\itemit{iv} None of $\EX,\GLX,\GRX,\FX,\FLX,\FRX$ have the semigroup Bergman property.
\eit
\end{thm}

\pf
(i).  For $y\in X$, define $\ve_y=\partnVI xx_{x\in X\sm\{y\}}$.  Let $Y\sub X$ be a countably infinite subset of $X$, and put $\Ga=\set{\ve_y}{y\in Y}$.  It suffices to show that $\Ga$ is not contained in a finitely generated subsemigroup of $\EX$.  To do so, suppose $\Ga\sub\la\Om\ra$, where $\Om\sub\EX$.  Fix some $y\in Y$, and consider an expression $\ve_y=\al_1\cdots\al_k$, where $\al_1,\ldots,\al_k\in\Om$.  Without loss of generality, we may assume that $\al_1\not=1$.  By Lemma \ref{lem:shdef}(iv), $\defect(\al_1)\leq\defect(\al_1\cdots\al_k)=\defect(\ve_y)=1$, so Theorem \ref{thm:IGPBX} gives $\sh(\al_1)=0$.  If $\defect(\al_1)=0$, then we would have $\al_1=1$, which we have excluded, so we must have $\defect(\al_1)=1$; together with $\sh(\al_1)=0$, it follows that $\al_1=\ve_z$ for some $z\in X$.  But then $z$ is an upper singleton of $\al_1$, and hence also of $\al_1\cdots\al_k=\ve_y$, so it follows that $z=y$, giving $\ve_y=\al_1\in\Om$.  We have shown that $\Ga\sub\Om$, and so $|\Om|\geq\aleph_0$, as required.

\pfitem{ii) and (iii}  Let $\Q$ denote any of $\GLX$, $\FLX$ or $\FX$; the cases in which $\Q$ is $\GRX$ or $\FRX$ are dual.  

Suppose first that $|X|=\aleph_n$, where $n\in\N$.  Lemmas \ref{lem:equiv}, \ref{lem:FFLFR} and \ref{lem:GGLGR}, show that $\G(\Q)=\S_X$ and that~${\Q\sm\S_X}$ is an ideal of $\Q$.  Theorems~\ref{thm:GLXSX},~\ref{thm:FXSX} and~\ref{thm:FLXSX} give
\[
\relrank{\GLX}{\S_X}=\relrank{\FX}{\S_X}=2n+4 \AND \relrank{\FLX}{\S_X}=3n+6.
\]
The stated formulae for $\SR(\Q)$ now follow from Lemma~\ref{lem:SRM} and Theorem~\ref{thm:Galvin}.

Suppose now that $|X|>\aleph_n$ for all $n\in\N$.  For each $n$, let $\al_n\in\Q$ be such that $\defect(\al_n)=\aleph_n$, and suppose $\Om\sub\Q$ is such that $\set{\al_n}{n\in\N}\sub\la \Om\ra$.  Let $n\in\N$ be arbitrary, and consider an expression $\al_n=\be_1\cdots\be_k$, where $\be_1,\ldots,\be_k\in \Om$.  
Corollary \ref{cor:tq}(ii) gives $\defect(\be_i)\geq\aleph_n$ for some $i$.  But then, since~$\Om\sub\FLX$, Lemma~\ref{lem:defFLX} gives $\defect(\be_i)\leq\defect(\be_1\cdots\be_k)=\defect(\al_n)=\aleph_n$, so that $\defect(\be_i)=\aleph_n$.  Thus, $\Om$ contains an element of defect~$\aleph_n$ for each $n\in\N$, and it follows that $|\Om|\geq\aleph_0$.  Thus, $\set{\al_n}{n\in\N}$ is not contained in any finitely generated subsemigroup of $\FX$, and so $\SR(\FX)=\infty$.

\pfitem{iv}  Let $\Q$ denote any of $\EX$, $\GLX$, $\FLX$ or $\FX$; the cases in which $\Q$ is $\GRX$ or $\FRX$ are dual.  
We claim that there exists a generating set $\Ga$ of $\Q$ such that every element of $\Q$ of finite defect has defect at most~$2$.  

Before we prove the claim, we that show the length function with respect to any such generating set $\Ga$ is unbounded.  To do so, let $n\in\N$ be arbitrary.  We must show that there exists $\al\in\Q$ such that any factorisation of $\al$ over $\Ga$ involves at least $n$ factors.
To do so, let $\al\in\Q$ be such that $\defect(\al)=2n$.  Consider an expression $\al=\be_1\cdots\be_k$, where $\be_1,\ldots,\be_k\in\Ga$.  If $\defect(\be_i)\geq\aleph_0$ for some $i$, then Lemma \ref{lem:defFLX} would give $2n=\defect(\be_1\cdots\be_k)\geq\defect(\be_i)\geq\aleph_0$, a contradiction.  Thus, each $\be_i$ has finite defect, and so, by assumption, we must have $\defect(\be_i)\leq2$ for each $i$.  Together with Lemma \ref{lem:shdef}(iv), this gives
\[
2n = \defect(\be_1\cdots\be_k) \leq \defect(\be_1)+\cdots+\defect(\be_k) \leq 2k,
\]
so that $k\geq n$.  That is, any factorisation of $\al$ over $\S_X\cup \Ga$ must involve at least $n$ factors.

It remains only to prove the above claim.  If $\Q$ is one of $\GLX$,~$\FX$ or $\FLX$, then Lemma \ref{lem:GLXSX1}, \ref{lem:U2} or \ref{lem:FLXSX4}, respectively, gives a subset $\Om$ of $\Q$ such that $\Ga=\S_X\cup\Om$ has the desired form.  It remains to prove the claim for $\EX$.  For $x\in X$, let $\ve_x$ be as in part (i).  For distinct $x,y\in X$, define $\eta_{xy} = \partn z{x,y}z{x,y}_{z\in X\sm\{x,y\}}$.  Let $\Om_1,\Om_2,\Om_3$ be as in the proof of Theorem~\ref{thm:IGPBX}, and put
\[
\Ga = \{1\} \cup \set{\ve_x}{x\in X} \cup \set{\eta_{xy}}{x,y\in X,\ x\not=y} \cup \Om_3.
\]
It follows quickly from the proof of \cite[Theorem 3.18]{DEG2017} that $\la \Ga\sm\Om_3\ra=\Om_1\cup\Om_2$.  In particular, Theorem~\ref{thm:IGPBX} gives $\EX=\Om_1\cup\Om_2\cup\Om_3\sub\la \Ga\ra\sub\EX$, so $\la \Ga\ra=\EX$.  Clearly $\Ga$ has the desired form.
\epf

\begin{rem}
If $X$ is finite and $|X|=n\geq3$, then $\GLX=\GRX=\S_X$, $\FLX=\FRX=\FX=\PB_X$ and $\EX$ are all finite, and so
\[
\SR(\S_X)=\rank(\S_X)=2 \COMMA \SR(\PB_X)=\rank(\PB_X)=4 \COMMA \SR(\EX)=\rank(\EX)=1+\tbinom{n+1}2=\tfrac{n^2+n+2}2.
\]
Indeed, the first of these is folklore, while the second and third are parts of \cite[Proposition 3.16 and Theorem~3.18]{DEG2017}.  It follows from Theorem \ref{thm:SBEF}(ii) that the Continuum Hypothesis is equivalent to the assertion that~$\mathcal F_{\mathbb R}$ has Sierpi\'nski rank $8$.
\end{rem}

\subsection*{Acknowledgements}

We thank the referee for a number of helpful suggestions.

\footnotesize
\def\bibspacing{-1.1pt}
\bibliography{biblio}
\bibliographystyle{abbrv}
\end{document}